# MRA-Wavelet subspace architecture for logic, probability, and symbolic sequence processing

DANIEL J. GREENHOE

*Abstract*: The linear subspaces of a *multiresolution analysis* (*MRA*) and the linear subspaces of the *wavelet analysis* induced by the MRA, together with the *set inclusion* relation $\subseteq$, form a very special lattice of subspaces which herein is called a *primorial lattice*. This paper introduces an operator $\mathbf{R}$ that extracts a set of $2^{N-1}$ element Boolean lattices from a $2^N$ element Boolean lattice. Used recursively, a sequence of Boolean lattices with decreasing order is generated—a structure that is similar to an *MRA*. A second operator, which is a special case of a "*difference operator*", is introduced that operates on consecutive Boolean lattices $\boldsymbol{L}_2^n$ and $\boldsymbol{L}_2^{n-1}$ to produce a sequence of *orthocomplemented lattices*. These two sequences, together with the subset ordering relation $\subseteq$, form a *primorial lattice* $\mathbb{P}$. A *logic* or *probability* constructed on a Boolean lattice $\boldsymbol{L}_2^N$ likewise induces a primorial lattice $\mathbb{P}$. Such a logic or probability can then be rendered at $N$ different "resolutions" by selecting any one of the $N$ Boolean lattices in $\mathbb{P}$ and at $N$ different "frequencies" by selecting any of the $N$ different orthocomplemented lattices in $\mathbb{P}$. Furthermore, $\mathbb{P}$ can be used for *symbolic sequence analysis* by projecting sequences of symbols onto the sublattices in $\mathbb{P}$ using one of three lattice projectors introduced. $\mathbb{P}$ can be used for *symbolic sequence processing* by judicious rejection and selection of projected sequences. Examples of symbolic sequences include sequences of logic values, sequences of probabilistic events, and genomic sequences (as used in "*genomic signal processing*").

*2010 Mathematics Subject Classification* 06B05,42C40,03B50 (primary); 03G05,03G10,03B47,03B60 (secondary)

*Keywords*: multiresolution analysis, MRA, wavelets, lattice theory, order theory, subspaces, genomic signal processing, genomic sequence processing, GSP, multi-valued logic, ortho logic

# Contents

# 1 Background: lattices

## 1.1 Order

### 1.1.1 Order relations

**Definition 1.1** [1] Let $X$ be a set.Let $2^{XX}$ be the set of all relations on $X$. A relation $\le$ is an **order relation** in $2^{XX}$ if

| | | | | | |
|---|---|---|---|---|---|
| 1. | $x \le x$ | $\forall x \in X$ | (*reflexive*) | and | preorder |
| 2. | $x \le y$ and $y \le z \implies x \le z$ | $\forall x,y,z \in X$ | (*transitive*) | and | preorder |
| 3. | $x \le y$ and $y \le x \implies x = y$ | $\forall x,y \in X$ | (*anti-symmetric*) | | |

An **ordered set** is the pair $(X, \le)$. The set $X$ is called the **base set** of $(X, \le)$. If $x \le y$ or $y \le x$, then elements $x$ and $y$ are said to be **comparable**, denoted $x \sim y$. Otherwise they are **incomparable**, denoted $x \| y$. The relation $\lneq$ is the relation $\le \setminus =$ ("less than but not equal to"), where $\setminus$ is the *set difference* operator, and $=$ is the equality relation.

**Definition 1.2** [2] Let $(X, \le)$ be an *ordered set* (Definition 1.1 page 3). Let $2^{XX}$ be the set of all relations on $X$. The relations $\ge$, $<$, $> \in 2^{XX}$ are defined as follows:

[1] [113], page 470, [12], page 1, [105], page 156, ⟨I, II, (1)⟩, [38], page 373, ⟨I–III⟩. An *order relation* is also called a **partial order relation**. An *ordered set* is also called a **partially ordered set** or **poset**.

[2] [139], page 2

$$x \geq y \overset{\text{def}}{\iff} y \leq x \qquad \forall x,y \in X$$
$$x \lneq y \overset{\text{def}}{\iff} x \leq y \text{ and } x \neq y \qquad \forall x,y \in X$$
$$x \gneq y \overset{\text{def}}{\iff} x \geq y \text{ and } x \neq y \qquad \forall x,y \in X$$

The relation $\geq$ is called the **dual** of $\leq$.

**Example 1.3**

| order relation | | dual order relation | |
|---|---|---|---|
| $\leq$ | (integer less than or equal to) | $\geq$ | (integer greater than or equal to) |
| $\subseteq$ | (subset) | $\supseteq$ | (super set) |
| $\mid$ | (divides) | | (divided by) |
| $\Longrightarrow$ | (implies) | $\Longleftarrow$ | (implied by) |

**Definition 1.4** [3] A relation $\leq$ is a **linear order relation** on $X$ if

1. $\leq$ is an *order relation* (Definition 1.1 page 3) and
2. $x \leq y$ or $y \leq x$ $\forall x,y \in X$ (*comparable*).

A **linearly ordered set** is the pair $(X, \leq)$.
A linearly ordered set is also called a **totally ordered set**, a **fully ordered set**, and a **chain**.

### 1.1.2 Representation

**Definition 1.5** [4] $y$ **covers** $x$ in the ordered set $(X, \leq)$ if

1. $x \leq y$ ($y$ is greater than $x$) and
2. $(x \leq z \leq y) \implies (z = x \text{ or } z = y)$ (there is no element between $x$ and $y$).

The case in which $y$ covers $x$ is denoted $x \prec y$.

An ordered set can be represented in any of three ways:

- *Hasse diagram* (Definition 1.6 page 4)
- a set of ordered pairs of *order relations* (Definition 1.1 page 3)
- a set of ordered pairs of *cover relations* (Definition 1.5 page 4)

**Definition 1.6** Let $(X, \leq)$ be an ordered pair. A diagram is a **Hasse diagram** of $(X, \leq)$ if it satisfies the following criteria:

- Each element in $X$ is represented by a dot or small circle.
- For each $x, y \in X$, if $x \prec y$, then $y$ appears at a higher position than $x$ and a line connects $x$ and $y$.

---

[3] [113], page 470, [133], page 410
[4] [14], page 445

**Example 1.7** Here are three ways of representing the ordered set $\left(2^{\{x,y\}}, \subseteq\right)$;

(1) **Hasse diagrams**: If two elements are comparable, then the lesser of the two is drawn lower on the page than the other with a line connecting them.

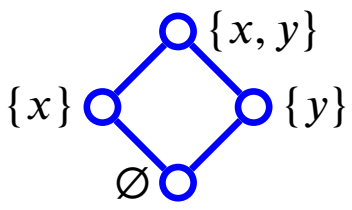

(2) Sets of ordered pairs specifying *order relations* (Definition 1.1 page 3):

$$\subseteq = \left\{ \begin{array}{llll} (\varnothing,\varnothing), & (\{x\},\{x\}), & (\{y\},\{y\}), & (\{x,y\},\{x,y\}), \\ (\varnothing,\{x\}), & (\varnothing,\{y\}), & (\varnothing,\{x,y\}), & (\{x\},\{x,y\}), (\{y\},\{x,y\}) \end{array} \right\}$$

(3) Sets of ordered pairs specifying *covering relations*:

$$< = \left\{ (\varnothing,\{x\}), \quad (\varnothing,\{y\}), \quad (\{x\},\{x,y\}), (\{y\},\{x,y\}) \right\}$$

### 1.1.3 Decomposition

**Definition 1.8** [5] The tupple $(Y, \olessthan)$ is a **subposet** of the ordered set $(X, \le)$ if

1. $Y \subseteq X$ ($Y$ is a subset of $X$) and
2. $\olessthan = (\le \cap Y^2)$ ($\olessthan$ is the relation $\le$ restricted to $Y \times Y$)

**Example 1.9**

Subposets of 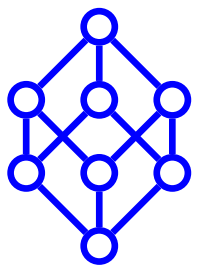 include 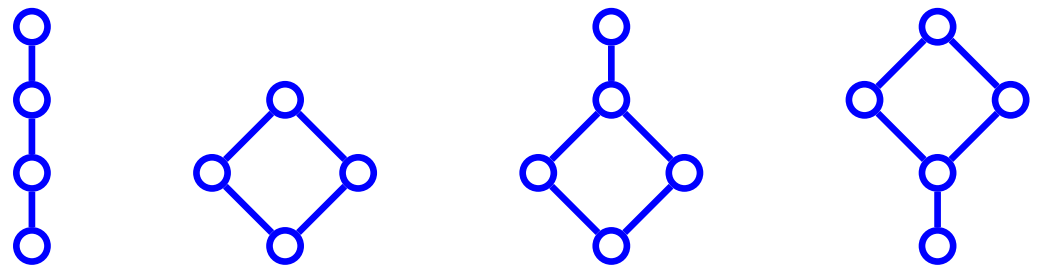

**Example 1.10** Let

$$(X, \le) \triangleq \Big(\{0, a, b, c, p, 1\}, \ \big\{ (0,0), (a,a), (b,b), (c,c), (p,p), (1,1), (0,a), (0,b), (0,c), (0,p), (0,1), (a,b), (a,c), (a,1), (p,1), (b,c), (b,1), (c,1), (p,1) \big\}\Big)$$

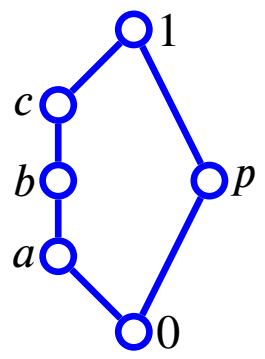

$$(Y, \olessthan) \triangleq \Big(\{0, a, c, p, 1\}, \ \big\{ (0,0), (a,a), (c,c), (p,p), (1,1), (0,a), (0,c), (0,p), (0,1), (a,c), (a,1), (p,1), (c,1), (p,1) \big\}\Big).$$

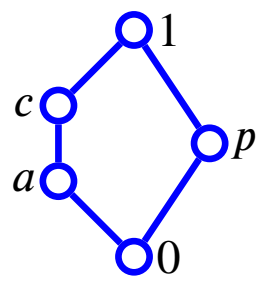

Then $(Y, \olessthan)$ is a subposet of $(X, \le)$ because $Y \subseteq X$ and $\olessthan = (\le \cap Y^2)$.

[5] [72], page 2

A *chain* is an ordered set in which every pair of elements is *comparable* (Definition 1.4 page 4). An *antichain* is just the opposite—it is an ordered set in which *no* pair of elements is comparable (next definition).

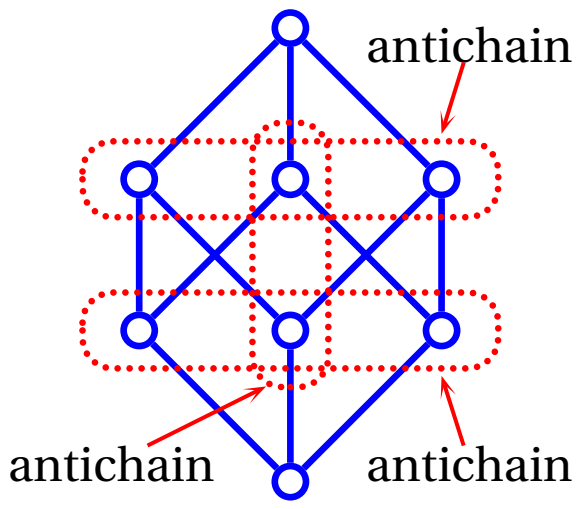

**Definition 1.11** [6] The subposet $(A, \olessthan)$ in the ordered set $(X, \le)$ is an **antichain** if all elements in $A$ are *incomparable* (Definition 1.1 page 3), such that

$$x \| y \qquad \forall x, y \in A$$

**Definition 1.12** [7] The **length** $\ell(\boldsymbol{L})$ of a *chain* (Definition 1.4 page 4) $\boldsymbol{L}$ with $N$ elements is $N-1$. The **length** of an *ordered set* (Definition 1.1 page 3) is the length of the longest chain in the ordered set. The **width** of an ordered set is the number of elements in the largest *antichain* in the ordered set.

**Theorem 1.13** (Dilworth's theorem) [8] *Let $(X, \le)$ be an ordered set.*

$$\left\{ \begin{array}{l} \text{WIDTH } N \text{ of } (X, \le) \\ \text{is FINITE} \end{array} \right\} \implies \left\{ \begin{array}{ll} 1. & \text{there exists a partition of } (X, \le) \text{ into } N \text{ chains} \quad \text{and} \\ 2. & \text{there does not exist any partition} \\ & \text{of } (X, \le) \text{ into less than } N \text{ chains} \end{array} \right\}$$

**Definition 1.14** [9] Let $X$ and $Y$ be disjoint sets. Let $\boldsymbol{P} \triangleq (X, \olessthan)$ and $\boldsymbol{Q} \triangleq (Y, \lessdot)$ be ordered sets on $X$ and $Y$. The **direct sum** of $\boldsymbol{P}$ and $\boldsymbol{Q}$ is defined as

$$\boldsymbol{P} + \boldsymbol{Q} \triangleq (X \cup Y, \le)$$

where $x \le y$ if

1. $x, y \in X$ and $x \olessthan y$ or
2. $x, y \in Y$ and $x \lessdot y$

The direct sum operation is also called the **disjoint union**. The notation $n\boldsymbol{P}$ is defined as

$$n\boldsymbol{P} \triangleq \underbrace{\boldsymbol{P} + \boldsymbol{P} + \cdots + \boldsymbol{P}}_{n-1 \text{ “+” operations}}.$$

**Definition 1.15** [10] Let $X$ and $Y$ be disjoint sets. Let $\boldsymbol{P} \triangleq (X, \olessthan)$ and $\boldsymbol{Q} \triangleq (Y, \lessdot)$ be ordered sets on $X$ and $Y$. The **direct product** of $\boldsymbol{P}$ and $\boldsymbol{Q}$ is defined as

$$\boldsymbol{P} \times \boldsymbol{Q} \triangleq (X \times Y, \le)$$

where $(x_1, y_1) \le (x_2, y_2)$ if $x_1 \olessthan x_2$ and $y_1 \olessthan y_2$.

[6] [72], page 2
[7] [72], page 2, [18], page 5
[8] [47], page 161, [48], [56], page 4
[9] [155], page 100
[10] [155], pages 100–101, [154], page 43

The direct product operation is also called the **cartesian product**. The order relation $\le$ is called a **coordinate wise** order relation. The notation $\mathbf{P}^n$ is defined as

$$\mathbf{P}^n \triangleq \underbrace{\mathbf{P} \times \mathbf{P} \times \cdots \times \mathbf{P}}_{n-1 \text{ ``}\times\text{'' operations}}.$$

**Definition 1.16** [11] Let $X$ and $Y$ be disjoint sets. Let $\mathbf{P} \triangleq (X, \circledcirc)$ and $\mathbf{Q} \triangleq (Y, \lessdot)$ be ordered sets on $X$ and $Y$. The **ordinal sum** of $\mathbf{P}$ and $\mathbf{Q}$ is defined as

$$\mathbf{P} \oplus \mathbf{Q} \triangleq (X \cup Y, \le)$$

where $x \le y$ if

1. $x, y \in X$ and $x \circledcirc y$ or
2. $x, y \in Y$ and $x \lessdot y$ or
3. $x \in X$ and $y \in Y$.

**Definition 1.17** [12] Let $X$ and $Y$ be disjoint sets. Let $\mathbf{P} \triangleq (X, \circledcirc)$ and $\mathbf{Q} \triangleq (Y, \lessdot)$ be ordered sets on $X$ and $Y$. The **ordinal product** of $\mathbf{P}$ and $\mathbf{Q}$ is defined as

$$\mathbf{P} \otimes \mathbf{Q} \triangleq (X \times Y, \le)$$

where $(x_1, y_1) \le (x_2, y_2)$ if $\left\{ \begin{array}{llll} 1. & x_1 \ne x_2 & \text{and} & x_1 \circledcirc x_2 \quad \text{or} \\ 2. & x_1 = x_2 & \text{and} & y_1 \lessdot y_2 \end{array} \right\}$

The order relation $\le$ is called a **lexicographical** order relation, **dictionary** order relation, or **alphabetic** order relation.

**Definition 1.18** [13] Let $\mathbf{P} \triangleq (X, \le)$ be an ordered set. Let $\ge$ be the dual order relation of $\le$. The **dual** of $\mathbf{P}$ is defined as $\mathbf{P}^* \triangleq (X, \ge)$

**Definition 1.19** [14] Let $X$ and $Y$ be disjoint sets. Let $\mathbf{P} \triangleq (X, \circledcirc)$ and $\mathbf{Q} \triangleq (Y, \lessdot)$ be ordered sets on $X$ and $Y$. $\mathbf{Q}^{\mathbf{P}} \triangleq \left( \left\{ \mathsf{f} \in Y^X \,\middle|\, \mathsf{f} \text{ is } \textit{order preserving} \right\}, \le \right)$
where $\mathsf{f} \le \mathsf{g}$ if $\mathsf{f}(x) \le \mathsf{g}(x) \quad \forall x \in X$. The order relation $\le$ is called a **pointwise order relation**.

**Theorem 1.20** (cardinal arithmetic) [15] *Let $\mathbf{P} \triangleq (X, \le)$ be an ordered set.*

| | | | | |
|---|---|---|---|---|
| 1. | $\mathbf{P} + \mathbf{Q}$ | $=$ | $\mathbf{Q} + \mathbf{P}$ | (COMMUTATIVE) |
| 2. | $\mathbf{P} \times \mathbf{Q}$ | $=$ | $\mathbf{Q} \times \mathbf{P}$ | (COMMUTATIVE) |
| 3. | $(\mathbf{P} + \mathbf{Q}) + \mathbf{R}$ | $=$ | $\mathbf{P} + (\mathbf{Q} + \mathbf{R})$ | (ASSOCIATIVE) |
| 4. | $(\mathbf{P} \times \mathbf{Q}) \times \mathbf{R}$ | $=$ | $\mathbf{P} \times (\mathbf{Q} \times \mathbf{R})$ | (ASSOCIATIVE) |
| 5. | $\mathbf{P} \times (\mathbf{Q} + \mathbf{R})$ | $=$ | $(\mathbf{P} \times \mathbf{Q}) + (\mathbf{P} \times \mathbf{R})$ | (DISTRIBUTIVE) |
| 6. | $\mathbf{R}^{\mathbf{P}+\mathbf{Q}}$ | $=$ | $\mathbf{R}^{\mathbf{P}} \times \mathbf{R}^{\mathbf{Q}}$ | |
| 7. | $\left(\mathbf{P}^{\mathbf{Q}}\right)^{\mathbf{R}}$ | $=$ | $\mathbf{P}^{\mathbf{Q} \times \mathbf{R}}$ | |

[11] [155], page 100
[12] [155], page 101, [154], page 44, [81], page 58, [82], page 54
[13] [155], page 101
[14] [155], page 101
[15] [155], page 102

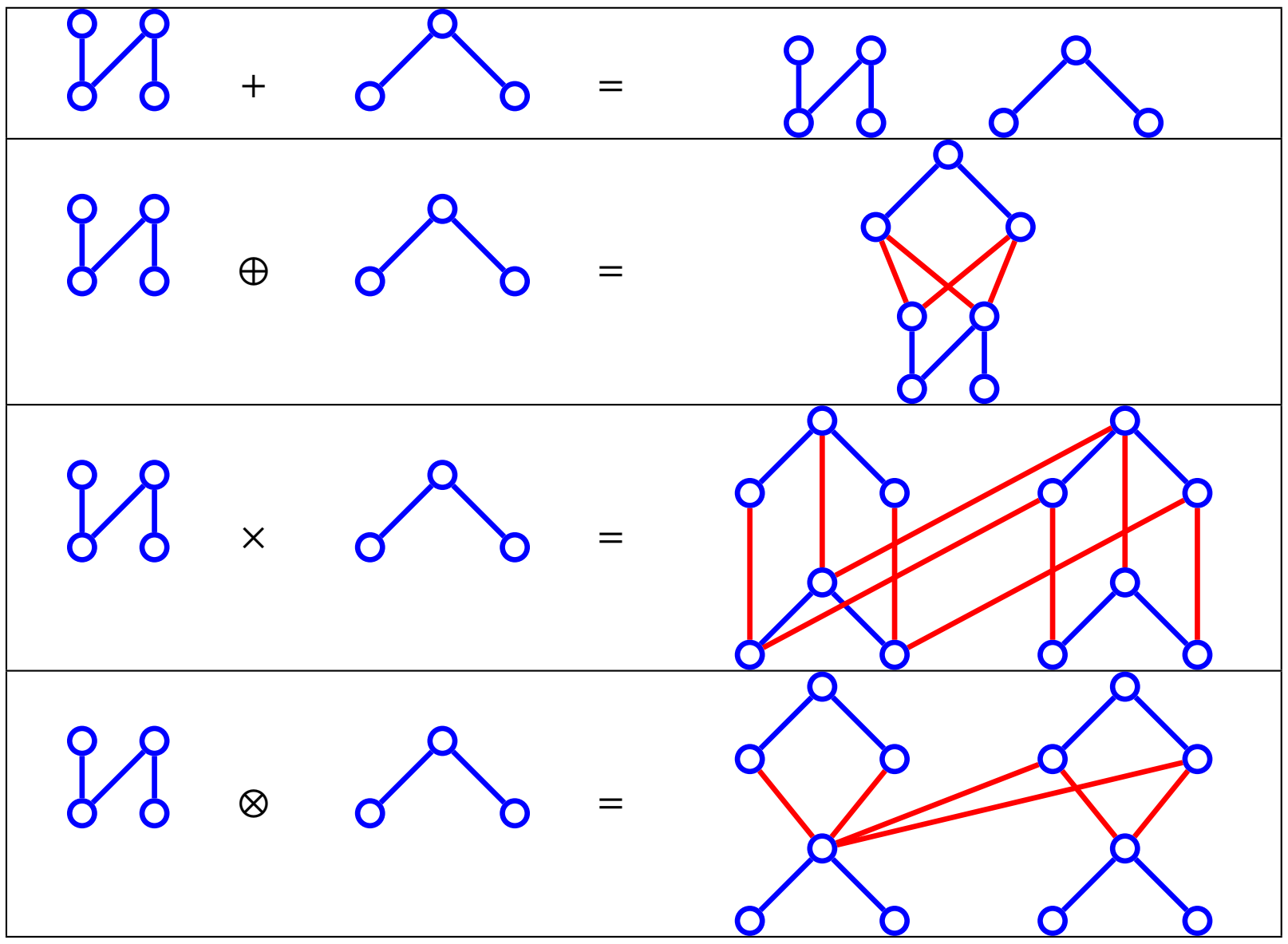

Figure 1: Operations on ordered sets (Example 1.23 page 8)

**Definition 1.21** The ordered set $\boldsymbol{L}_1$ is defined as $(\{x\}, \leq)$, for some value $x$.
It is illustrated by the Hasse diagram to the right.

**Definition 1.22** The ordered set $\boldsymbol{L}_2$ is defined as $\boldsymbol{L}_2 \triangleq \boldsymbol{L}_1^2$.
It is illustrated by the Hasse diagram to the right.

### 1.1.4 Decomposition examples

**Example 1.23** Figure 1 (page 8) illustrates the four ordered set operations $+$, $\times$, $\oplus$, and $\otimes$.

**Example 1.24** [16]The ordered set $n\boldsymbol{L}_1$ is the *anti-chain* with $n$ elements.
The ordered set $4\boldsymbol{L}_1$ is illustrated to the right.

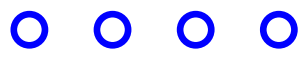

**Example 1.25** The ordered set $\boldsymbol{L}_1^n$ is the *chain* with $n$ elements.
The ordered set $\boldsymbol{L}_1^4$ is illustrated to the right.

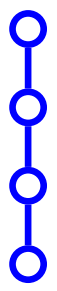

Examples of the *Boolean lattices* (Definition 1.69 page 18) $\boldsymbol{L}_2^1$, $\boldsymbol{L}_2^2$, $\boldsymbol{L}_2^3$, $\boldsymbol{L}_2^4$ and $\boldsymbol{L}_2^5$ are illustrated in Example 1.74 (page 21).

[16] [155], page 100

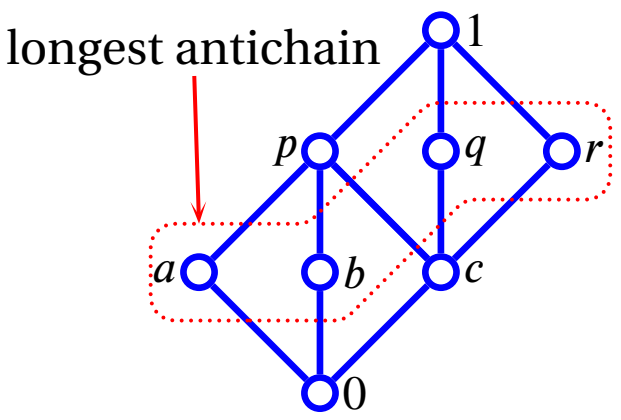

| | | | | |
|---|---|---|---|---|
| partition 1: | $\{0,a\}$, | $\{b,p\}$, | $\{c,q\}$, | $\{r,1\}$ |
| partition 2: | $\{0,a,p\}$, | $\{b\}$, | $\{c,q\}$, | $\{r,1\}$ |
| partition 3: | $\{0,a,p,1\}$, | $\{b\}$, | $\{c,q\}$, | $\{r\}$ |
| partition 4: | $\{0,b,p,1\}$, | $\{a\}$, | $\{c,q\}$, | $\{r\}$ |
| partition 5: | $\{0,c,r,1\}$, | $\{a,p\}$, | $\{b\}$, | $\{q\}$ |
| partition 6: | $\{0,c,q,1\}$, | $\{a,p\}$, | $\{b\}$, | $\{r\}$ |

examples of partitions of chains

Figure 2: Lattice of width 4 and examples of minimal order partitions of chains (see Example 1.26 page 9)

**Example 1.26** [17] The longest *antichain* (Definition 1.11 page 6) in the lattice illustrated in Figure 2 (page 9) has 4 elements giving this ordered set a *width* (Definition 1.12 page 6) of 4. The longest chain also has 4 elements, giving the ordered set a *length* (Definition 1.12 page 6) of 3. By *Dilworth's theorem* (Theorem 1.13 page 6), the smallest partition consists of four *chains* (Definition 1.4 page 4). Examples of such minimal order partitions those listed in Figure 2.

**Definition 1.27** Let $(X, \le)$ be an ordered set and $2^X$ the power set of $X$. For any set $A \in 2^X$, $c$ is an **upper bound** of $A$ in $(X, \le)$ if

1. $x \le c \quad \forall x \in A.$

An element $b$ is the **least upper bound**, or **LUB**, of $A$ in $(X, \le)$ if

2. $b$ and $c$ are *upper bounds* of $A$ $\implies$ $b \le c$.

The least upper bound of the set $A$ is denoted $\bigvee A$. It is also called the **supremum** of $A$, which is denoted $\sup A$. The **join** $x \vee y$ of $x$ and $y$ is defined as $x \vee y \triangleq \bigvee \{x, y\}$.

**Definition 1.28** Let $(X, \le)$ be an ordered set and $2^X$ the power set of $X$. For any set $A \in 2^X$, $p$ is a **lower bound** of $A$ in $(X, \le)$ if

1. $p \le x \quad \forall x \in A.$

An element $a$ is the **greatest lower bound**, or **GLB**, of $A$ in $(X, \le)$ if

2. $a$ and $p$ are *lower bounds* of $A$ $\implies$ $p \le a$.

The greatest lower bound of the set $A$ is denoted $\bigwedge A$. It is also called the **infimum** of $A$, which is denoted $\inf A$. The **meet** $x \wedge y$ of $x$ and $y$ is defined as $x \wedge y \triangleq \bigwedge \{x, y\}$.

**Proposition 1.29** *Let* $(X, \vee, \wedge; \le)$ *be an* ORDERED SET *(Definition 1.1 page 3).*

$$x \le y \iff \left\{ \begin{array}{lrcll} 1. & x \wedge y & = & x & \text{and} \\ 2. & x \vee y & = & y & \end{array} \right\} \quad \forall x,y \in X$$

[17] [56], page 4

**Proposition 1.30** *Let $2^X$ be the* POWER SET *of a set $X$.*

$$A \subseteq B \implies \left\{ \begin{array}{lllll} 1. & \bigvee A & \leq & \bigvee B & \text{and} \\ 2. & \bigwedge A & \leq & \bigwedge B & \end{array} \right\} \qquad \forall A,B \in 2^X$$

## 1.2 Lattices

### 1.2.1 Definition

The structure available in an *ordered set* (Definition 1.1 page 3) tends to be insufficient to ensure "well-behaved" mathematical systems. This situation is greatly remedied if every pair of elements in the ordered set has both a *least upper bound* and a *greatest lower bound* (Definition 1.28 page 9) in the set; in this case, that ordered set is a *lattice* (next definition). Gian-Carlo Rota (1932–1999) has illustrated the advantage of lattices over simple ordered sets by pointing out that the *ordered set of partitions of an integer* "is fraught with pathological properties", while the *lattice* of partitions of a set "remains to this day rich in pleasant surprises".[18]

**Definition 1.31** [19] An algebraic structure $\boldsymbol{L} \triangleq (X, \vee, \wedge; \leq)$ is a **lattice** if

1. $(X, \leq)$ is an ordered set ($(X, \leq)$ is a partially or totally ordered set) and
2. $\exists x \vee y \in X \quad \forall x, y \in X$ (every pair of elements in $X$ has a *least upper bound* in $X$) and
3. $\exists x \wedge y \in X \quad \forall x, y \in X$ (every pair of elements in $X$ has a *greatest lower bound* in $X$).

The algebraic structure $\boldsymbol{L}^* \triangleq (X, \ovee, \owedge; \geq)$ is the **dual** lattice of $\boldsymbol{L}$, where $\ovee$ and $\owedge$ are determined by $\geq$. The *lattice* $\boldsymbol{L}$ is *linear* if $(X, \leq)$ is a *chain* (Definition 1.4 page 4).

**Theorem 1.32** [20] $(X, \vee, \wedge; \leq)$ *is a* LATTICE $\iff$

$$\left\{ \begin{array}{rcl|rcl|lll} x \vee x &=& x & x \wedge x &=& x & \forall x \in X & (\text{IDEMPOTENT}) & \text{and} \\ x \vee y &=& y \vee x & x \wedge y &=& y \wedge x & \forall x,y \in X & (\text{COMMUTATIVE}) & \text{and} \\ (x \vee y) \vee z &=& x \vee (y \vee z) & (x \wedge y) \wedge z &=& x \wedge (y \wedge z) & \forall x,y,z \in X & (\text{ASSOCIATIVE}) & \text{and} \\ x \vee (x \wedge y) &=& x & x \wedge (x \vee y) &=& x & \forall x,y \in X & (\text{ABSORPTIVE}). & \end{array} \right\}$$

**Lemma 1.33** [21] *Let $\boldsymbol{L} \triangleq (X, \vee, \wedge; \leq)$ be* LATTICE *(Definition 1.31 page 10).*

$$x \leq y \iff x = x \wedge y \qquad \forall x,y \in \boldsymbol{L}$$

✎PROOF:

[18] [148], page 1440, ⟨(illustration)⟩, [147], page 498, ⟨partitions of a set⟩
[19] [113], page 473, [17], page 16, [133], [14], page 442, [116], page 1
[20] [113], pages 473–475, ⟨LEMMA 1, THEOREM 4⟩, [23], pages 4–7, [16], pages 795–796, [133], page 409, ⟨($\alpha$)⟩, [14], page 442, [38], pages 371–372, ⟨(1)–(4)⟩
[21] [88]

(1) Proof for $\Longrightarrow$ case: by left hypothesis and definition of $\wedge$ (Definition 1.28 page 9).

(2) Proof for $\Longleftarrow$ case: by right hypothesis and definition of $\wedge$ (Definition 1.28 page 9).

✏

**Proposition 1.34** (Monotony laws) [22] *Let* $(X, \vee, \wedge; \leq)$ *be a lattice.*

$$\left\{\begin{array}{lcll} a & \leq & b & \text{and} \\ x & \leq & y & \end{array}\right\} \implies \left\{\begin{array}{lcll} a \wedge x & \leq & b \wedge y & \text{and} \\ a \vee x & \leq & b \vee y & \end{array}\right\}$$

**Theorem 1.35** (Minimax inequality) [23] *Let* $(X, \vee, \wedge; \leq)$ *be a lattice.*

$$\underbrace{\bigvee_{i=1}^{m} \bigwedge_{j=1}^{n} x_{ij}}_{\text{maxmini: largest of the smallest}} \leq \underbrace{\bigwedge_{j=1}^{n} \bigvee_{i=1}^{m} x_{ij}}_{\text{minimax: smallest of the largest}} \qquad \forall x_{ij} \in X$$

Special cases of the minimax inequality include three distributive *inequalities* (next theorem). If for some lattice any *one* of these inequalities is an *equality*, then *all three* are *equalities* (Theorem 1.54 page 15); and in this case, the lattice is a called a *distributive* lattice (Definition 1.53 page 15).

**Theorem 1.36** (distributive inequalities) [24] $(X, \vee, \wedge; \leq)$ *is a lattice* $\implies$

$$\left\{\begin{array}{rcll} x \wedge (y \vee z) & \geq & (x \wedge y) \vee (x \wedge z) & \forall x,y,z \in X \quad (\text{JOIN SUPER-DISTRIBUTIVE}) \quad \text{and} \\ x \vee (y \wedge z) & \leq & (x \vee y) \wedge (x \vee z) & \forall x,y,z \in X \quad (\text{MEET SUB-DISTRIBUTIVE}) \quad \text{and} \\ (x \wedge y) \vee (x \wedge z) \vee (y \wedge z) & \leq & (x \vee y) \wedge (x \vee z) \wedge (y \vee z) & \forall x,y,z \in X \quad (\text{MEDIAN INEQUALITY}). \end{array}\right.$$

Besides the distributive property, another consequence of the minimax inequality is the *modularity inequality* (next theorem). A lattice in which this inequality becomes equality is said to be *modular* (Definition 1.47 page 14).

**Theorem 1.37** (Modular inequality) [25] *Let* $(X, \vee, \wedge; \leq)$ *be a* LATTICE *(Definition 1.31 page 10).*

$$x \leq y \implies x \vee (y \wedge z) \leq y \wedge (x \vee z)$$

Theorem 1.32 (page 10) gives 4 necessary and sufficient pairs of properties for a structure $(X, \vee, \wedge; \leq)$ to be a *lattice*. However, these 4 pairs are actually *overly* sufficient (they are not *independent*), as demonstrated next.

[22] [68], page 39, [50], pages 97–99, [78], ⟨§4.2⟩

[23] [17], pages 19–20

[24] [36], page 85, [72], page 38, [14], page 444, [105], page 157, [125], page 13, ⟨terminology⟩

[25] [17], page 19, [23], page 11, [38], page 374

**Theorem 1.38** [26]
$(X, \vee, \wedge; \leq)$ *is a lattice* $\iff$

$$\left\{\begin{array}{rcl|rcl|lll} x\vee y &=& y\vee x & x\wedge y &=& y\wedge x & \forall x,y\in X & (\text{COMMUTATIVE}) & \text{and}\\ (x\vee y)\vee z &=& x\vee(y\vee z) & (x\wedge y)\wedge z &=& x\wedge(y\wedge z) & \forall x,y,z\in X & (\text{ASSOCIATIVE}) & \text{and}\\ x\vee(x\wedge y) &=& x & x\wedge(x\vee y) &=& x & \forall x,y\in X & (\text{ABSORPTIVE}) & \end{array}\right\}$$

### 1.2.2 Bounded lattices

Let $\boldsymbol{L} \triangleq (X, \vee, \wedge; \leq)$ be a lattice. By the definition of a *lattice* (Definition 1.31 page 10), the *upper bound* $(x \vee y)$ and *lower bound* $(x \wedge y)$ of any two elements in $X$ is also in $X$. But what about the upper and lower bounds of the entire set $X$ ($\bigvee X$ and $\bigwedge X$) (Definition 1.27 page 9, Definition 1.28 page 9)? If both of these are in $X$, then the lattice $\boldsymbol{L}$ is said to be *bounded* (next definition). All *finite* lattices are bounded (next proposition). However, not all lattices are bounded—for example, the lattice $(\mathbb{Z}, \leq)$ (the lattice of integers with the standard integer ordering relation) is *unbounded*.

**Definition 1.39** Let $\boldsymbol{L} \triangleq (X, \vee, \wedge; \leq)$ be a lattice. Let $\bigvee X$ be the least upper bound of $(X, \leq)$ and let $\bigwedge X$ be the greatest lower bound of $(X, \leq)$.

$\boldsymbol{L}$ is **upper bounded** if $\left(\bigvee X\right) \in X$.
$\boldsymbol{L}$ is **lower bounded** if $\left(\bigwedge X\right) \in X$.
$\boldsymbol{L}$ is **bounded** if $\boldsymbol{L}$ is both upper and lower bounded.

A *bounded* lattice is optionally denoted $(X, \vee, \wedge, 0, 1; \leq)$, where $0 \triangleq \bigwedge X$ and $1 \triangleq \bigvee X$.

**Proposition 1.40** *Let* $\boldsymbol{L} \triangleq (X, \vee, \wedge; \leq)$ *be a lattice.*

$$\{\boldsymbol{L} \text{ is FINITE}\} \implies \{\boldsymbol{L} \text{ is BOUNDED}\}$$

**Proposition 1.41** [27] *Let* $\boldsymbol{L} \triangleq (X, \vee, \wedge; \leq)$ *be a lattice with* $\bigvee X \triangleq 1$ *and* $\bigwedge X \triangleq 0$.

$$\{\boldsymbol{L} \text{ is BOUNDED}\} \implies \left\{\begin{array}{rcllll} x\vee 1 &=& 1 & \forall x\in X & (\textit{upper bounded}) & \textit{and}\\ x\wedge 0 &=& 0 & \forall x\in X & (\textit{lower bounded}) & \textit{and}\\ x\vee 0 &=& x & \forall x\in X & (\textit{join-identity}) & \textit{and}\\ x\wedge 1 &=& x & \forall x\in X & (\textit{meet-identity}) & \end{array}\right\}$$

**Definition 1.42** [28] Let $\boldsymbol{L} \triangleq (X, \vee, \wedge, 0, 1; \leq)$ be a *bounded lattice* (Definition 1.39 page 12). The **height** $\mathsf{h}(x)$ of a point $x \in \boldsymbol{L}$ is the *least upper bound* of the *lengths* (Definition 1.12 page 6) of all the *chains* that have 0 and in which $x$ is the *least upper bound*. The **height** $\mathsf{h}(\boldsymbol{L})$ of the lattice $\boldsymbol{L}$ is defined as

$$\mathsf{h}(\boldsymbol{L}) \triangleq \mathsf{h}(1).$$

---

26 [136], pages 7–8, [12], page 5, [120], page 24, [77], ⟨Theorem 1.22⟩, [78], ⟨§4.4⟩
27 [77], ⟨§1.2.2⟩, [78], ⟨§4.5⟩
28 [18], page 5

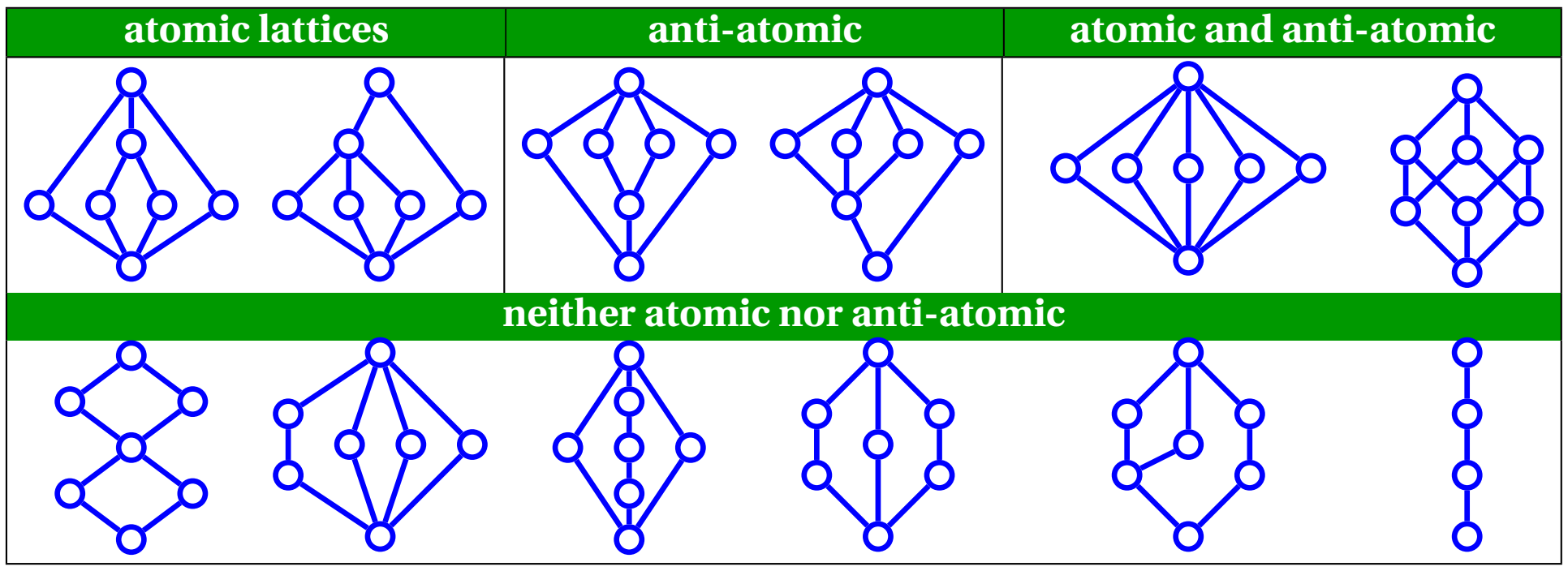

Figure 3: Selected *atomic*, *anti-atomic*, and neither atomic nor anti-atomic lattices (see Example 1.45 page 13)

**Example 1.43** The *height* of the lattice illustrated in Figure 2 (page 9) is 3 because

$$\begin{aligned}
\mathsf{h}(\boldsymbol{L}) &\triangleq \mathsf{h}(1) \\
&\triangleq \bigvee \left\{ \ell(\boldsymbol{C}) \,\middle|\, \boldsymbol{C} \text{ is a } \textit{chain} \text{ in } \boldsymbol{L} \text{ containing both 0 and 1} \right\} \\
&= \bigvee \left\{ \ell(\{0,a,p,1\},\leq),\ \ell(\{0,b,p,1\},\leq),\ \ell(\{0,c,p,1\},\leq),\ \ell(\{0,c,q,1\},\leq), \right. \\
&\qquad \left. \ell(\{0,c,r,1\},\leq), \right\} \\
&= \bigvee \{4-1, 4-1, 4-1, 4-1, 4-1\} \\
&= \bigvee \{3,3,3,3,3\} \\
&= 3
\end{aligned}$$

### 1.2.3 Atomic lattices

**Definition 1.44** [29] Let $\boldsymbol{L} \triangleq (X, \vee, \wedge, 0, 1; \leq)$ be a *bounded lattice* (Definition 1.39 page 12).
$x$ is an **atom** of $\boldsymbol{L}$ if $x$ *covers* (Definition 1.5 page 4) 0.
$x$ is an **anti-atom** of $\boldsymbol{L}$ if $x$ is *covered by* 1.
$\boldsymbol{L}$ is **atomic** if every $x \in X \backslash 0$ can be represented as joins of atoms of $\boldsymbol{L}$.
$\boldsymbol{L}$ is **anti-atomic** if every $x \in X \backslash 1$ can be represented as meets of anti-atoms of $\boldsymbol{L}$.

**Example 1.45** Figure 3 (page 13) illustrates some examples of lattices that are *atomic*, *anti-atomic*, both, and neither.

[29] [108], page 178, [16], page 800, ⟨see footnote ‡⟩

### 1.2.4 Modular Lattices

**Definition 1.46** [30] Let $(X, \vee, \wedge; \le)$ be a lattice. Let $2^{XX}$ be the set of all *relations* in $X^2$. The **modularity** relation $Ⓜ \in 2^{XX}$ and the **dual modularity** relation $Ⓜ^* \in 2^{XX}$ are defined as

$$x Ⓜ y \overset{\text{def}}{\iff} \left\{(x,y) \in X^2 \,|\, a \le y \implies y \wedge (x \vee a) = (y \wedge x) \vee a \quad \forall a \in X\right\}$$
$$x Ⓜ^* y \overset{\text{def}}{\iff} \left\{(x,y) \in X^2 \,|\, a \ge y \implies y \vee (x \wedge a) = (y \vee x) \wedge a \quad \forall a \in X\right\}.$$

A pair $(x,y) \in Ⓜ$ is alternatively denoted as $(x,y)Ⓜ$, and is called a **modular** pair. A pair $(x,y) \in Ⓜ^*$ is alternatively denoted as $(x,y)Ⓜ^*$, and is called a **dual modular** pair. A pair $(x,y)$ that is *not* a modular pair ($(x,y) \notin Ⓜ$) is denoted $x \not{Ⓜ} y$. A pair $(x,y)$ that is *not* a dual modular pair is denoted $x \not{Ⓜ}^* y$.

Modular lattices are a generalization of *distributive lattices* (Definition 1.53 page 15) in that all distributive lattices are modular, but not all modular lattices are distributive (Example 1.61 page 16, Example 1.62 page 17).

**Definition 1.47** [31] A lattice $(X, \vee, \wedge; \le)$ is **modular** if $\quad x Ⓜ y \qquad \forall x, y \in X$.

**Theorem 1.48** [32] *Let* $\mathbf{L} \triangleq (X, \vee, \wedge; \le)$ *be a lattice.*

$$\begin{array}{lclcll}
\mathbf{L} \text{ is MODULAR} & \iff & \{x \le y & \implies & x \vee (z \wedge y) = (x \vee z) \wedge y\} & \forall x,y,z \in X \\
& \iff & x \vee [(x \vee y) \wedge z] & = & (x \vee y) \wedge (x \vee z) & \forall x,y,z \in X \\
& \iff & x \wedge [(x \wedge y) \vee z] & = & (x \wedge y) \vee (x \wedge z) & \forall x,y,z \in X
\end{array}$$

**Definition 1.49** (N5 lattice/pentagon) [33] The **N5 lattice** is the ordered set $(\{0, a, b, p, 1\}, \le)$ with cover relation
$\lessdot = \{(0,a), (a,b), (b,1), (p,1), (0,p)\}$.
The N5 lattice is also called the **pentagon**. The N5 lattice is illustrated by the Hasse diagram to the right.

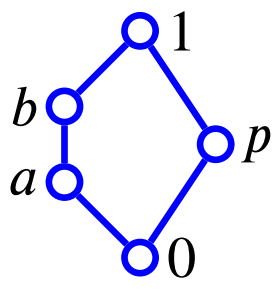

**Theorem 1.50** [34] *Let* $\mathbf{L}$ *be a* LATTICE *(Definition 1.31 page 10).*
$\mathbf{L}$ *is* MODULAR *(Definition 1.47 page 14)* $\iff$ $\mathbf{L}$ *does* NOT *contain the* N5 LATTICE *(Definition 1.49 page 14).*

**Theorem 1.51** [35] *Let* $\mathbf{A} \triangleq (X, \vee, \wedge; \le)$ *be an algebraic structure.*

$$\left\{\begin{array}{lcll}
(x \wedge y) \vee (x \wedge z) & = & [(z \wedge x) \vee y] \wedge x & \forall x,y,z \in X \quad \textit{and} \\
[x \vee (y \vee z)] \wedge z & = & z & \forall x,y,z \in X
\end{array}\right\} \iff \left\{\begin{array}{l}\mathbf{A} \textit{ is a} \\ \textit{modular lattice}\end{array}\right\}$$

[30] [157], page 11, [116], page 1, ⟨Definition (1.1)⟩, [117], page 248
[31] [18], page 82, [116], page 3, ⟨Definition (1.7)⟩
[32] [136], page 39, [133], page 413, ⟨(2)⟩, [78], ⟨Theorem 5.1⟩
[33] [12], pages 12–13, [38], pages 391–392, ⟨(44) and (45)⟩
[34] [23], page 11, [71], page 70, [38], ⟨cf Stern 1999 page 10⟩, [78], ⟨Theorem 5.1⟩
[35] [136], pages 42–43, [145]

Examples of *modular lattice*s are provided in Example 1.61 (page 16) and Example 1.62 (page 17).

### 1.2.5 Distributive Lattices

**Definition 1.52** [36] Let $(X, \vee, \wedge; \leq)$ be a *lattice* (Definition 1.31 page 10). Let $2^{XXX}$ be the set of all *relation*s in $X^3$. The **distributivity** relation $\mathbb{D} \in 2^{XXX}$ and the **dual distributivity** relation $\mathbb{D}^* \in 2^{XXX}$ are defined as

$$\mathbb{D} \triangleq \left\{(x,y,z) \in X^3 \,|\, x \wedge (y \vee z) = (x \wedge y) \vee (x \wedge z)\right\} \quad \text{(each } (x,y,z) \text{ is } \textit{disjunctive distributive}\text{) and}$$

$$\mathbb{D}^* \triangleq \left\{(x,y,z) \in X^3 \,|\, x \vee (y \wedge z) = (x \vee y) \wedge (x \vee z)\right\} \quad \text{(each } (x,y,z) \text{ is } \textit{conjunctive distributive}\text{).}$$

A triple $(x,y,z) \in \mathbb{D}$ is alternatively denoted as $(x,y,z)\,\mathbb{D}$, and is a **distributive** triple. A triple $(x,y,z) \in \mathbb{D}^*$ is alternatively denoted as $(x,y,z)\,\mathbb{D}^*$, and is a **dual distributive** triple.

**Definition 1.53** [37] A lattice $(X, \vee, \wedge; \leq)$ is **distributive** if $\quad (x,y,z) \in \mathbb{D} \quad \forall x,y,z \in X$

Not all lattices are *distributive*. But if a lattice $\boldsymbol{L}$ does happen to be distributive (Definition 1.53 page 15)—that is all triples in $\boldsymbol{L}$ satisfy the *distributive* property (Definition 1.53 page 15)—then all triples in $\boldsymbol{L}$ also satisfy the *dual distributive* property, as well as another property called the *median property*. The converses also hold (next theorem).

**Theorem 1.54** [38] *Let* $\boldsymbol{L} \triangleq (X, \vee, \wedge; \leq)$ *be a* LATTICE *(Definition 1.31 page 10).*
$\boldsymbol{L}$ *is* DISTRIBUTIVE *(Definition 1.53 page 15)*

$$\iff x \wedge (y \vee z) = (x \wedge y) \vee (x \wedge z) \qquad \forall x,y,z \in X \qquad \text{(DISJUNCTIVE DISTRIBUTIVE)}$$

$$\iff x \vee (y \wedge z) = (x \vee y) \wedge (x \vee z) \qquad \forall x,y,z \in X \qquad \text{(CONJUNCTIVE DISTRIBUTIVE)}$$

$$\iff (x \vee y) \wedge (x \vee z) \wedge (y \vee z) = (x \wedge y) \vee (x \wedge z) \vee (y \wedge z) \qquad \forall x,y,z \in X \qquad \text{(MEDIAN PROPERTY)}$$

**Definition 1.55** (M3 lattice/diamond) [39] The **M3 lattice** is the ordered set $(\{0,p,q,r,1\}, \leq)$ with covering relation

$$<= \{(p,1), (q,1), (r,1), (0,p), (0,q), (0,r)\}.$$

The M3 lattice is also called the **diamond**, and is illustrated by the Hasse diagram to the right.

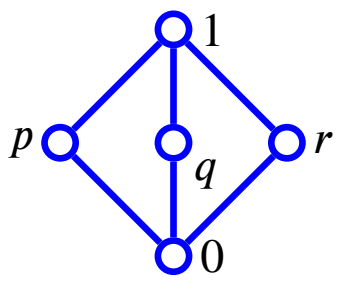

[36] [116], page 15, ⟨Definition 4.1⟩, [62], page 67, [130], page 32, ⟨Definition 5.1⟩, [37], page 314, ⟨*disjunctive distributive* and *conjunctive distributive* functions⟩

[37] [23], page 10, [17], page 133, [133], page 414, ⟨*arithmetic axiom*⟩, [14], page 453, [9], page 48, ⟨Definition II.5.1⟩

[38] [49], page 237, [23], page 10, [133], page 416, ⟨(7),(8), Theorem 3⟩, [134], ⟨cf Gratzer 2003 page 159⟩, [153], page 286, ⟨cf Birkhoff(1948)p.133⟩, [105], ⟨cf Birkhoff(1948)p.133⟩, [78], ⟨Theorem 6.1⟩

[39] [12], pages 12–13, [105], page 157, ⟨ $p_1 \equiv x$, $p_2 \equiv y$, $p_3 \equiv z$, $g \equiv 1$, $0 \equiv 0$ ⟩

**Lemma 1.56** [40]

$$\left\{ \begin{array}{l} \mathbf{L} \text{ is an} \\ \textbf{M3 lattice} \end{array} \right\} \implies \left\{ \begin{array}{lll} 1. & \mathbf{L} \text{ is NOT distributive} & \text{(Definition 1.53 page 15)} \quad \text{and} \\ 2. & \mathbf{L} \text{ IS modular} & \text{(Definition 1.47 page 14)} \end{array} \right\}$$

**Theorem 1.57** (Birkhoff distributivity criterion) [41] *Let* $\mathbf{L} \triangleq (X, \vee, \wedge; \leq)$ *be a* LATTICE.

$$\mathbf{L} \text{ is DISTRIBUTIVE} \iff \left\{ \begin{array}{l} \mathbf{L} \text{ does } \textbf{not} \text{ contain N5 as a sublattice} \quad \text{and} \\ \mathbf{L} \text{ does } \textbf{not} \text{ contain M3 as a sublattice} \end{array} \right.$$

Distributive lattices are a special case of modular lattices. That is, all distributive lattices are modular, but not all modular lattices are distributive (next theorem). An example is the *M3 lattice*—it is modular, but yet it is not *distributive*.

**Theorem 1.58** [42] *Let* $(X, \vee, \wedge; \leq)$ *be a lattice.*

$$\{(X, \vee, \wedge; \leq) \text{ is DISTRIBUTIVE}\} \quad \overset{\Longrightarrow}{\not\Longleftarrow} \quad \{(X, \vee, \wedge; \leq) \text{ is MODULAR}\}$$

**Theorem 1.59** [43] *Let* $\mathbf{L} \triangleq (X, \vee, \wedge; \leq)$ *be a* LATTICE *(Definition 1.31 page 10)*.

$$\left\{ \begin{array}{lll} 1. & \mathbf{L} \text{ is DISTRIBUTIVE} & \text{and} \\ 2. & x \vee a = x \vee b & \text{and} \\ 3. & x \wedge a = x \wedge b & \end{array} \right\} \implies \{a = b\} \qquad \forall x, a, b \in X$$

**Proposition 1.60** [44] *Let* $X_n$ *be a finite set with order* $n = |X_n|$. *Let* $l_n$ *be the number of unlabeled lattices on* $X_n$, $m_n$ *the number of unlabeled modular lattices on* $X_n$, *and* $d_n$ *the number of unlabeled distributive lattices on* $X_n$.

| $n$ | 0 | 1 | 2 | 3 | 4 | 5 | 6 | 7 | 8 | 9 | 10 | 11 | 12 | 13 | 14 |
|---|---|---|---|---|---|---|---|---|---|---|---|---|---|---|---|
| $l_n$ | 1 | 1 | 1 | 1 | 2 | 5 | 15 | 53 | 222 | 1078 | 5994 | 37622 | 262776 | 2018305 | 16873364 |
| $m_n$ | 1 | 1 | 1 | 1 | 2 | 4 | 8 | 16 | 34 | 72 | 157 | 343 | 766 | 1718 | 3899 |
| $d_n$ | 1 | 1 | 1 | 1 | 2 | 3 | 5 | 8 | 15 | 26 | 47 | 82 | 151 | 269 | 494 |

**Example 1.61** [45] There are a total of 5 unlabeled lattices on a five element set. Of these, 3 are *distributive* (Proposition 1.60 page 16, and thus also *modular*), one is *modular* but *non-*

[40] [17], page 6, [23], page 11, [105], page 157, ⟨cf Salii1988 p. 37⟩
[41] [23], page 12, [17], page 134, [19] [78], ⟨Theorem 6.2⟩
[42] [17], page 134, [23], page 11 [77], ⟨Theorem 1.37⟩, [78], ⟨§6.2.3⟩
[43] [113], pages 484–485
[44] [2] ⟨http://oeis.org/A006966⟩, [2] ⟨http://oeis.org/A006982⟩, [2] ⟨http://oeis.org/A006981⟩, [84], ⟨$l_n$⟩, [54], page 17, ⟨$d_n$⟩, [160]
[45] [54], pages 4–5, [78], ⟨Example 6.2⟩

*distributive*, and one is *non-distributive* (and *non-modular*).

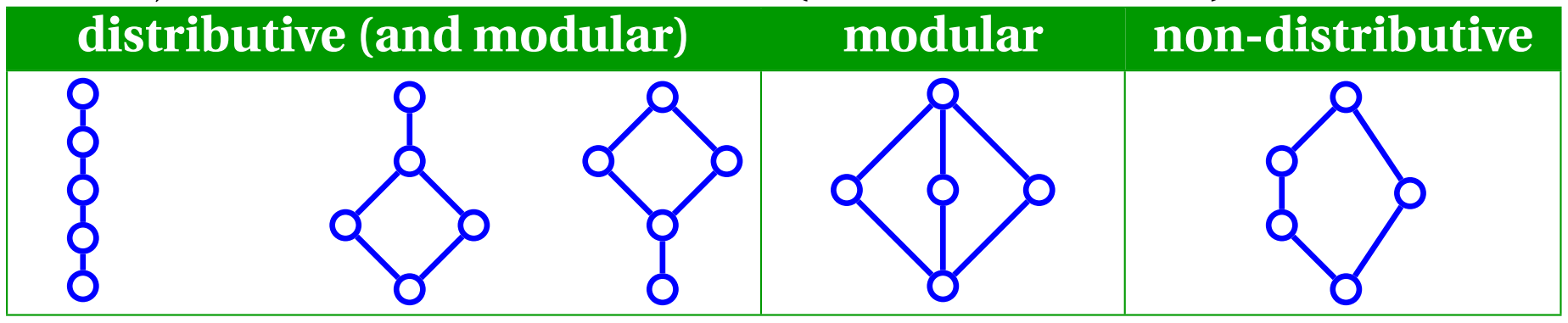

**Example 1.62** [46] There are a total of 15 unlabeled lattices on a six element set. Of these, 5 are *distributive* (Proposition 1.60 page 16, and *modular*), 3 are *modular* but *non-distributive*, and 7 are *non-distributive* (and *non-modular*).

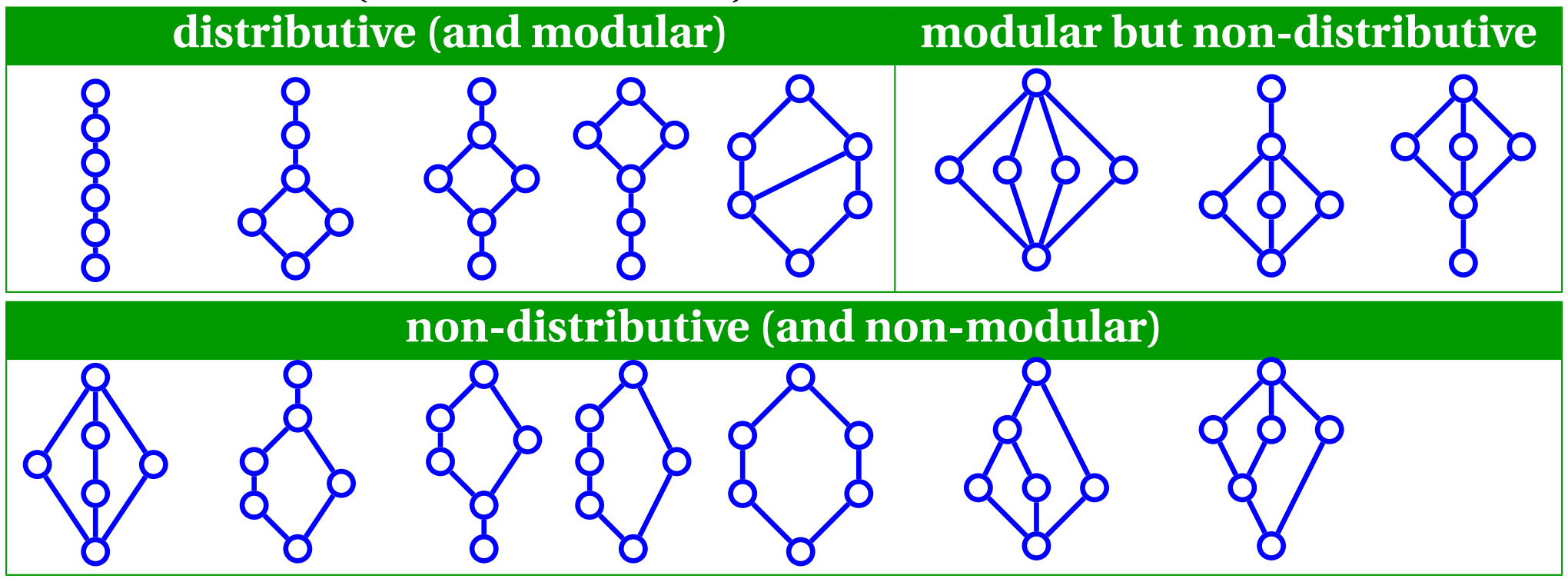

### 1.2.6 Complemented lattices

**Definition 1.63** [47] Let $\boldsymbol{L} \triangleq (X, \vee, \wedge, 0, 1; \leq)$ be a *bounded lattice* (Definition 1.39 page 12). An element $x' \in X$ is a **complement** of an element $x$ in $\boldsymbol{L}$ if

1. $x \wedge x' = 0$ (*non-contradiction*) and
2. $x \vee x' = 1$ (*excluded middle*).

An element $x'$ in $\boldsymbol{L}$ is the *unique complement* of $x$ in $\boldsymbol{L}$ if $x'$ is a *complement* of $x$ and $y'$ is a *complement* of $x \implies x' = y'$. $\boldsymbol{L}$ is **complemented** if every element in $X$ has a complement in $X$. $\boldsymbol{L}$ is **uniquely complemented** if every element in $X$ has a unique complement in $X$. A complemented lattice that is *not* uniquely complemented is **multiply complemented**.

**Example 1.64** Here are some examples:

[46] [78], ⟨Example 5.6⟩
[47] [157], page 9, [17], page 23

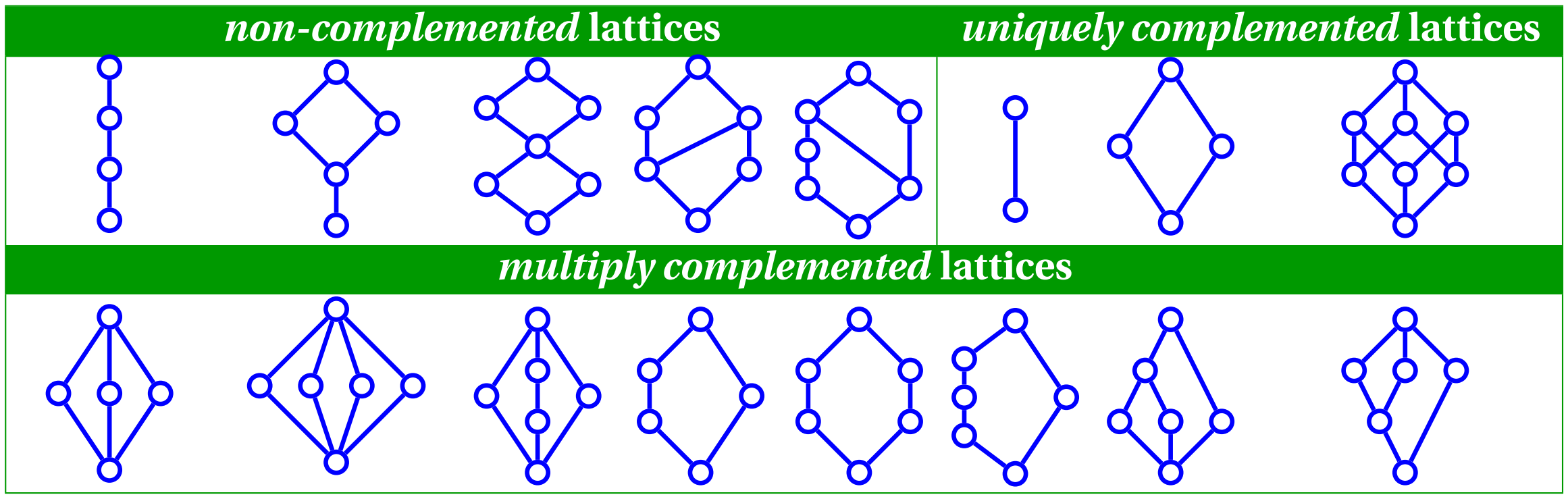

**Example 1.65** Of the 53 unlabeled lattices on a 7 element set, 0 are *uniquely complemented*, 17 are *multiply complemented*, and 36 are *non-complemented*.

Theorem 1.66 (next) is a landmark theorem in mathematics.

**Theorem 1.66** [48] *For every lattice* $\boldsymbol{L}$*, there exists a lattice* $\boldsymbol{U}$ *such that*

1. $\boldsymbol{L} \subseteq \boldsymbol{U}$ *(*$\boldsymbol{L}$ *is a sublattice of* $\boldsymbol{U}$*)* *and*
2. $\boldsymbol{U}$ *is* UNIQUELY COMPLEMENTED.

**Corollary 1.67** [49] *Let* $\boldsymbol{L} \triangleq (X, \vee, \wedge; \le)$ *be a lattice.*

$$\left\{\begin{array}{ll} 1. & \boldsymbol{L} \text{ is DISTRIBUTIVE} \quad \textit{and} \\ 2. & \boldsymbol{L} \text{ is COMPLEMENTED} \end{array}\right\} \begin{array}{c}\Longrightarrow \\ \not\Longleftarrow\end{array} \{\boldsymbol{L} \text{ is UNIQUELY COMPLEMENTED}\}$$

**Theorem 1.68** (Huntington properties) [50] *Let* $\boldsymbol{L}$ *be a lattice.*

$$\left\{\begin{array}{l} \boldsymbol{L} \text{ is} \\ \text{UNIQUELY} \\ \text{COMPLEMENTED} \end{array}\right\} \textit{and} \underbrace{\left\{\begin{array}{ll} \boldsymbol{L} \text{ is MODULAR} & \textit{or} \\ \boldsymbol{L} \text{ is ATOMIC} & \textit{or} \\ \boldsymbol{L} \text{ is ORTHOCOMPLEMENTED} & \textit{or} \\ \boldsymbol{L} \text{ has FINITE WIDTH} & \textit{or} \\ \boldsymbol{L} \text{ is DE MORGAN} & \end{array}\right\}}_{\text{HUNTINGTON PROPERTIES}} \Longrightarrow \left\{\begin{array}{l} \boldsymbol{L} \text{ is} \\ \text{DISTRIBUTIVE} \end{array}\right\}$$

### 1.2.7 Boolean lattices

**Definition 1.69** [51] A *lattice* (Definition 1.31 page 10) $\boldsymbol{L}$ is **Boolean** if

1. $\boldsymbol{L}$ is *bounded* (Definition 1.39 page 12) and
2. $\boldsymbol{L}$ is *distributive* (Definition 1.53 page 15) and
3. $\boldsymbol{L}$ is *complemented* (Definition 1.63 page 17).

---

[48] [46], page 123, [151], page 51, [72], page 378, ⟨Corollary 3.8⟩
[49] [113], page 488, [151], page 30, ⟨Theorem 10⟩
[50] [146], page 103, [3], page 79, [151], page 40, [46], page 123, [73], page 698
[51] [113], page 488, [97]

In this case, $\boldsymbol{L}$ is a **Boolean algebra** or a **Boolean lattice**.
In this paper, a *Boolean lattice* with $2^N$ elements is sometimes denoted $\boldsymbol{L}_2^N$.

The next theorem presents the classic properties of any Boolean algebra. The first 4 pairs of properties are true for any lattice (Theorem 1.32 page 10). The *bounded*, *distributive*, and *complemented* properties are true by definition of a *Boolean lattice* (Definition 1.69 page 18).

**Theorem 1.70** (classic 10 Boolean properties) [52] *Let* $\boldsymbol{A} \triangleq (X, \vee, \wedge, 0, 1; \leq)$ *be an algebraic structure. In the event that* $\boldsymbol{A}$ *is a* BOUNDED LATTICE *(Definition 1.39 page 12), let* $x'$ *represent a* COMPLEMENT *(Definition 1.63 page 17) of an element* $x$ *in* $\boldsymbol{A}$.
$\boldsymbol{A}$ *is a* ***Boolean algebra*** $\iff$ $\forall x, y, z \in X$

| disjunctive properties | conjunctive properties | property name | |
|---|---|---|---|
| $x \vee x = x$ | $x \wedge x = x$ | (IDEMPOTENT) | *and* |
| $x \vee y = y \vee x$ | $x \wedge y = y \wedge x$ | (COMMUTATIVE) | *and* |
| $x \vee (y \vee z) = (x \vee y) \vee z$ | $x \wedge (y \wedge z) = (x \wedge y) \wedge z$ | (ASSOCIATIVE) | *and* |
| $x \vee (x \wedge y) = x$ | $x \wedge (x \vee y) = x$ | (ABSORPTIVE) | *and* |
| $x \vee 1 = 1$ | $x \wedge 0 = 0$ | (BOUNDED) | *and* |
| $x \vee 0 = x$ | $x \wedge 1 = x$ | (IDENTITY) | *and* |
| $x \vee (y \wedge z) = (x \vee y) \wedge (x \vee z)$ | $x \wedge (y \vee z) = (x \wedge y) \vee (x \wedge z)$ | (DISTRIBUTIVE) | *and* |
| $x \vee x' = 1$ | $x \wedge x' = 0$ | (COMPLEMENTED) | *and* |
| $(x \vee y)' = x' \wedge y'$ | $(x \wedge y)' = x' \vee y'$ | (DE MORGAN) | *and* |
| $(x')' = x$ | | (INVOLUTORY) | |

**Proposition 1.71** (Huntington's fourth set) [53] *Let* $\boldsymbol{A} \triangleq (X, \vee, \wedge; \leq)$ *be an* ALGEBRAIC STRUCTURE. $\boldsymbol{A}$ *is a* ***Boolean algebra*** $\iff$

$$\left\{\begin{array}{lllll}
1. & x \vee x = x & \forall x \in X & \text{(IDEMPOTENT)} & \text{and} \\
2. & x \vee y = y \vee x & \forall x,y \in X & \text{(COMMUTATIVE)} & \text{and} \\
3. & (x \vee y) \vee z = x \vee (y \vee z) & \forall x,y,z \in X & \text{(ASSOCIATIVE)} & \text{and} \\
4. & \left(x' \vee y'\right)' \vee \left(x' \vee y\right)' = x & \forall x,y \in X. & \text{(HUNTINGTON'S AXIOM)} &
\end{array}\right\}$$

## 1.3 Orthocomplemented Lattices

*Orthocomplemented lattice*s (Definition 1.72 page 20) are a kind of generalization of *Boolean algebra*s. The relationship between lattices of several types, including orthocomplemented and Boolean lattices, is stated in Theorem 1.86 (page 26) and illustrated in Figure 4 (page 20).

[52] [89], pages 292–293, ⟨"1st set"⟩, [90], page 280, ⟨"4th set"⟩, [113], page 488, [68], page 10, [124], pages 20–21, [153], [167], pages 35–37
[53] [90], page 280, ⟨"4th set"⟩

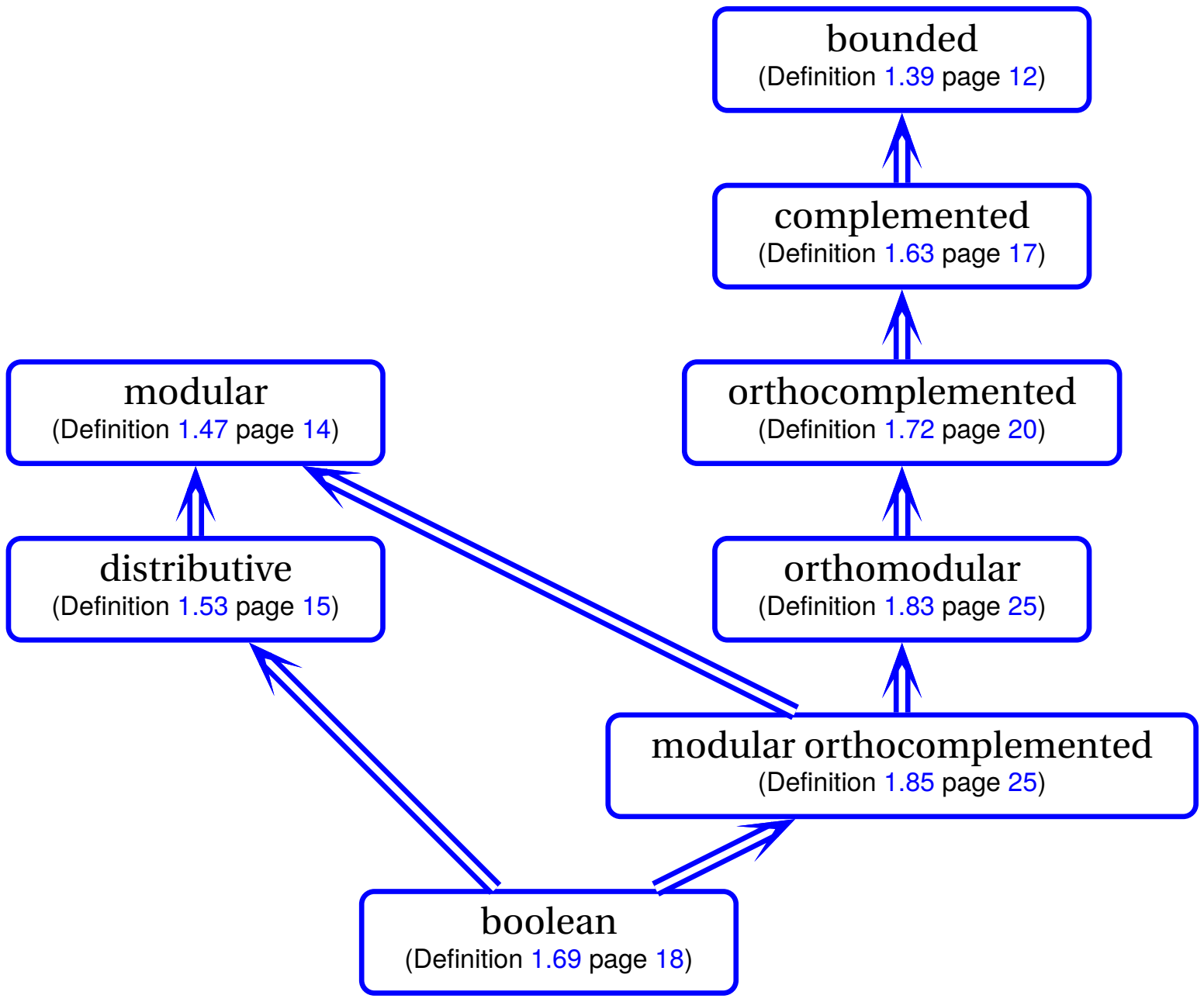

Figure 4: relationships between selected lattice types (see Theorem 1.86 page 26)

### 1.3.1 Definition

**Definition 1.72** [54] Let $\boldsymbol{L} \triangleq (X, \vee, \wedge, 0, 1; \leq)$ be a *bounded lattice* (Definition 1.39 page 12). An element $x^\perp \in X$ is an **orthocomplement** of an element $x \in X$ if

1. $x^{\perp\perp} = x \quad \forall x \in X$ (*involutory*) and
2. $x \wedge x^\perp = 0 \quad \forall x \in X$ (*non-contradiction*) and
3. $x \leq y \implies y^\perp \leq x^\perp \quad \forall x, y \in X$ (*antitone*).

The lattice $\boldsymbol{L}$ is **orthocomplemented** ($\boldsymbol{L}$ is an **orthocomplemented lattice**) if every element $x$ in $X$ has an *orthocomplement*. The elements $\{x, y\}$ are **orthocomplemented pairs** in $\boldsymbol{L}$ if $y = x^\perp$.

**Definition 1.73** [55]
The $\mathbf{O_6}$ **lattice** is the ordered set $\left(\{0, p, q, p^\perp, q^\perp, 1\}, \leq\right)$ with cover relation
$\lessdot = \left\{(0,p), (0,q), (p,q^\perp), (q,p^\perp), (p^\perp,1), (q^\perp,1)\right\}$.
The $O_6$ lattice is illustrated by the Hasse diagram to the right.

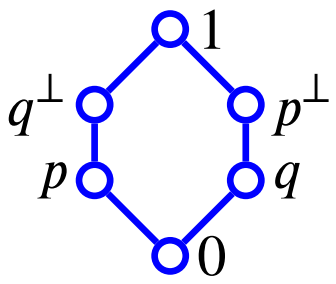

[54] [157], page 11, [12], page 28, [98], page 16, [79], page 76, [112], page 3, [20], page 830, ⟨L71–L73⟩

[55] [98], page 22, [88], page 50, [12], page 33, [157], page 12. The $O_6$ *lattice* is also called the **Benzene ring** or the **hexagon**.

**Example 1.74** [56] There are a total of 10 **orthocomplemented lattices** with 8 elements or less. These 10, along with 3 other orthocomplemented lattices with 10 elements, are illustrated next:

Lattices that are **orthocomplemented** but *non-orthomodular* and hence also *not modular orthocomplemented* and *non-Boolean*:

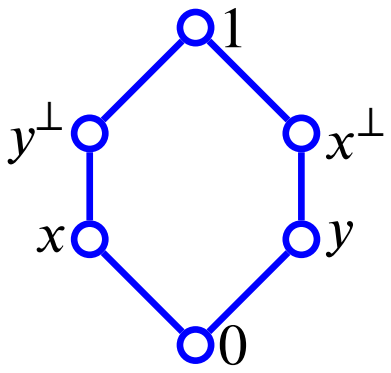

1. $O_6$ *lattice*

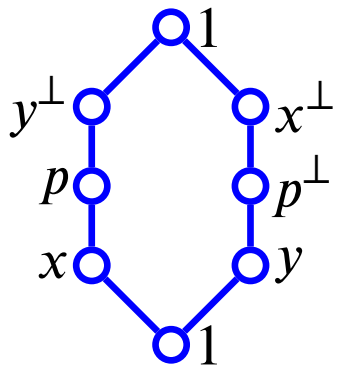

2. $O_8$ *lattice*

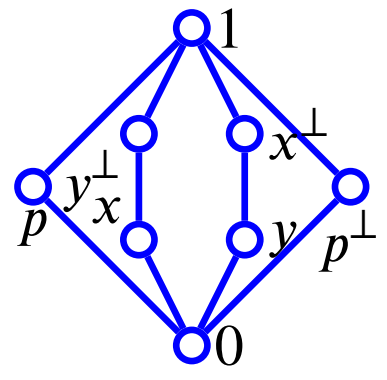

3.

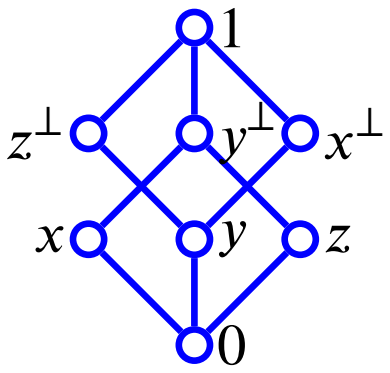

4.

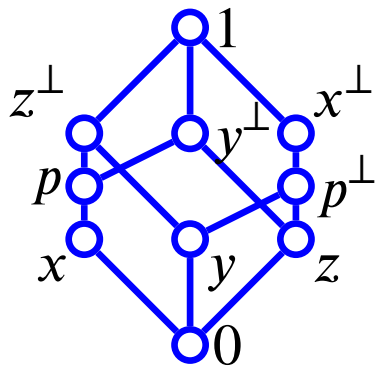

5.

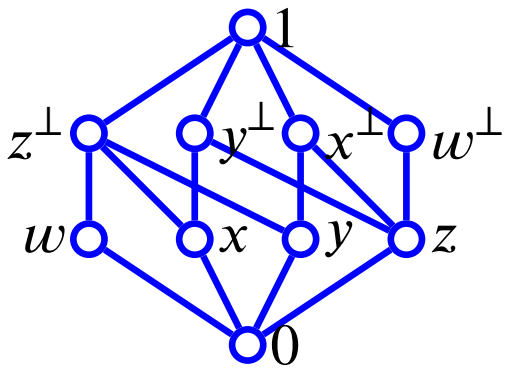

6.

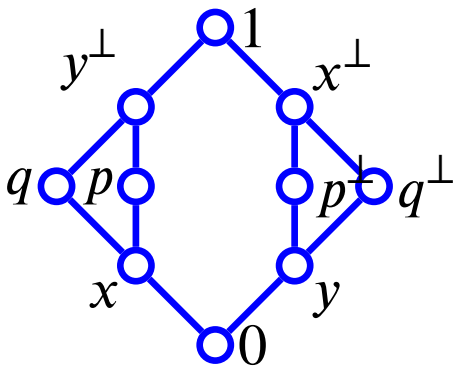

7.

Lattices that are **orthocomplemented** and **orthomodular** but *not modular orthocomplemented* and hence also *non-Boolean*:

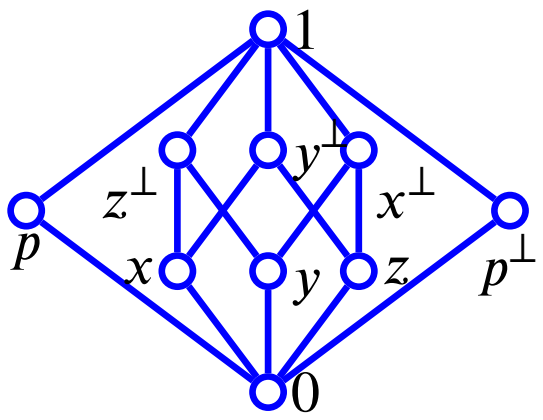

8.

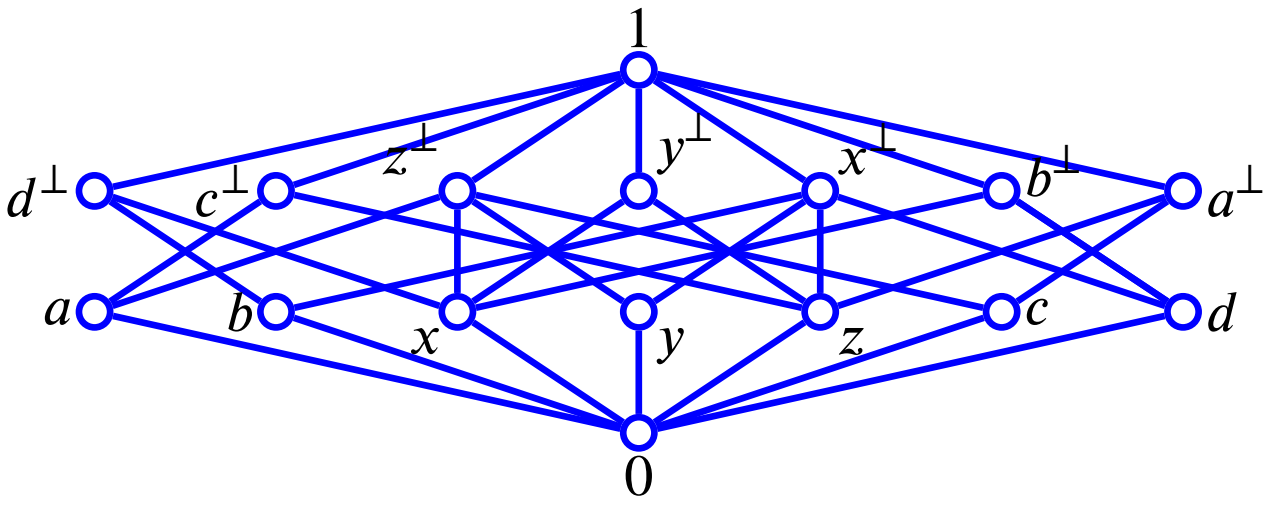

9.

Lattices that are **orthocomplemented**, **orthomodular**, and **modular orthocomplemented** but *non-Boolean*:

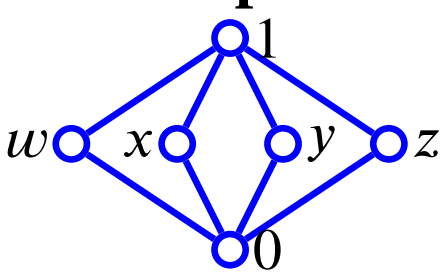

10. $M_4$ *lattice*

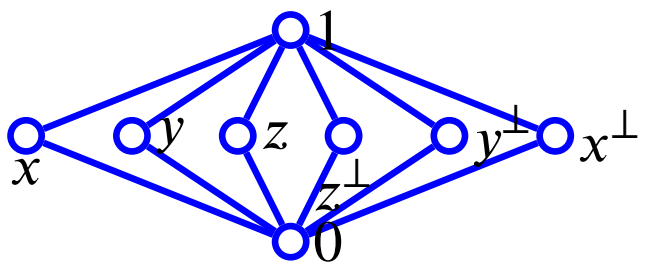

11. $M_6$ *lattice*

Lattices that are **orthocomplemented**, **orthomodular**, **modular orthocomplemented** and **Boolean**:

[56] [12], pages 33–42, [117], page 250, [98], page 24, ⟨Figure 3.2⟩, [157], page 12, [88], page 50

12. $L_1$ *lattice*

13. $L_2$ *lattice*

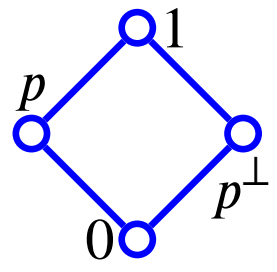

14. $L_2^2$ *lattice*

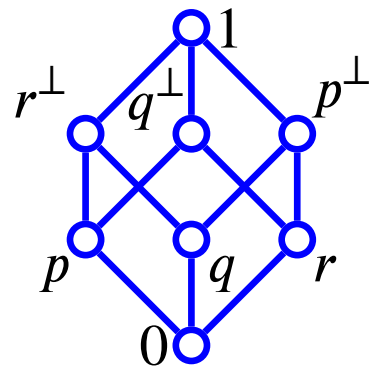

15. $L_2^3$ *lattice*

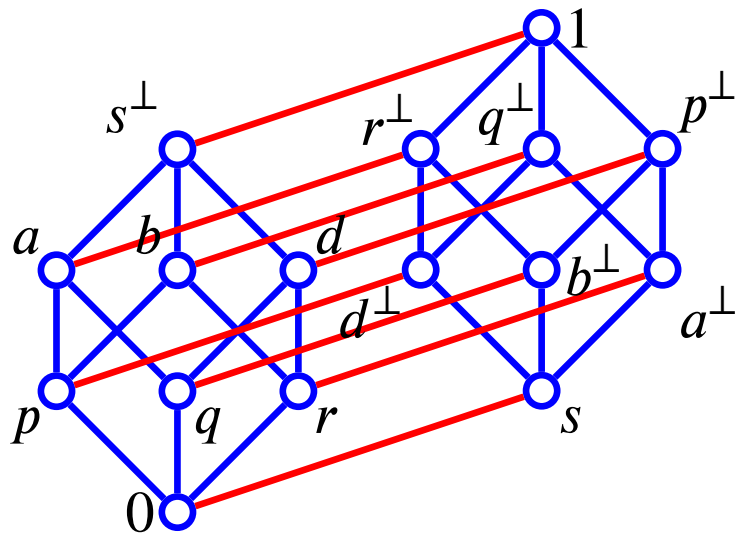

16. $L_2^4$ *lattice*

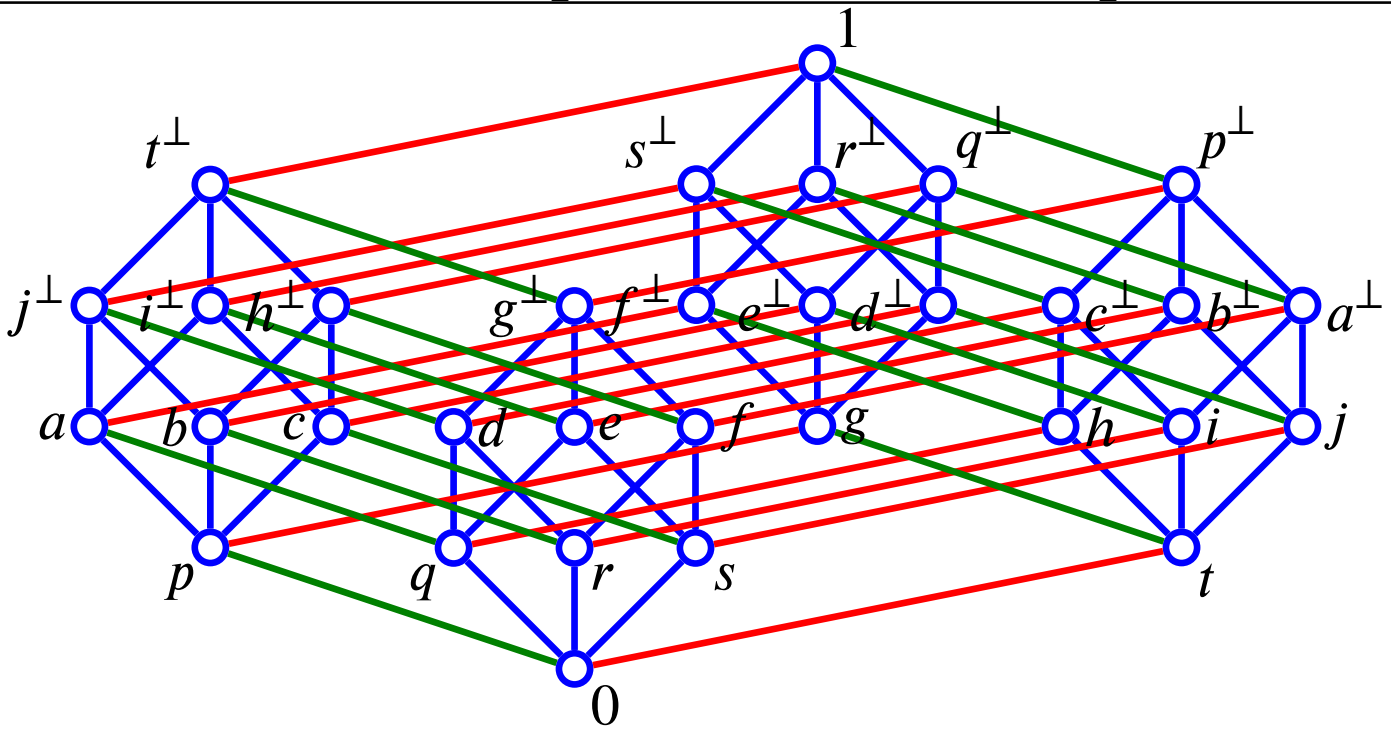

17. $L_2^5$ *lattice*

**Example 1.75** The structure $\left(2^{\mathbb{R}^N}, +, \cap, \varnothing, \boldsymbol{H}\,;\subseteq\right)$ is an **orthocomplemented lattice** where

- $\mathbb{R}^N$ is an **Euclidean space** with dimension $N$
- $2^{\mathbb{R}^N}$ is the set of all subspaces of $\mathbb{R}^N$
- $\boldsymbol{V} + \boldsymbol{W}$ is the *Minkowski sum* of subspaces $\boldsymbol{V}$ and $\boldsymbol{W}$
- $\boldsymbol{V} \cap \boldsymbol{W}$ is the *intersection* of subspaces $\boldsymbol{V}$ and $\boldsymbol{W}$.

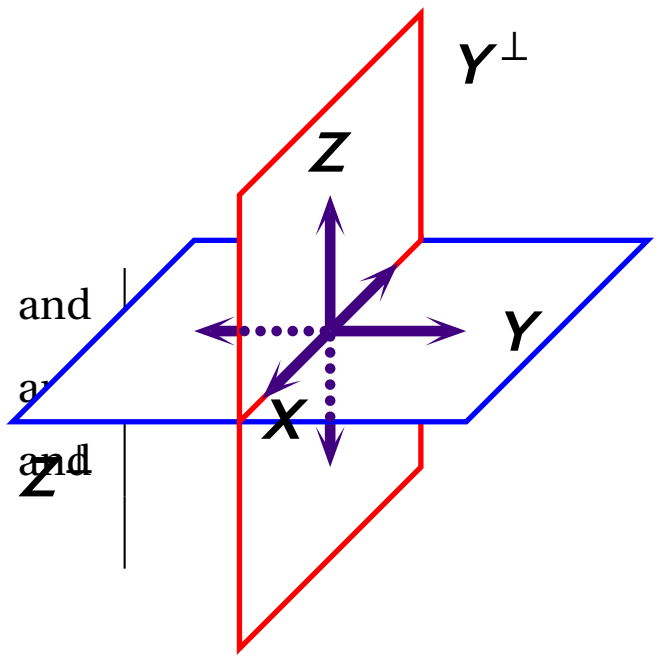

**Example 1.76** The structure $(2^{\boldsymbol{H}}, \oplus, \cap, \varnothing, \boldsymbol{H}\,;\subseteq)$ is an **orthocomplemented lattice** where $\boldsymbol{H}$ is a **Hilbert space**, $2^{\boldsymbol{H}}$ is the set of all closed subspaces of $\boldsymbol{H}$, $\boldsymbol{X} + \boldsymbol{Y}$ is the *Minkowski sum* of subspaces $\boldsymbol{X}$ and $\boldsymbol{Y}$, $\boldsymbol{X} \oplus \boldsymbol{Y} \triangleq (\boldsymbol{X} + \boldsymbol{Y})^-$ is the *closure* of $\boldsymbol{X} + \boldsymbol{Y}$, and $\boldsymbol{X} \cap \boldsymbol{Y}$ is the *intersection* of subspaces $\boldsymbol{X}$ and $\boldsymbol{Y}$.

### 1.3.2 Properties

**Theorem 1.77** [57] *Let $x^\perp$ be the* ORTHOCOMPLEMENT (Definition 1.72 page 20) *of an element $x$ in a* BOUNDED LATTICE $\boldsymbol{L} \triangleq (X, \vee, \wedge, 0, 1\,; \leq)$.

[57] [12], pages 30–31, [20], page 830, ⟨L74⟩, [29], page 37, ⟨3B.13. Theorem⟩

$$\left.\begin{array}{l}\textit{L is}\\ \textit{ortho-}\\ \textit{complemented}\end{array}\right\} \implies \left\{\begin{array}{llcllll}(1). & 0^\perp &=& 1 & & (\text{BOUNDARY CONDITION}) & \text{and}\\ (2). & 1^\perp &=& 0 & & (\text{BOUNDARY CONDITION}) & \text{and}\\ (3). & (x \vee y)^\perp &=& x^\perp \wedge y^\perp & \forall x,y\in X & (\text{DISJUNCTIVE DE MORGAN}) & \text{and}\\ (4). & (x \wedge y)^\perp &=& x^\perp \vee y^\perp & \forall x,y\in X & (\text{CONJUNCTIVE DE MORGAN}) & \text{and}\\ (5). & x \vee x^\perp &=& 1 & \forall x\in X & (\text{EXCLUDED MIDDLE}). & \end{array}\right.$$

✎PROOF: Let $x^\perp \triangleq \neg x$, where $\neg$ is an *ortho negation* function (Definition 2.14 page 29). Then this theorem follows directly from Theorem 2.21 (page 30). ✏

**Corollary 1.78** *Let $\mathbf{L} \triangleq (X, \vee, \wedge, 0, 1; \leq)$ be a* LATTICE *(Definition 1.31 page 10).*

$$\left\{\begin{array}{l}\textbf{\textit{L is orthocomplemented}}\\ \text{(Definition 1.72 page 20)}\end{array}\right\} \implies \left\{\begin{array}{l}\textbf{\textit{L is complemented}}\\ \text{(Definition 1.63 page 17)}\end{array}\right\}$$

✎PROOF: This follows directly from the definition of *orthocomplemented lattices* (Definition 1.72 page 20) and *complemented lattices* (Definition 1.63 page 17). ✏

**Example 1.79**

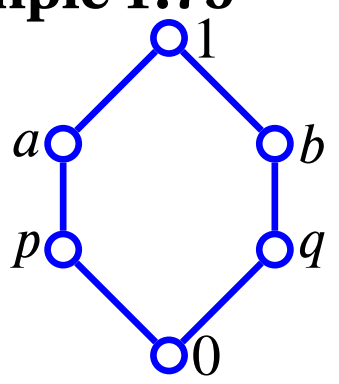

The $O_6$ *lattice* (Definition 1.73 page 20) illustrated to the left is both **orthocomplemented** (Definition 1.72 page 20) and **multiply complemented** (Definition 1.63 page 17). The lattice illustrated to the right is **multiply complemented**, but is **non-orthocomplemented**.

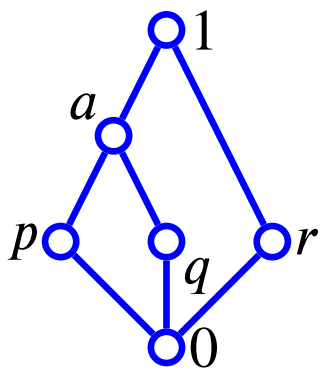

✎PROOF:

(1) Proof that $O_6$ *lattice* is multiply complemented: $b$ and $q$ are both *complements* of $p$.

(2) Proof that the right side lattice is multiply complemented: $a$, $p$, and $q$ are all *complements* of $r$.

✏

### 1.3.3 Restrictions resulting in Boolean algebras

**Proposition 1.80** [58] *Let $\mathbf{L} = (X, \vee, \wedge, 0, 1; \leq)$ be a* BOUNDED LATTICE *(Definition 1.39 page 12).*

$$\left\{\begin{array}{lll}1. & \textbf{\textit{L is orthocomplemented}} & \text{(Definition 1.72 page 20)} \quad \text{and}\\ 2. & \textbf{\textit{L is distributive}} & \text{(Definition 1.53 page 15)}\end{array}\right\} \implies \left\{\begin{array}{l}\textbf{\textit{L is Boolean}}\\ \text{(Definition 1.69 page 18)}\end{array}\right\}$$

[58] [98], page 22

✎PROOF:

$$\left\{\begin{array}{ll} \mathbf{L} \text{ is } \textit{orthocomplemented} & \text{and} \\ \mathbf{L} \text{ is } \textit{distributive} & \end{array}\right\} \implies \left\{\begin{array}{ll} \mathbf{L} \text{ is } \textit{complemented} & \text{and} \\ \mathbf{L} \text{ is } \textit{distributive} & \end{array}\right\} \quad \text{by Corollary 1.78}$$

$$\implies \left\{ \mathbf{L} \text{ is } \textbf{Boolean} \right\} \quad \text{by Definition 1.69}$$

✏

The *center* of an *orthocomplemented lattice* is defined later, but here is a characterization involving it now anyways.

**Proposition 1.81** *Let* $\mathbf{L} = (X, \vee, \wedge, 0, 1; \leq)$ *be a* LATTICE *(Definition 1.31 page 10).*

$$\left\{\begin{array}{lll} 1. & \mathbf{L} \textit{ is } \textbf{\textit{orthocomplemented}} & \text{(Definition 1.72 page 20)} \quad \textit{and} \\ 2. & \textit{Every } x \in \mathbf{L} \textit{ is in the } \textbf{\textit{center}} \textit{ of } \mathbf{L} & \text{(Definition 3.15 page 37)} \end{array}\right\} \iff \left\{ \mathbf{L} \textit{ is } \textbf{\textit{Boolean}} \right\}$$

✎PROOF:

(1) Proof that (1,2) $\implies$ *Boolean*: $\mathbf{L}$ is *Boolean* because it satisfies *Huntington's Fourth Set* (Proposition 1.71 page 19), as demonstrated by the following ...

(a) Proof that $x \vee x = x$ (*idempotent*): $\mathbf{L}$ is a *lattice* (by definition of $\mathbf{L}$), and all lattices are *idempotent* (Definition 1.31 page 10).

(b) Proof that $x \vee y = y \vee x$ (*commutative*): $\mathbf{L}$ is a *lattice* (by definition of $\mathbf{L}$), and all lattices are *commutative* (Definition 1.31 page 10).

(c) Proof that $(x \vee y) \vee z = x \vee (y \vee z)$ (*associative*): $\mathbf{L}$ is a *lattice* (by definition of $\mathbf{L}$), and all lattices are *associative* (Definition 1.31 page 10).

(d) Proof that $(x^\perp \vee y^\perp)^\perp \vee (x^\perp \vee y)^\perp = x$ (*Huntington's axiom*):

$$\begin{aligned} & (x^\perp \vee y^\perp)^\perp \vee (x^\perp \vee y)^\perp \\ &= (x^\perp \perp \wedge y^\perp \perp) \vee (x^\perp \perp \wedge y^\perp) && \text{by } \textit{de Morgan} \text{ property (Theorem 1.77 page 22)} \\ &= (x \wedge y) \vee (x \wedge y^\perp) && \text{by } \textit{involution} \text{ property (Definition 1.72 page 20)} \\ &= x && \text{by def. of } \textit{center} \text{ (Definition 3.15 page 37)} \end{aligned}$$

(2) Proof that (1) $\impliedby$ *Boolean*:

(a) Proof that $x \vee x^\perp = 1$: by definition of *Boolean algebras* (Definition 1.69 page 18).

(b) Proof that $x \wedge x^\perp = 0$: by definition of *Boolean algebras* (Definition 1.69 page 18).

(c) Proof that $x^{\perp\perp} = x$: by *involutory* property of *Boolean algebra* (Theorem 1.70 page 19).

(d) Proof that $x \leq y \implies y^\perp \leq x^\perp$:

$$\begin{aligned} y^\perp \leq x^\perp \iff y^\perp &= y^\perp \wedge x^\perp && \text{by Lemma 1.33 page 10} \\ \iff y^{\perp\perp} &= \left(y^\perp \wedge x^\perp\right)^\perp \\ \iff y^{\perp\perp} &= y^{\perp\perp} \vee x^{\perp\perp} && \text{by } \textit{de Morgan} \text{ property (Theorem 1.70 page 19)} \\ \iff y &= y \vee x && \text{by } \textit{involutory} \text{ property (Theorem 1.70 page 19)} \\ \iff y &= y && \text{by } x \leq y \text{ hypothesis} \end{aligned}$$

(3) Proof that (2) $\Longleftarrow$ *Boolean*: for all $x, y \in \boldsymbol{L}$

$$
\begin{aligned}
(x \wedge y) \vee (x \wedge y^{\perp}) &= [(x \wedge y) \vee x] \wedge [(x \wedge y) \vee y^{\perp}] && \text{by distributive property (Theorem 1.70 page 19)} \\
&= x \wedge [(x \wedge y) \vee y^{\perp}] && \text{by absorptive property (Theorem 1.70 page 19)} \\
&= x \wedge [(x \vee y^{\perp}) \wedge (y \vee y^{\perp})] && \text{by distributive property (Theorem 1.70 page 19)} \\
&= x \wedge (x \vee y^{\perp}) \wedge 1 && \text{by complement property (Theorem 1.70 page 19)} \\
&= x && \text{by absorptive property (Theorem 1.70 page 19)} \\
&\Longrightarrow x \copyright y \quad \forall x, y \in \boldsymbol{L} && \text{by Definition 3.9 page 36} \\
&\Longrightarrow x \text{ is in the center of } \boldsymbol{L} && \text{by Definition 3.15 page 37}
\end{aligned}
$$

✏

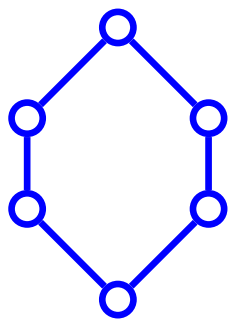

**Example 1.82** The $O_6$ *lattice* (Definition 1.73 page 20) illustrated to the left is **orthocomplemented** (Definition 1.72 page 20) but **non-join-distributive** (Definition 1.53 page 15),and hence *non-Boolean*. The lattice illustrated to the right is **orthocomplemented** *and* **distributive** and hence also **Boolean** (Proposition 1.80 page 23).

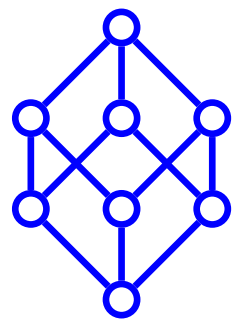

### 1.3.4 Orthomodular lattices

**Definition 1.83** [59] Let $\boldsymbol{L} \triangleq (X, \vee, \wedge, 0, 1; \leq)$ be a *bounded lattice* (Definition 1.39 page 12).
$\boldsymbol{L}$ is **orthomodular** if

1. $\boldsymbol{L}$ is *orthocomplemented* and
2. $x \leq y \implies x \vee (x^{\perp} \wedge y) = y \quad \forall x, y \in X$ (*orthomodular identity*)

**Theorem 1.84** [60] *Let* $\boldsymbol{L} = (X, \vee, \wedge, 0, 1; \leq)$ *be an algebraic structure.*

$$
\left\{ \begin{array}{ll} \boldsymbol{L} \text{ is an } \textbf{orthomodular lattice} & \text{and} \\ \underbrace{(x \wedge y^{\perp})^{\perp} = y \vee (x^{\perp} \wedge y^{\perp})}_{\text{Elkan's law}} & \forall x, y \in X \end{array} \right\} \implies \left\{ \begin{array}{l} \boldsymbol{L} \text{ is a} \\ \textbf{Boolean algebra} \\ \text{(Definition 1.69 page 18)} \end{array} \right\}
$$

**Definition 1.85** Let $\boldsymbol{L} \triangleq (X, \vee, \wedge, 0, 1; \leq)$ be a *bounded lattice* (Definition 1.39 page 12).
$\boldsymbol{L}$ is a **modular orthocomplemeted lattice** if

1. $\boldsymbol{L}$ is **orthocomplemented** (Definition 1.72 page 20) and
2. $\boldsymbol{L}$ is **modular** (Definition 1.47 page 14)

[59] [98], page 22, [110], page 90, [91]
[60] [144], page 72

**Theorem 1.86** [61] *Let* $\boldsymbol{L}$ *be a lattice.*

$$\{\boldsymbol{L} \text{ is BOOLEAN}\} \implies \{\boldsymbol{L} \text{ is MODULAR ORTHOCOMPLEMENTED}\} \quad \text{(Definition 1.85 page 25)}$$
$$\implies \{\boldsymbol{L} \text{ is ORTHOMODULAR}\} \quad \text{(Definition 1.83 page 25)}$$
$$\implies \{\boldsymbol{L} \text{ is ORTHOCOMPLEMENTED}\} \quad \text{(Definition 1.72 page 20)}$$

# 2 Background: functions on lattices

## 2.1 Valuations

**Definition 2.1** [62] Let $\boldsymbol{L} \triangleq (X, \vee, \wedge; \leq)$ be a *lattice* (Definition 1.31 page 10).
A function $\mathrm{v} \in \mathbb{R}^X$ is a **valuation** on $\boldsymbol{L}$ if

$$\mathrm{v}(x \vee y) + \mathrm{v}(x \wedge y) = \mathrm{v}(x) + \mathrm{v}(y) \quad \forall x,y \in X$$

**Proposition 2.2** *Let* $\mathrm{v} \in \mathbb{R}^X$ *be a* FUNCTION *on a* LATTICE $\boldsymbol{L} \triangleq (X, \vee, \wedge; \leq)$ *(Definition 1.31 page 10).*

$$\{\boldsymbol{L} \text{ is LINEAR (Definition 1.31 page 10)}\} \implies \{\mathrm{v} \text{ is a VALUATION (Definition 2.1 page 26)}\}$$

✎PROOF: Let $x, y \in X$ such that $x \leq y$ or $y \lneq x$.

$$\mathrm{v}(x \vee y) + \mathrm{v}(x \wedge y) = \mathrm{v}(x) + \mathrm{v}(y) \qquad \text{because } \boldsymbol{L} \text{ is } \textit{linear}$$

✏

**Example 2.3** [63] Consider the *real valued lattice* $\boldsymbol{L} \triangleq (\mathbb{R}, \max, \min; \leq)$.
The *absolute value* function $|\cdot|$ is a *valuation* on $\boldsymbol{L}$.

✎PROOF: $\boldsymbol{L}$ is *linear* (Definition 1.31 page 10), so v is a *valuation* by Proposition 2.2 (page 26). ✏

**Definition 2.4** [64] Let $X$ be a set and $\mathbb{R}^{\vdash}$ the set of non-negative real numbers.
A function $\mathsf{d} \in \mathbb{R}^{\vdash X \times X}$ is a **metric** on $X$ if

1. $\mathsf{d}(x, y) \geq 0 \quad \forall x,y \in X$ (*non-negative*) and
2. $\mathsf{d}(x, y) = 0 \iff x = y \quad \forall x,y \in X$ (*nondegenerate*) and
3. $\mathsf{d}(x, y) = \mathsf{d}(y, x) \quad \forall x,y \in X$ (*symmetric*) and
4. $\mathsf{d}(x, y) \leq \mathsf{d}(x, z) + \mathsf{d}(z, y) \quad \forall x,y,z \in X$ (*subadditive/ triangle inequality*).[65]

A **metric space** is the pair $(X, \mathsf{d})$. A *metric* is also called a **distance function**.

[61] [98], page 32, ⟨20.⟩, [94], page 57
[62] [93], page 127, [18], page 230, ⟨Definition X.1(V1)⟩, [22], page 58, ⟨Exercise 4.25⟩, [43], page 105, ⟨(8.1.1)⟩, [41], page 143, ⟨§10.3⟩, [42], page 193, ⟨§10.3⟩
[63] [101], page 119, ⟨§5.7⟩
[64] [45], page 28, [31], page 21, [82], page 109, [64], [63], page 30

**Definition 2.5** [66] Let $(X, \mathsf{d})$ be a *metric space* (Definition 2.4 page 26).

| | |
|---|---|
| An **open ball** centered at $x$ with radius $r$ | is the set $\mathsf{B}(x, r) \triangleq \{y \in X \mid \mathsf{d}(x, y) \lneq r\}$. |
| A **closed ball** centered at $x$ with radius $r$ | is the set $\overline{\mathsf{B}}(x, r) \triangleq \{y \in X \mid \mathsf{d}(x, y) \leq r\}$. |
| A **unit ball** centered at $x$ | is the set $\mathsf{B}(x, 1)$. |
| A **closed unit ball** centered at $x$ | is the set $\overline{\mathsf{B}}(x, 1)$. |

**Theorem 2.6** [67] *Let* $\mathrm{v} \in \mathbb{R}^{X}$ *be a function on a* LATTICE $\boldsymbol{L} \triangleq (X, \vee, \wedge; \leq)$ *(Definition 1.31 page 10).*

$$\left.\begin{array}{llll} 1. & \mathrm{v}(x \vee y) + \mathrm{v}(x \wedge y) = \mathrm{v}(x) + \mathrm{v}(y) & \forall x,y \in X & (\text{VALUATION}) \text{ and} \\ 2. & x \leq y \implies \mathrm{v}(x) \leq \mathrm{v}(y) & \forall x,y \in X & (\text{ISOTONE}) \end{array}\right\} \implies \left\{\begin{array}{l} \mathsf{d}(x, y) \triangleq \\ \mathrm{v}(x \vee y) - \mathrm{v}(x \wedge y) \\ \text{is a METRIC on } \boldsymbol{L} \end{array}\right.$$

**Definition 2.7** [68] Let v be a *valuation* (Definition 2.1 page 26) on a *lattice* $\boldsymbol{L} \triangleq (X, \vee, \wedge; \leq)$ (Definition 1.31 page 10). Let $\mathsf{d}(x, y)$ be the *metric* defined in Theorem 2.6 (page 27).
The pair $(\boldsymbol{L}, \mathsf{d})$ is called a *metric lattice*.

For *finite modular* lattices, the *height* function $\mathsf{h}(x)$ (Definition 1.42 page 12) can serve as the isotone valuation that induces a metric (next proposition).

**Proposition 2.8** [69] *Let* $\mathsf{h}(x)$ *be the* HEIGHT *(Definition 1.42 page 12) of a point* $x$ *in a* BOUNDED LATTICE *(Definition 1.39 page 12)* $\boldsymbol{L} \triangleq (X, \vee, \wedge, 0, 1; \leq)$.

$$\left\{ 1.\ \boldsymbol{L} \text{ is MODULAR and } 2.\ \boldsymbol{L} \text{ is FINITE} \right\}$$

$$\implies \left\{\begin{array}{llll} 1. & \mathsf{h}(x \vee y) + \mathsf{h}(x \wedge y) = \mathsf{h}(x) + \mathsf{h}(y) & \forall x,y \in X & (\text{VALUATION}) \text{ and} \\ 2. & x \lneq y \implies \mathsf{h}(x) \lneq \mathsf{h}(y) & \forall x,y \in X & (\text{POSITIVE}) \end{array}\right\}$$

$$\implies \left\{\begin{array}{llll} 1. & \mathsf{h}(x \vee y) + \mathsf{h}(x \wedge y) = \mathsf{h}(x) + \mathsf{h}(y) & \forall x,y \in X & (\text{VALUATION}) \text{ and} \\ 2. & x \leq y \implies \mathsf{h}(x) \leq \mathsf{h}(y) & \forall x,y \in X & (\text{ISOTONE}) \end{array}\right\}$$

**Theorem 2.9** [70] *Let* v *be a* VALUATION *(Definition 2.1 page 26) on a* LATTICE $\boldsymbol{L} \triangleq (X, \vee, \wedge; \leq)$ *(Definition 1.31 page 10). Let* $\mathsf{d}(x, y)$ *be the* METRIC *defined in Theorem 2.6 (page 27).*

$$\left\{\begin{array}{c} (\boldsymbol{L}, \mathsf{d}) \text{ is a METRIC LATTICE} \\ \text{(Definition 2.7 page 27)} \end{array}\right\} \implies \left\{\begin{array}{c} \boldsymbol{L} \text{ is MODULAR} \\ \text{(Definition 1.47 page 14)} \end{array}\right\}$$

---

65 [55], ⟨Book I Proposition 20⟩
66 [5], page 35
67 [43], page 105, ⟨(8.1.2)⟩, [18], pages 230–231
68 [43], page 105, [18], page 231, ⟨§X.2⟩
69 [18], page 230
70 [18], page 232, ⟨Theorem X.2⟩, [43], pages 105–106, [22], page 58, ⟨Exercise 4.25⟩

**Example 2.10**
The function h on the *Boolean* (and thus also *modular*) lattice $\boldsymbol{L}_2^3$ illustrated to the right is a *valuation* (Definition 2.1 page 26) that is *positive* (and thus also *isotone*, Proposition 2.8 page 27). Therefore

$$\mathsf{d}(x,y) \triangleq \mathsf{h}(x \vee y) - \mathsf{h}(x \wedge y) \qquad \forall x,y \in X$$

is a *metric* (Definition 2.7 page 27) on $\boldsymbol{L}_2^3$. For example,

$$\mathsf{d}(b,q) \triangleq \mathsf{h}(b \vee q) - \mathsf{h}(b \wedge q) = \mathsf{h}(1) - \mathsf{h}(0) = 3 - 0 = 3\,.$$

The *closed unit ball* centered at $b$ (Definition 2.5 page 27) and illustrated with solid dots to the right is

$$\mathsf{B}(b,1) \triangleq \{x \in X \mid \mathsf{d}(b,x) \le 1\} = \{b,p,r,0\}$$

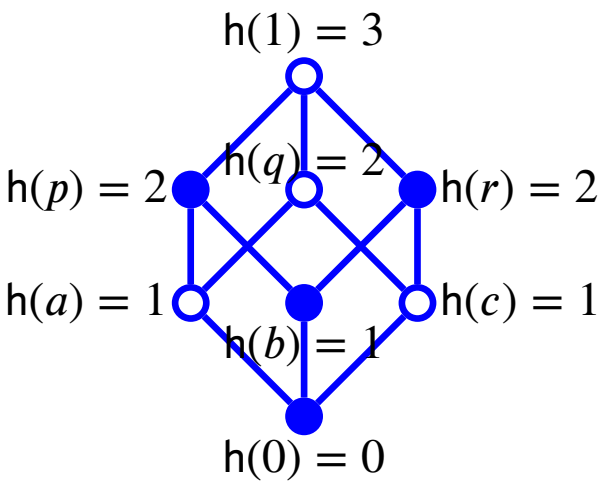

**Example 2.11**
The *height* function h (Definition 1.42 page 12) on the *orthocomplemented* but *non-modular* lattice $O_6$ illustrated to the right is *not* a *valuation* because for example

$$\mathsf{h}(a \vee c) + \mathsf{h}(a \wedge c) = \mathsf{h}(1) + \mathsf{h}(0) = 3 + 0 = 3 \neq 2 = 1 + 1 = \mathsf{h}(a) + \mathsf{h}(b).$$

Moreover, we might expect the "distance" from $a$ to $c$ to be 2. However, if we attempt to use $\mathsf{h}(x)$ to define a metric on $O_6$, then we get

$$\mathsf{d}(a,c) \triangleq \mathsf{h}(a \vee c) - \mathsf{h}(a \wedge c) = \mathsf{h}(1) - \mathsf{h}(0) = 3 - 0 = 3 \neq 2.$$

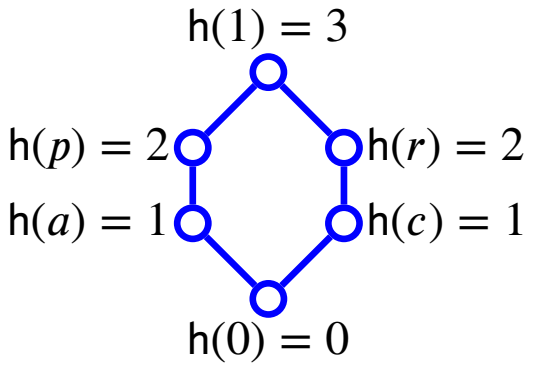

## 2.2 Negation

### 2.2.1 Definitions

**Definition 2.12** [71] Let $\boldsymbol{L} \triangleq (X, \vee, \wedge, 0, 1; \le)$ be a *bounded lattice* (Definition 1.39 page 12).
A *function* $\neg \in X^X$ is a **subminimal negation** on $\boldsymbol{L}$ if [72]

$$x \le y \implies \neg y \le \neg x \quad \forall x,y \in X \quad (antitone).$$

**Definition 2.13** [73] Let $\boldsymbol{L} \triangleq (X, \vee, \wedge, 0, 1; \le)$ be a *bounded lattice* (Definition 1.39 page 12).

[71] [51], pages 4–6, [52], pages 24–26, ⟨2 THE KITE OF NEGATIONS⟩

[72] In the context of natural language, D. Devidi has argued that, *subminimal negation* (Definition 2.12 page 28) is "difficult to take seriously as" a negation. For further details see [40], page 511, [39], page 568, [77], ⟨§2.1.1⟩, [78], ⟨§11.1⟩

[73] [51], pages 4–6, [52], pages 24–26, ⟨2 THE KITE OF NEGATIONS⟩, [161], PAGE 4, ⟨1.6 INTUITIONISM. (B)⟩, [162], PAGE 11, ⟨DEFINITION 16⟩, [70], PAGE 21, ⟨DEFINITION 3.3⟩, [132], PAGE 50, ⟨DEFINITION 2.26⟩, [131], PAGES 98–99, ⟨5.4 NEGATIONS⟩, [10], PAGES 155–156, ⟨(N1) $\neg 0 = 1$ AND $\neg 1 = 0$, (N3) $\neg\neg x = x$⟩

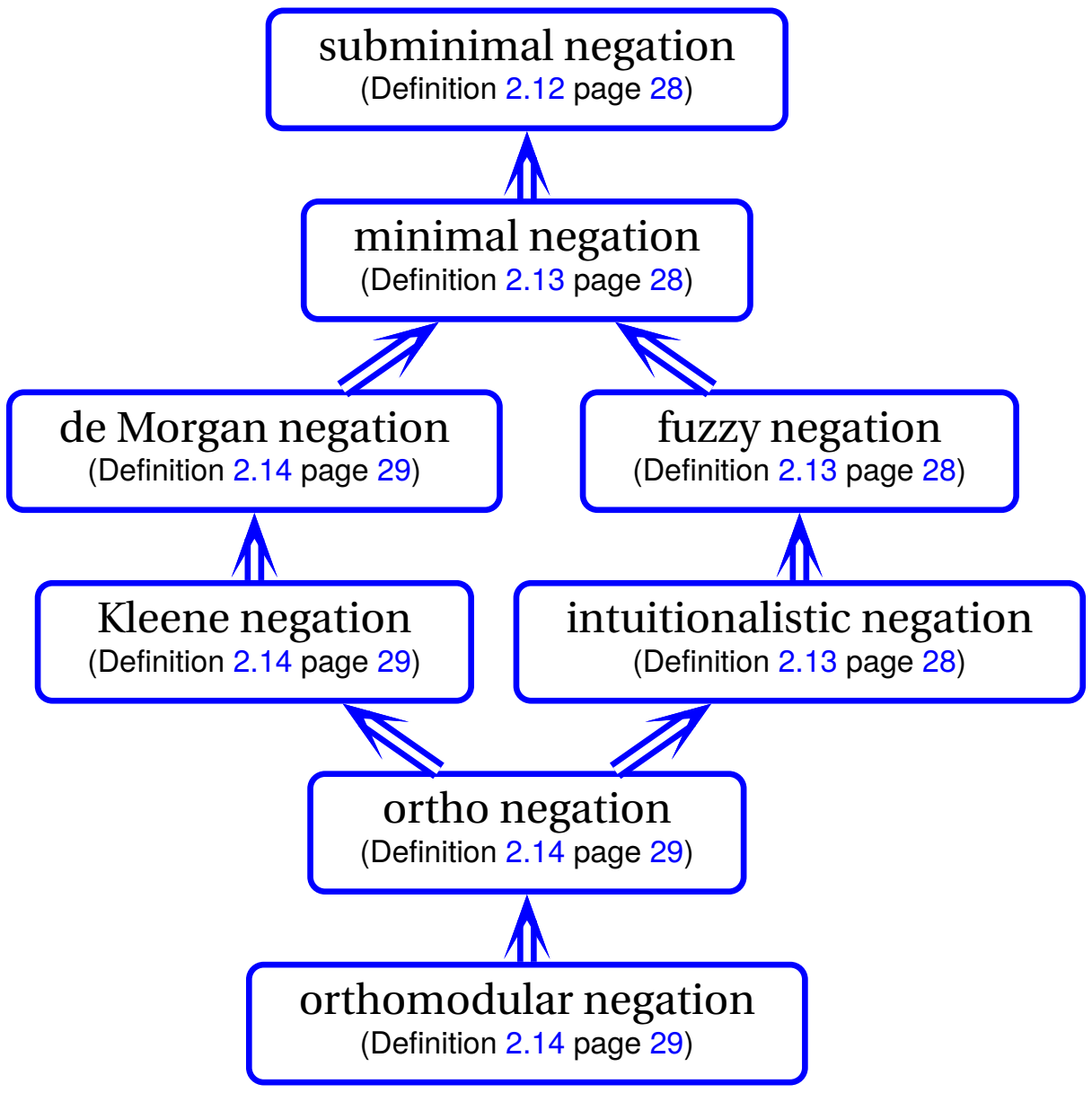

Figure 5: lattice of negations

A *function* $\neg \in X^X$ is a **negation**, or **minimal negation**, on $\boldsymbol{L}$ if

1. $x \le y \implies \neg y \le \neg x \quad \forall x,y \in X$ (*antitone*) and
2. $x \le \neg\neg x \quad \forall x \in X$ (*weak double negation*).

A *minimal negation* $\neg$ is an **intuitionistic negation** on $\boldsymbol{L}$ if

3. $x \wedge \neg x = 0 \quad \forall x,y \in X$ (*non-contradiction*).

A *minimal negation* $\neg$ is a **fuzzy negation** on $\boldsymbol{L}$ if

4. $\neg 1 = 0$ (*boundary condition*).

**Definition 2.14** [74] Let $\boldsymbol{L} \triangleq (X, \vee, \wedge, 0, 1; \le)$ be a *bounded lattice* (Definition 1.39 page 12).

A *minimal negation* $\neg$ is a **de Morgan negation** on $\boldsymbol{L}$ if

5. $x = \neg\neg x \quad \forall x \in X$ (*involutory*).

A *de Morgan negation* $\neg$ is a **Kleene negation** on $\boldsymbol{L}$ if

6. $x \wedge \neg x \le y \vee \neg y \quad \forall x,y \in X$ (*Kleene condition*).

A *de Morgan negation* $\neg$ is an **ortho negation** on $\boldsymbol{L}$ if

7. $x \wedge \neg x = 0 \quad \forall x,y \in X$ (*non-contradiction*).

A *de Morgan negation* $\neg$ is an **orthomodular negation** on $\boldsymbol{L}$ if

8. $x \wedge \neg x = 0 \quad \forall x,y \in X$ (*non-contradiction*) and
9. $x \le y \implies x \vee (x^\perp \wedge y) = y \quad \forall x,y \in X$ (*orthomodular*).

[74] [52], pages 24–26, ⟨2 The Kite of Negations⟩, [96], page 283, [98], page 22, [110], page 90, [91]

**Remark 2.15** [75] The *Kleene condition* is a weakened form of the *non-contradiction* and *excluded middle* properties in the sense $\underbrace{x \wedge \neg x = 0}_{\text{non-contradiction}} \leq \underbrace{1 = y \vee \neg y}_{\text{excluded middle}}$.

**Definition 2.16** Let $\boldsymbol{L} \triangleq (X, \vee, \wedge, \neg, 0, 1; \leq)$ be a *bounded lattice* (Definition 1.39 page 12) with a function $\neg \in X^X$. If $\neg$ is a *negation* (Definition 2.13 page 28), then $\boldsymbol{L}$ is a **lattice with negation**.

### 2.2.2 Properties of negations

**Theorem 2.17** [76] *Let* $\neg \in X^X$ *be a function on a* BOUNDED LATTICE $\boldsymbol{L} \triangleq (X, \vee, \wedge, 0, 1; \leq)$.

$$\left\{\begin{array}{l}\neg \textit{ is a}\\ \text{FUZZY NEGATION}\end{array}\right\} \implies \left\{ \neg 0 = 1 \quad \text{(BOUNDARY CONDITION)} \right\}$$

**Theorem 2.18** [77] *Let* $\neg \in X^X$ *be a function on a* BOUNDED LATTICE $\boldsymbol{L} \triangleq (X, \vee, \wedge, 0, 1; \leq)$.

$$\left\{\begin{array}{l}\neg \textit{ is an}\\ \text{INTUITIONISTIC NEGATION}\end{array}\right\} \implies \left\{\begin{array}{llcll}(a) & \neg 1 & = & 0 & \text{(BOUNDARY CONDITION)} \quad \textit{and}\\ (b) & \neg 0 & = & 1 & \text{(BOUNDARY CONDITION)} \quad \textit{and}\\ (c) & \multicolumn{4}{l}{\neg \textit{ is a } \text{FUZZY NEGATION}}\end{array}\right\}$$

**Theorem 2.19** [78] *Let* $\neg \in X^X$ *be a function on a* BOUNDED LATTICE $\boldsymbol{L} \triangleq (X, \vee, \wedge, 0, 1; \leq)$.

$$\left\{\begin{array}{l}\neg \textit{ is a}\\ \textit{minimal}\\ \textit{negation}\end{array}\right\} \implies \left\{\begin{array}{lclll}\neg x \vee \neg y & \leq & \neg(x \wedge y) & \forall x,y \in X & \text{(CONJUNCTIVE DE MORGAN INEQUALITY)} \quad \textit{and}\\ \neg(x \vee y) & \leq & \neg x \wedge \neg y & \forall x,y \in X & \text{(DISJUNCTIVE DE MORGAN INEQUALITY)}\end{array}\right\}$$

**Theorem 2.20** [79] *Let* $\neg \in X^X$ *be a function on a* BOUNDED LATTICE $\boldsymbol{L} \triangleq (X, \vee, \wedge, 0, 1; \leq)$.

$$\left.\begin{array}{l}\neg \textit{ is a}\\ \textit{de Morgan negation}\end{array}\right\} \implies \left\{\begin{array}{lclll}\neg(x \vee y) & = & \neg x \wedge \neg y & \forall x,y \in X & \text{(DISJUNCTIVE DE MORGAN)} \quad \textit{and}\\ \neg(x \wedge y) & = & \neg x \vee \neg y & \forall x,y \in X & \text{(CONJUNCTIVE DE MORGAN)}\end{array}\right.$$

**Theorem 2.21** [80] *Let* $\neg \in X^X$ *be a function on a* BOUNDED LATTICE $\boldsymbol{L} \triangleq (X, \vee, \wedge, 0, 1; \leq)$.

$$\left\{\begin{array}{l}\neg \textit{ is an}\\ \textit{ortho}\\ \textit{negation}\end{array}\right\} \implies \left\{\begin{array}{lrclll}1. & \neg 0 & = & 1 & & \text{(BOUNDARY CONDITION)} \quad \textit{and}\\ 2. & \neg 1 & = & 0 & & \text{(BOUNDARY CONDITION)} \quad \textit{and}\\ 3. & \neg(x \vee y) & = & \neg x \wedge \neg y & \forall x,y \in X & \text{(DISJUNCTIVE DE MORGAN)} \quad \textit{and}\\ 4. & \neg(x \wedge y) & = & \neg x \vee \neg y & \forall x,y \in X & \text{(CONJUNCTIVE DE MORGAN)} \quad \textit{and}\\ 5. & x \vee \neg x & = & 1 & \forall x \in X & \text{(EXCLUDED MIDDLE)} \quad \textit{and}\\ 6. & x \wedge \neg x & \leq & y \vee \neg y & \forall x,y \in X & \text{(KLEENE CONDITION).}\end{array}\right\}$$

---

[75] [26], page 78
[76] [77], ⟨§2.1.2⟩, [78], ⟨§11.2⟩
[77] [77], ⟨§2.1.2⟩, [78], ⟨§11.2⟩
[78] [77], ⟨§2.1.2⟩, [78], ⟨§11.2⟩
[79] [77], ⟨§2.1.2⟩, [78], ⟨§11.2⟩
[80] [77], ⟨§2.1.2⟩, [78], ⟨§11.2⟩

## 2.3 Projections

**Definition 2.22** [81] Let $\mathbf{L} \triangleq (X, \vee, \wedge, 0, 1; \leq)$ be an *orthocomplemented lattice* (Definition 1.72 page 20). A function $\phi_x \in X^X$ is a **Sasaki projection** on $x \in X$ if $\phi_x(y) \triangleq (y \vee x^\perp) \wedge x$.
The *Sasaki projections* $\phi_x$ and $\phi_y$ are **permutable** if $\phi_x \circ \phi_y(u) = \phi_y \circ \phi_x(u) \quad \forall u \in X$.

**Proposition 2.23** *Let $\phi_x(y)$ be the* SASAKI PROJECTION OF $y$ ONTO $x$ *(Definition 2.22 page 31) in an* ORTHOCOMPLEMENTED LATTICE $\mathbf{L} \triangleq (X, \vee, \wedge, 0, 1; \leq)$.

*(1).* $x \leq y \implies \phi_x(y) = x \quad \forall x,y \in X$
*(2).* $y \leq x \implies y \leq \phi_x(y) \leq x \quad \forall x,y \in X$
*(3).* $y \leq x$ *and* $\mathbf{L}$ *is* BOOLEAN $\implies \phi_x(y) = y \quad \forall x,y \in X$

✎PROOF:

| | |
|---|---|
| (1) $\implies \phi_x(y) \triangleq (y \vee x^\perp) \wedge x$ | by definition of *Sasaki projection* (Definition 2.22 page 31) |
| $= 1 \wedge x$ | by $x \leq y$ hypothesis and Proposition 3.1 page 34 |
| $= x$ | by property of bounded lattices (Proposition 1.41 page 12) |
| (2) $\implies \boxed{y} = y \wedge x$ | by $y \leq x$ hypothesis |
| $\leq (y \vee x^\perp) \wedge x$ | by definition of $\vee$ (Definition 1.27 page 9) |
| $= \boxed{\phi_x(y)}$ | by definition of *Sasaki projection* (Definition 2.22 page 31) |
| $\leq (y \vee x^\perp) \wedge x$ | by definition of *Sasaki projection* (Definition 2.22 page 31) |
| $\leq \boxed{x}$ | by definition of $\wedge$ (Definition 1.28 page 9) |
| (3) $\implies \phi_x(y) = (y \vee x^\perp) \wedge x$ | by definition of *Sasaki projection* (Definition 2.22 page 31) |
| $= (y \wedge x) \vee (x^\perp \wedge x)$ | by *distributive* property of *Boolean lattices* (Theorem 1.70 page 19) |
| $= (y \wedge x) \vee 0$ | by *non-contradiction* of *Boolean lattices* (Theorem 1.70 page 19) |
| $= (y \wedge x)$ | by *boundary* property of *bounded lattices* (Proposition 1.41 page 12) |
| $= y$ | by $y \leq x$ hypothesis and definition of $\wedge$ (Definition 1.28 page 9) |

✏

**Proposition 2.24** *Let $\phi_x(y)$ be the* SASAKI PROJECTION OF $y$ ONTO $x$ *(Definition 2.22 page 31) in an* ORTHOCOMPLEMENTED LATTICE $(X, \vee, \wedge, 0, 1; \leq)$.

*(1).* $\phi_0(y) = 0 \quad \forall y \in X$
*(2).* $\phi_x(0) = 0 \quad \forall x \in X$
*(3).* $\phi_1(y) = 1 \quad \forall y \in X$
*(4).* $\phi_x(1) = x \quad \forall x \in X$
*(5).* $\phi_x(x^\perp) = 0 \quad \forall x \in X$

[81] [127], pages 158–159, ⟨equation (S)⟩, [152], page 300, ⟨Def.5.1, cf Foulis 1962⟩, [98], page 117

✎PROOF:

$$\begin{aligned}
\phi_0(y) &= 0 && \text{because } 0 \le y \text{ and by Proposition 2.23 page 31}\\
\phi_x(0) &\triangleq \left(0 \vee x^\perp\right) \wedge x && \text{by definition of \textit{Sasaki projection} (Definition 2.22 page 31)}\\
&= x^\perp \wedge x && \text{by property of bounded lattices (Proposition 1.41 page 12)}\\
&= 0 && \text{by definition of \textit{orthocomplemented} (Definition 1.72 page 20)}\\
\phi_1(y) &\triangleq \left(y \vee 1^\perp\right) \wedge 1 && \text{by definition of \textit{Sasaki projection} (Definition 2.22 page 31)}\\
&= (y \vee 0) \wedge 1 && \text{by \textit{boundary condition} (Theorem 2.21 page 30)}\\
&= y \wedge 1 && \text{by property of bounded lattices (Proposition 1.41 page 12)}\\
&= 1 && \text{by property of bounded lattices (Proposition 1.41 page 12)}\\
\phi_x(1) &= x && \text{because } x \le 1 \text{ and by Proposition 2.23 page 31}\\
\phi_x\left(x^\perp\right) &\triangleq \left(x^\perp \vee x^\perp\right) \wedge x && \text{by definition of \textit{Sasaki projection} (Definition 2.22 page 31)}\\
&= x^\perp \wedge x && \text{by \textit{idempotency} of lattices (Theorem 1.32 page 10)}\\
&= 0 && \text{by \textit{non-contradiction} prop. of \textit{orthocomplemented lattice} (Definition 1.72 page 20)}
\end{aligned}$$

✎

**Example 2.25** Here are some examples of projections in the $O_6$ *lattice* onto the element $x$:

$$\begin{array}{lcllllllllll}
\phi_p(q) & \triangleq & (q & \vee & p^\perp) & \wedge & p & = & p^\perp & \wedge & p & = & 0 & \text{(because } p \perp q)\\
\phi_p(p^\perp) & \triangleq & (p^\perp & \vee & p^\perp) & \wedge & p & = & p^\perp & \wedge & p & = & 0 & \text{(because } p \perp p^\perp)\\
\phi_p(q^\perp) & \triangleq & (q^\perp & \vee & p^\perp) & \wedge & p & = & 1 & \wedge & p & = & p & \text{(because } p \le q^\perp)\\
\phi_{q^\perp}(p) & \triangleq & (p & \vee & q) & \wedge & q^\perp & = & 1 & \wedge & q^\perp & = & q^\perp & \text{(because } q^\perp \le 1)\\
\phi_p(1) & \triangleq & (1 & \vee & p^\perp) & \wedge & p & = & 1 & \wedge & p & = & p & \text{(because } p \le 1)\\
\phi_p(0) & \triangleq & (0 & \vee & p^\perp) & \wedge & p & = & p^\perp & \wedge & p & = & 0 & \text{(because } p \perp 0)
\end{array}$$

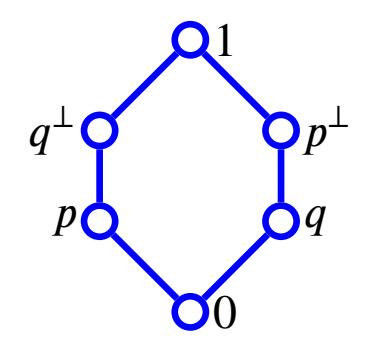

**Example 2.26**

Let $\mathbb{R}^3$ be the *3-dimensional Euclidean space* (Example 1.75 page 22) with subspaces $\boldsymbol{Z}$ and $\boldsymbol{V}$. Then the projection operator $P_{\boldsymbol{Z}^\perp}$ onto $\boldsymbol{Z}^\perp$ is a *sasaki projection* $\phi_{\boldsymbol{Z}^\perp}$. In particular

$$\begin{aligned}
\mathbf{P}_{\boldsymbol{Z}^\perp}\boldsymbol{V} &\triangleq \phi_{\boldsymbol{Z}^\perp}(\boldsymbol{V})\\
&\triangleq \left(\boldsymbol{V} + \boldsymbol{Z}^{\perp\perp}\right) \cap \boldsymbol{Z}^\perp\\
&= (\boldsymbol{V} + \boldsymbol{Z}) \cap \boldsymbol{Z}^\perp
\end{aligned}$$

as illustrated to the right.

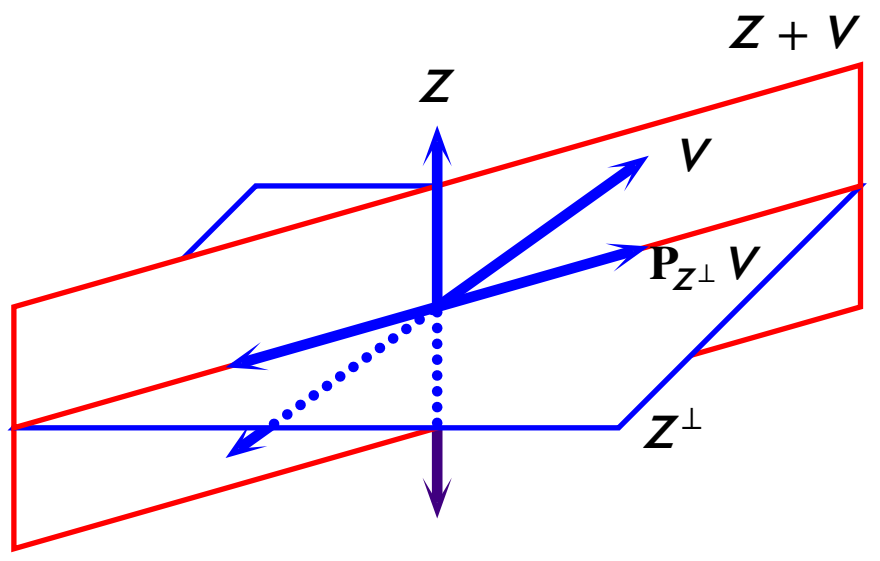

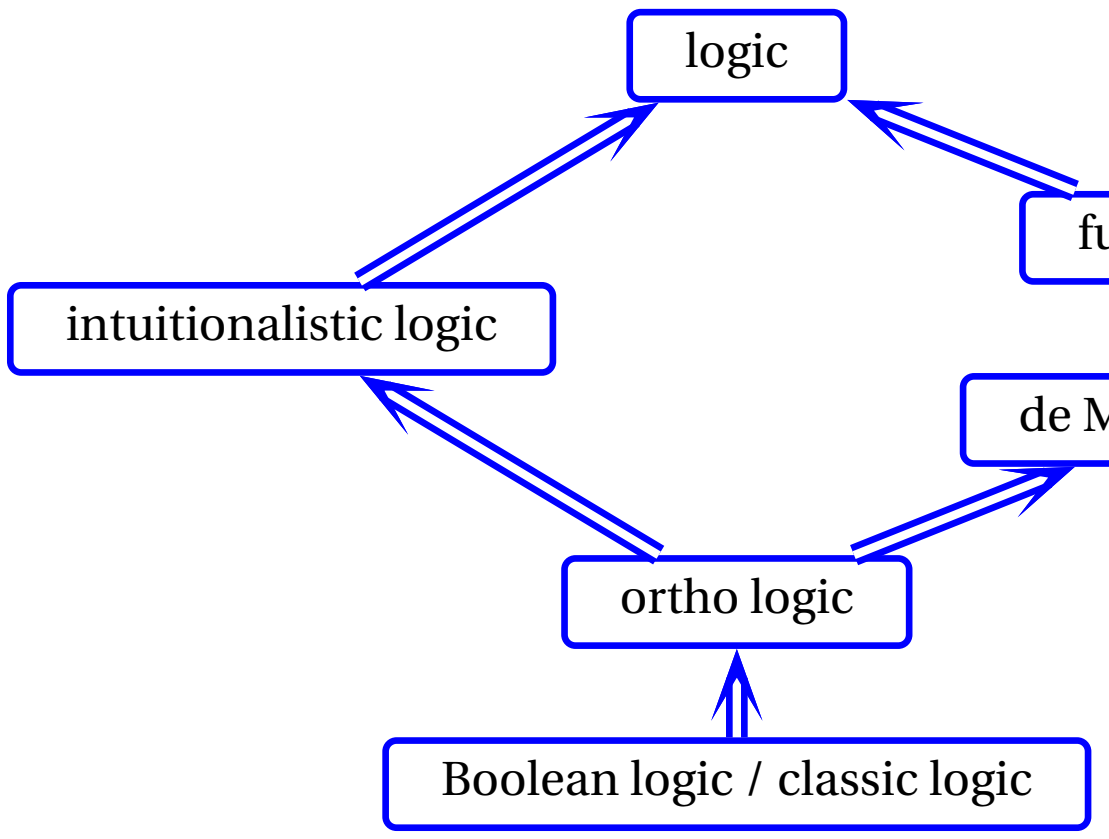

Figure 6: lattice of logics

## 2.4 Logics

**Definition 2.27** [82] Let $\rightarrow$ be an *implication* function defined on a *lattice with negation* $\boldsymbol{L} \triangleq (X, \vee, \wedge, \neg, 0, 1 ; \leq)$ (Definition 2.16 page 30).

| | |
|---|---|
| $(X, \vee, \wedge, \neg, 0, 1 ; \leq, \rightarrow)$ is a **logic** | if $\neg$ is a *minimal negation.* |
| $(X, \vee, \wedge, \neg, 0, 1 ; \leq, \rightarrow)$ is a **fuzzy logic** | if $\neg$ is a *fuzzy negation.* |
| $(X, \vee, \wedge, \neg, 0, 1 ; \leq, \rightarrow)$ is an **intuitionalistic logic** | if $\neg$ is an *intuitionalistic negation.* |
| $(X, \vee, \wedge, \neg, 0, 1 ; \leq, \rightarrow)$ is a **de Morgan logic** | if $\neg$ is a *de Morgan negation.* |
| $(X, \vee, \wedge, \neg, 0, 1 ; \leq, \rightarrow)$ is a **Kleene logic** | if $\neg$ is a *Kleene negation.* |
| $(X, \vee, \wedge, \neg, 0, 1 ; \leq, \rightarrow)$ is an **ortho logic** | if $\neg$ is an *ortho negation.* |
| $(X, \vee, \wedge, \neg, 0, 1 ; \leq, \rightarrow)$ is a **Boolean logic** | if $\neg$ is an *ortho negation* and $\boldsymbol{L}$ is *Boolean.* |

For examples and a definition of *implication*, see [77], ⟨§3.1⟩.

# 3 Background: relations on lattices

The relations in this section are typically defined on an *orthocomplemented lattice* (Definition 1.72 page 20). Here, some relations are generalized to a *lattice with negation* (Definition 2.16 page 30). A *lattice* (Definition 1.31 page 10) with an *ortho negation* successfully defined on it is an *orthocomplemented lattice* (Definition 1.72 page 20). In many cases, these relations only work

[82] [159], page 136, ⟨Definition 2.1⟩, [162], page 11, ⟨Definition 16⟩, [77], ⟨§3.1⟩

well on an *orthocomplemented lattice*, and thus many results are restricted to orthocomplemented lattices.

## 3.1 Orthogonality

**Proposition 3.1** *Let* $(X, \vee, \wedge, 0, 1; \leq)$ *be an* ORTHOCOMPLEMENTED LATTICE *(Definition 1.72 page 20).*

$$x \leq y \implies \left\{\begin{array}{ccccc} x^\perp & \vee & y & = & 1 \quad \textit{and} \\ x & \wedge & y^\perp & = & 0 \end{array}\right\} \qquad \forall x,y \in X$$

✎PROOF:

$$\begin{array}{rlll} x \leq y & \implies & x \vee x^\perp \leq y \vee x^\perp & \text{by \textit{monotone} property of \textit{lattices} (Proposition 1.34 page 11)} \\ & \implies & 1 \leq y \vee x^\perp & \text{by \textit{excluded middle} property (Definition 1.72 page 20)} \\ & \implies & x^\perp \vee y = 1 & \text{by \textit{upper bounded} property of \textit{bounded lattices} (Definition 1.39 page 12)} \\ x \leq y & \implies & x \wedge y^\perp \leq y \wedge y^\perp & \text{by \textit{monotone} property of \textit{lattices} (Proposition 1.34 page 11)} \\ & \implies & x \wedge y^\perp \leq 0 & \text{by \textit{non-contradiction} property (Definition 1.72 page 20)} \\ & \implies & x \wedge y^\perp = 0 & \text{by \textit{lower bounded} property of \textit{bounded lattices} (Definition 1.39 page 12)} \end{array}$$

✏

**Definition 3.2** [83] Let $(X, \vee, \wedge, \neg, 0, 1; \leq)$ be a *lattice with negation* (Definition 2.16 page 30).
The **orthogonality** relation $\perp \in 2^{XX}$ is defined as

$$x \perp y \overset{\text{def}}{\iff} x \leq \neg y$$

If $x \perp y$, we say that $x$ is **orthogonal** to $y$.

**Lemma 3.3** *Let* $(X, \vee, \wedge, \neg, 0, 1; \leq)$ *be a* LATTICE WITH NEGATION *(Definition 2.16 page 30).*

$$\left\{ x \perp y \quad (\text{ORTHOGONAL Definition 3.2 page 34}) \right\} \implies \left\{ y \perp x \quad (\text{SYMMETRIC}) \right\}$$

✎PROOF:

$$\begin{array}{rlll} x \perp y & \implies & x \leq \neg y & \text{by definition of } \perp \text{ (Definition 3.2 page 34)} \\ & \implies & (\neg\neg y) \leq \neg x & \text{by \textit{antitone} property (Definition 1.72 page 20)} \\ & \implies & y \leq \neg x & \text{by \textit{weak double negation} property of \textit{negation} (Definition 2.13 page 28)} \\ & \implies & y \perp x & \text{by definition of } \perp \text{ (Definition 3.2 page 34)} \end{array}$$

✏

83 [157], page 12, [112], page 3

**Lemma 3.4** [84] *Let* $(X, \vee, \wedge, 0, 1; \leq)$ *be an* ORTHOCOMPLEMENTED LATTICE *(Definition 1.72 page 20)*.

$$\underbrace{x \perp y}_{\text{ORTHOGONAL (Definition 3.2 page 34)}} \implies \left\{ \begin{array}{llcll} 1. & x \wedge y & = & 0 & \text{and} \\ 2. & x^{\perp} \vee y^{\perp} & = & 1 & \end{array} \right\}$$

**Remark 3.5** In an *orthocomplemented lattice* $\boldsymbol{L}$, the *orthogonality* relation $\perp$ is in general *non-associative*. That is,

$$\left\{ \begin{array}{ccll} x & \perp & y & \text{and} \\ y & \perp & z & \end{array} \right\} \not\Longrightarrow x \perp z$$

✎PROOF: Consider the $\boldsymbol{L}_2^4$ *Boolean lattice* in Example 1.74 (page 21).

- $a^{\perp} \perp p$ because $a^{\perp} \leq p^{\perp}$.
- $p \perp r$ because $p \leq r^{\perp}$.
- But yet $a^{\perp}$ is *not* orthogonal to $r$ because $a^{\perp} \not\leq r^{\perp}$. ✎

**Example 3.6** In the $O_6$ *lattice* (Definition 1.73 page 20), there are a total of $\binom{6}{2} = \frac{6!}{(6-2)!2!} = \frac{6\times 5}{2} = 15$ distinct unordered (the $\perp$ relation is *symmetric* by Lemma 3.3 page 34 so the order doesn't matter) pairs of elements.

Of these 15 pairs, 8 are orthogonal to each other, and 0 is orthogonal to itself, making a total of 9 orthogonal pairs:

| | | |
|---|---|---|
| $x \perp y$ | $x \perp 0$ | $y^{\perp} \perp 0$ |
| $x \perp x^{\perp}$ | $y \perp 0$ | $1 \perp 0$ |
| $y \perp y^{\perp}$ | $x^{\perp} \perp 0$ | $0 \perp 0$ |

**Example 3.7** In lattice 5 of Example 1.74 (page 21), there are a total of $\binom{10}{2} = \frac{10!}{(10-2)!2!} = \frac{10\times 9}{2} = 45$ distinct unordered pairs of elements.

Of these 45 pairs, 18 are orthogonal to each other, and 0 is orthogonal to itself, making a total of 19 orthogonal pairs:

| | | | |
|---|---|---|---|
| $p \perp p^{\perp}$ | $x \perp x^{\perp}$ | $y \perp z$ | $x^{\perp} \perp 0$ |
| $p \perp x^{\perp}$ | $x \perp y$ | $y \perp 0$ | $y^{\perp} \perp 0$ |
| $p \perp y$ | $x \perp z$ | $z \perp z^{\perp}$ | $z^{\perp} \perp 0$ |
| $p \perp z$ | $x \perp 0$ | $z \perp 0$ | $0 \perp 0$ |
| $p \perp 0$ | $y \perp y^{\perp}$ | $p^{\perp} \perp 0$ | |

**Example 3.8** In the $\mathbb{R}^3$ **Euclidean space** illustrated in Example 1.75 (page 22),

$$\begin{array}{llll} \boldsymbol{X} \subseteq \boldsymbol{Y}^{\perp} \implies \boldsymbol{X} \perp \boldsymbol{Y} & & \boldsymbol{Y} \subseteq \boldsymbol{X}^{\perp} \implies \boldsymbol{Y} \perp \boldsymbol{X} \\ \boldsymbol{X} \subseteq \boldsymbol{Z}^{\perp} \implies \boldsymbol{X} \perp \boldsymbol{Z} & & \boldsymbol{Y} \subseteq \boldsymbol{Z}^{\perp} \implies \boldsymbol{Y} \perp \boldsymbol{Z} \\ \boldsymbol{X} \wedge \boldsymbol{Y} = \boldsymbol{X} \wedge \boldsymbol{Z} = \boldsymbol{Y} \wedge \boldsymbol{Z} = \boldsymbol{0} \end{array}$$

[84] [87], page 67, [78], ⟨Lemma 13.2⟩

## 3.2 Commutativity

The *commutes* relation is defined next. Motivation for the name "commutes" is provided by Proposition 3.14 (page 36) which shows that if $x$ commutes with $y$ in a lattice $\boldsymbol{L}$, then $x$ and $y$ commute in the *Sasaki projection* $\phi_x(y)$ on $\boldsymbol{L}$.

**Definition 3.9** [85] Let $\boldsymbol{L} \triangleq (X, \vee, \wedge, \neg, 0, 1; \le)$ be a *lattice with negation* (Definition 2.16 page 30). The **commutes** relation © is defined as

$$x © y \quad \overset{\text{def}}{\iff} \quad x = (x \wedge y) \vee (x \wedge \neg y) \qquad \forall x,y \in X,$$

in which case we say, "$x$ **commutes** with $y$ in $\boldsymbol{L}$".
That is, © is a relation in $2^{XX}$ such that

$$© \triangleq \left\{(x,y) \in X^2 \,|\, x = (x \wedge y) \vee (x \wedge \neg y)\right\}$$

**Proposition 3.10** [86] *Let* $\boldsymbol{L} \triangleq (X, \vee, \wedge, 0, 1; \le)$ *be an* ORTHOCOMPLEMENTED LATTICE.

| | | | | | | | |
|---|---|---|---|---|---|---|---|
| $x © 0$ | *and* $0 © x$ | $\forall x \in X$ | | $x © y$ | $\iff$ | $x © y^\perp$ | $\forall x,y \in X$ |
| $x © 1$ | *and* $1 © x$ | $\forall x \in X$ | | $x \le y$ | $\implies$ | $x © y$ | $\forall x,y \in X$ |
| $x © x$ | | $\forall x \in X$ | | $x \perp y$ | $\implies$ | $x © y$ | $\forall x,y \in X$ |

**Definition 3.11** Let © be the *commutes* relation (Definition 3.9 page 36) on a *lattice with negation* $\boldsymbol{L} \triangleq (X, \vee, \wedge, \neg, 0, 1; \le)$ (Definition 2.16 page 30). $\boldsymbol{L}$ is **symmetric** if

$$x © y \implies y © x \qquad \forall x,y \in X$$

In general, the commutes relation is not *symmetric*. But Proposition 3.12 (next) describes some conditions under which it *is* symmetric.

**Proposition 3.12** [87] *Let* $(X, \vee, \wedge, 0, 1; \le)$ *be an* ORTHOCOMPLEMENTED LATTICE *(Definition 1.72 page 20).*

$$\underbrace{\{x © y \implies y © x\}}_{© \text{ is SYMMETRIC at } (x,y) \;(1)} \iff \left\{x \le y \implies y = x \vee \left(x^\perp \wedge y\right)\right\} \quad \text{(ORTHOMODULAR IDENTITY)} \quad (2)$$

$$\iff \left\{x \le y \implies x = y \wedge \left(x \vee y^\perp\right)\right\} \quad (x = \phi_y(x) \text{ (SASAKI PROJECTION)}) \quad (3)$$

$$\iff \left\{y = (x \wedge y) \vee \left[y \wedge (x \wedge y)^\perp\right]\right\} \quad (4)$$

$$\iff \left\{x = (x \vee y) \wedge \left[x \vee (x \vee y)^\perp\right]\right\} \quad (5)$$

**Theorem 3.13** [88] *Let* $\boldsymbol{L} \triangleq (X, \vee, \wedge, 0, 1; \le)$ *be an* ORTHOCOMPLEMENTED LATTICE *(Definition 1.72 page 20).*

$$\{x © c \quad \forall x \in X\} \iff \left\{\boldsymbol{L} \text{ is ISOMORPHIC to } [0, c] \times \left[0, c^\perp\right]\right\}$$
$$\text{with isomorphism } \theta(x) \triangleq \left([0, c], \left[0, c^\perp\right]\right).$$

**Proposition 3.14** [89] *Let* $(X, \vee, \wedge, 0, 1; \le)$ *be an* ORTHOMODULAR *lattice.*

[85] [98], page 20, [88], page 79, ⟨A. Commutativity⟩, [115], page 227, ⟨Hilfssatz (Lemma) XII.1.2⟩, [152], page 301, ⟨Def.5.2, cf Foulis 1962⟩, [15], page 833, ⟨"$a = (a \cap x) \cup (a \cap x')$"⟩
[86] [87], page 67, [78], ⟨Proposition 13.2⟩
[87] [87], page 68, [127], page 158, [78], ⟨Proposition 13.3⟩
[88] [98], page 20, [114]
[89] [62], page 66, [152], ⟨cf Foulis 1962⟩

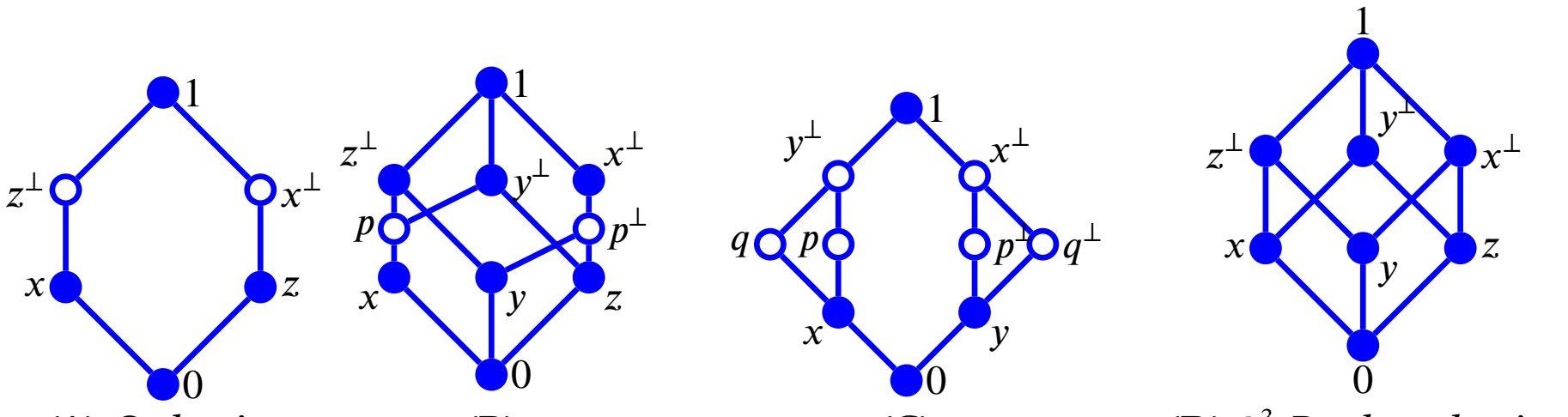

Figure 7: Lattices with centers marked with solid dots (see Example 3.17 page 37)

$$x \copyright y \quad \Longleftrightarrow \quad \phi_x(y) = \phi_y(x) = x \wedge y \qquad \forall x,y \in X$$

## 3.3 Center

An element in an *orthocomplemented lattice* (Definition 1.72 page 20) is in the *center* of the lattice if that element *commutes* (Definition 3.9 page 36) with every other element in the lattice (next definition). *All* the elements of an *orthocomplemented lattice* are in the *center* if and only if that lattice is *Boolean* (Proposition 1.81 page 24).

**Definition 3.15** [90] Let © be the *commutes* relation (Definition 3.9 page 36) on a *lattice with negation* $\mathbf{L} \triangleq (X, \vee, \wedge, \neg, 0, 1; \le)$ (Definition 2.16 page 30). The **center** of $\mathbf{L}$ is defined as

$$\{x \in X \mid x \copyright y \quad \forall y \in X\}$$

**Proposition 3.16** *Let* $\mathbf{L} \triangleq (X, \vee, \wedge, 0, 1; \le)$ *be an* ORTHOCOMPLEMENTED LATTICE *(Definition 1.72 page 20). The elements* 0 *and* 1 *are in the* **center** *of* $\mathbf{L}$.

✎PROOF: This follows directly from Definition 3.9 (page 36) and Proposition 3.10 (page 36). ✏

**Example 3.17** The **center**s of the lattices in Figure 7 (page 37) are illustrated with solid dots. Note that in the case of the Boolean lattice in (D), every dot is in the center (Proposition 1.81 page 24).

## 3.4 D-Posets

**Definition 3.18** [91] Let 1 be the *upper bound* of an *ordered set* $(X, \le)$.
An operation \ is a **difference** on $(X, \le)$ if

[90] [88], page 80
[91] [104], page 22,24, ⟨DEFINITIONS 1,2⟩

1. $x \leq y \implies y \backslash x \leq y \quad \forall x,y \in X$ and
2. $x \leq y \implies y \backslash (y \backslash x) = x \quad \forall x,y \in X$ and
3. $x \leq y \leq z \implies z \backslash y \leq z \backslash x \quad \forall x,y,z \in X$ and
4. $x \leq y \leq z \implies (z \backslash x) \backslash (z \backslash y) = y \backslash x \quad \forall x,y,z \in X$ .

The structure $(X, \leq, \backslash, 1)$ is called a **D-poset**.

**Proposition 3.19** [92] *Let $X$ be a* SET.

$$\left\{ \begin{array}{l} (X, \leq, \backslash, 1) \text{ is a} \\ \text{D-POSET} \\ \text{(Definition 3.18 page 37)} \end{array} \right\} \implies \left\{ \begin{array}{lllll} 1. & x \leq y \leq z & \implies & y \backslash x \leq z \backslash x & \forall x,y,z \in X \text{ and} \\ 2. & x \leq y \leq z & \implies & x \leq z \backslash (y \backslash x) & \forall x,y,z \in X \text{ and} \\ 3. & x \leq y \leq z & \implies & (z \backslash x) \backslash (y \backslash x) = z \backslash y & \forall x,y,z \in X \text{ and} \\ 4. & x \leq y \leq z & \implies & [z \backslash (y \backslash x)] \backslash x = z \backslash y & \forall x,y,z \in X \text{ .} \end{array} \right\}$$

**Example 3.20** [93] The structure $(\mathbb{R}^+, -, \leq)$ is a *D-poset* where $\mathbb{R}^+$ is the set of positive real numbers, $-$ is the standard subtraction operation on $\mathbb{R}$, and $\leq$ is the standard ordering relation on $\mathbb{R}^+$.

**Example 3.21** [94] The structure $(2^X, \backslash, \subseteq)$ is a *D-poset* where $2^X$ is the *power set* of a set $X$, $\backslash$ is the *set difference operator*, and $\subseteq$ is the *set inclusion* relation.

# 4 Background: MRA-wavelet analysis

## 4.1 Transversal Operators

**Definition 4.1** [95]

1. **T** is the **translation operator** on $\mathbb{C}^{\mathbb{C}}$ defined as
$$\mathbf{T}_\tau \mathsf{f}(x) \triangleq \mathsf{f}(x - \tau) \quad \text{and} \quad \mathbf{T} \triangleq \mathbf{T}_1 \qquad \forall \mathsf{f} \in \mathbb{C}^{\mathbb{C}}$$
2. **D** is the **dilation operator** on $\mathbb{C}^{\mathbb{C}}$ defined as
$$\mathbf{D}_\alpha \mathsf{f}(x) \triangleq \mathsf{f}(\alpha x) \quad \text{and} \quad \mathbf{D} \triangleq \sqrt{2}\mathbf{D}_2 \qquad \forall \mathsf{f} \in \mathbb{C}^{\mathbb{C}}$$

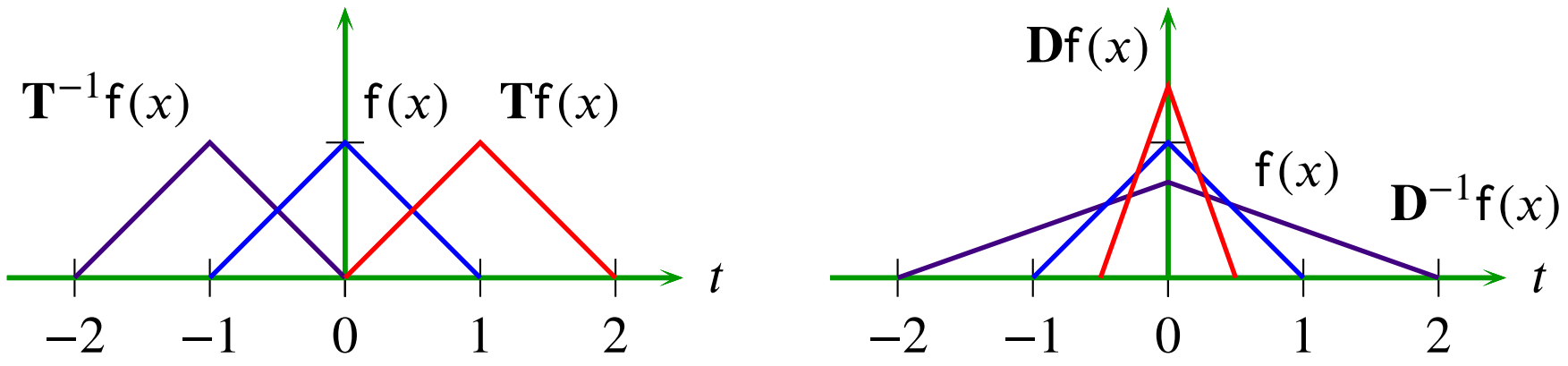

[92] [104], page 23, ⟨PROPOSITION 1.⟩
[93] [104], page 22, ⟨Example 1⟩
[94] [104], page 24, ⟨Example 4⟩
[95] [163], pages 79–80, ⟨Definition 3.39⟩, [27], pages 41–42, [168], page 18, ⟨Definitions 2.3,2.4⟩, [100], page A-21, [8], page 473, [135], page 260, [11], page , [83], page 250, ⟨Notation 9.4⟩, [25], page 74, [69], page 639, [34], page 81, [33], page 2, [75], page 2

**Proposition 4.2** [96] *Let* $\mathbf{T}$ *be the* TRANSLATION OPERATOR *(Definition 4.1 page 38).*

$$\sum_{n\in\mathbb{Z}} \mathbf{T}^n \mathsf{f}(x) = \sum_{n\in\mathbb{Z}} \mathbf{T}^n \mathsf{f}(x+1) \qquad \forall \mathsf{f}\in\mathbb{R}^\mathbb{R} \qquad \left(\sum_{n\in\mathbb{Z}} \mathbf{T}^n\mathsf{f}(x) \text{ is PERIODIC with period } 1\right)$$

**Proposition 4.3** [97] *Let* $\mathbf{T}$ *and* $\mathbf{D}$ *be as defined in Definition 4.1 page 38.*

$\mathbf{T}$ *has an inverse* $\mathbf{T}^{-1}$ *in* $\mathbb{C}^\mathbb{C}$ *expressed by the relation*

$$\mathbf{T}^{-1}\mathsf{f}(x) = \mathsf{f}(x+1) \qquad \forall \mathsf{f}\in\mathbb{C}^\mathbb{C} \qquad \text{(translation operator inverse).}$$

$\mathbf{D}$ *has an inverse* $\mathbf{D}^{-1}$ *in* $\mathbb{C}^\mathbb{C}$ *expressed by the relation*

$$\mathbf{D}^{-1}\mathsf{f}(x) = \frac{\sqrt{2}}{2}\mathsf{f}\left(\tfrac{1}{2}x\right) \qquad \forall \mathsf{f}\in\mathbb{C}^\mathbb{C} \qquad \text{(dilation operator inverse).}$$

**Proposition 4.4** [98] *Let* $\mathbf{T}$ *and* $\mathbf{D}$ *be as defined in Definition 4.1 page 38. Let* $\mathbf{D}^0 = \mathbf{T}^0 \triangleq \mathbf{I}$ *be the* IDENTITY OPERATOR.

$$\mathbf{D}^j\mathbf{T}^n\mathsf{f}(x) = 2^{j/2}\mathsf{f}\left(2^jx-n\right) \qquad \forall j,n\in\mathbb{Z}, \mathsf{f}\in\mathbb{C}^\mathbb{C}$$

**Example 4.5** (linear functions) [99] Let $\mathbf{T}$ be the *translation operator* (Definition 4.1 page 38). Let $\mathcal{L}(\mathbb{C},\mathbb{C})$ be the set of all *linear* functions in $\boldsymbol{L}^2_\mathbb{R}$.

1. $\{x, \mathbf{T}x\}$ is a *basis* for $\mathcal{L}(\mathbb{C},\mathbb{C})$ and
2. $\mathsf{f}(x) = \mathsf{f}(1)x - \mathsf{f}(0)\mathbf{T}x \qquad \forall \mathsf{f}\in\mathcal{L}(\mathbb{C},\mathbb{C})$

✎PROOF: By left hypothesis, $\mathsf{f}$ is *linear*; so let $\mathsf{f}(x)\triangleq ax+b$

$$\begin{aligned}
\mathsf{f}(1)x - \mathsf{f}(0)\mathbf{T}x &= \mathsf{f}(1)x - \mathsf{f}(0)(x-1) && \text{by Definition 4.1 page 38}\\
&= (ax+b)|_{x=1}\,x - (ax+b)|_{x=0}\,(x-1) && \text{by left hypothesis and definition of } \mathsf{f}\\
&= (a+b)x - b(x-1)\\
&= ax+bx-bx+b\\
&= ax+b\\
&= \mathsf{f}(x) && \text{by left hypothesis and definition of } \mathsf{f}
\end{aligned}$$

✎

**Example 4.6** (Cardinal Series) Let $\mathbf{T}$ be the *translation operator* (Definition 4.1 page 38). The *Paley-Wiener* class of functions $\boldsymbol{PW}^2_\sigma$ are those functions which are "*bandlimited*" with respect to their Fourier transform. The cardinal series forms an orthogonal basis for such a space. The *Fourier coefficients* for a projection of a function $\mathsf{f}$ onto the Cardinal series basis elements is particularly simple—these coefficients are samples of $\mathsf{f}(x)$ taken at regular intervals. In fact, one could represent the coefficients using inner product notation with

---

[96] [75], page 3
[97] [75], page 3
[98] [75], page 4
[99] [86], page 2

the *Dirac delta distribution* $\delta$ as follows:

$$\langle \mathsf{f}(x) \mid \mathbf{T}^n\delta(x)\rangle \triangleq \int_{\mathbb{R}} \mathsf{f}(x)\delta(x-n)\,\mathrm{d}x \triangleq \mathsf{f}(n)$$

1. $\left\{ \mathbf{T}^n \frac{\sin(\pi x)}{\pi x} \middle| n\in\mathbb{N} \right\}$ is a *basis* for $\boldsymbol{PW}^2_\sigma$ and
2. $\mathsf{f}(x) = \underbrace{\sum_{n=1}^{\infty} \mathsf{f}(n)\mathbf{T}^n \frac{\sin(\pi x)}{\pi x}}_{\text{Cardinal series}} \qquad \forall \mathsf{f}\in\boldsymbol{PW}^2_\sigma,\ \sigma\le\frac{1}{2}$

**Example 4.7** (Fourier Series)

1. $\left\{ \mathbf{D}_n e^{ix} \middle| n\in\mathbb{Z} \right\}$ is a *basis* for $\boldsymbol{L}(0,\,2\pi)$ and
2. $\mathsf{f}(x) = \frac{1}{\sqrt{2\pi}} \sum_{n\in\mathbb{Z}} \alpha_n \mathbf{D}_n e^{ix} \qquad \forall x\in(0,2\pi), \mathsf{f}\in\boldsymbol{L}(0,2\pi)$ where
3. $\alpha_n \triangleq \frac{1}{\sqrt{2\pi}} \int_0^{2\pi} \mathsf{f}(x)\mathbf{D}_n e^{-ix}\,\mathrm{d}x \qquad \forall \mathsf{f}\in\boldsymbol{L}(0,2\pi)$

**Example 4.8** (Fourier Transform)

1. $\left\{ \mathbf{D}_\omega e^{ix} \middle| \omega\in\mathbb{R} \right\}$ is a *basis* for $\boldsymbol{L}^2_{\mathbb{R}}$ and
2. $\mathsf{f}(x) = \frac{1}{\sqrt{2\pi}} \int_{\mathbb{R}} \tilde{\mathsf{f}}(\omega)\mathbf{D}_x e^{i\omega}\,\mathrm{d}\omega \qquad \forall \mathsf{f}\in\boldsymbol{L}^2_{\mathbb{R}}$ where
3. $\tilde{\mathsf{f}}(\omega) \triangleq \frac{1}{\sqrt{2\pi}} \int_{\mathbb{R}} \mathsf{f}(x)\mathbf{D}_\omega e^{-ix}\,\mathrm{d}x \qquad \forall \mathsf{f}\in\boldsymbol{L}^2_{\mathbb{R}}$

**Example 4.9** (Gabor Transform) [100]

1. $\left\{ \left(\mathbf{T}_\tau e^{-\pi x^2}\right)\left(\mathbf{D}_\omega e^{ix}\right) \middle| \tau,\omega\in\mathbb{R} \right\}$ is a *basis* for $\boldsymbol{L}^2_{\mathbb{R}}$ and
2. $\mathsf{f}(x) = \int_{\mathbb{R}} \mathsf{G}(\tau,\omega)\,\mathbf{D}_x e^{i\omega}\,\mathrm{d}\omega \qquad \forall x\in\mathbb{R}, \mathsf{f}\in\boldsymbol{L}^2_{\mathbb{R}}$ where
3. $\mathsf{G}(\tau,\omega) \triangleq \int_{\mathbb{R}} \mathsf{f}(x)\left(\mathbf{T}_\tau e^{-\pi x^2}\right)\left(\mathbf{D}_\omega e^{-ix}\right)\mathrm{d}x \qquad \forall x\in\mathbb{R}, \mathsf{f}\in\boldsymbol{L}^2_{\mathbb{R}}$

**Example 4.10** (wavelets) Let $\psi(x)$ be a *mother wavelet*.

1. $\left\{ \mathbf{D}^k\mathbf{T}^n\psi(x) \middle| k,n\in\mathbb{Z} \right\}$ is a *basis* for $\boldsymbol{L}^2_{\mathbb{R}}$ and
2. $\mathsf{f}(x) = \sum_{k\in\mathbb{Z}}\sum_{n\in\mathbb{Z}} \alpha_{k,n}\mathbf{D}^k\mathbf{T}^n\psi(x) \qquad \forall \mathsf{f}\in\boldsymbol{L}^2_{\mathbb{R}}$ where
3. $\alpha_n \triangleq \int_{\mathbb{R}} \mathsf{f}(x)\mathbf{D}^k\mathbf{T}^n\psi^*(x)\,\mathrm{d}x \qquad \forall \mathsf{f}\in\boldsymbol{L}^2_{\mathbb{R}}$

100 📚 [143], ⟨Chapter 3⟩

📚 [61], page 32, ⟨Definition 1.69⟩

## 4.2 The Structure of Wavelets

In Fourier analysis, *continuous* dilations (Definition 4.1 page 38) of the *complex exponential* form a *basis* for the *space of square integrable functions* $\boldsymbol{L}^2_{\mathbb{R}}$ such that

$$\boldsymbol{L}^2_{\mathbb{R}} = \operatorname{span}\left\{\mathbf{D}_\omega e^{ix} \,|\, \omega\in\mathbb{R}\right\}.$$

In Fourier series analysis , *discrete* dilations of the complex exponential form a basis for $\boldsymbol{L}^2_{\mathbb{R}}(0,\, 2\pi)$ such that

$$\boldsymbol{L}^2_{\mathbb{R}}(0,\, 2\pi) = \operatorname{span}\left\{\mathbf{D}_j e^{ix} \middle| j\in\mathbb{Z}\right\}.$$

In Wavelet analysis, for some *mother wavelet* (Definition 4.18 page 47) $\psi(x)$,

$$\boldsymbol{L}^2_{\mathbb{R}} = \operatorname{span}\left\{\mathbf{D}_\omega \mathbf{T}_\tau \psi(x) \,|\, \omega, \tau \in \mathbb{R}\right\}.$$

However, the ranges of parameters $\omega$ and $\tau$ can be much reduced to the countable set $\mathbb{Z}$ resulting in a *dyadic* wavelet basis such that for some mother wavelet $\psi(x)$,

$$\boldsymbol{L}^2_{\mathbb{R}} = \operatorname{span}\left\{\mathbf{D}^j \mathbf{T}^n \psi(x) \,|\, j, n \in \mathbb{Z}\right\}.$$

Wavelets that are both *dyadic* and *compactly supported* have the attractive feature that they can be easily implemented in hardware or software by use of the *Fast Wavelet Transform* (Figure 10 page 49).

In 1989, Stéphane G. Mallat introduced the *Multiresolution Analysis* (MRA, Definition 4.12 page 43) method for wavelet construction. The MRA has since become the dominate wavelet construction method. Moreover, P.G. Lemarié has proved that all wavelets with *compact support* are generated by an MRA.[101]

The MRA is an **analysis** of the linear space $\boldsymbol{L}^2_{\mathbb{R}}$. An analysis of a linear space $\boldsymbol{X}$ is any sequence $\left(\!\left(\boldsymbol{V}_j\right)\!\right)_{j\in\mathbb{Z}}$ of linear subspaces of $\boldsymbol{X}$. The partial or complete reconstruction of $\boldsymbol{X}$ from $\left(\!\left(\boldsymbol{V}_j\right)\!\right)_{j\in\mathbb{Z}}$ is a **synthesis**.[102] Some analyses are completely *characterized* by a *transform*. For example, a Fourier analysis is a sequence of subspaces with sinusoidal bases. Examples of subspaces in a Fourier analysis include $\boldsymbol{V}_1 = \operatorname{span}\left\{e^{ix}\right\}$, $\boldsymbol{V}_{2.3} = \operatorname{span}\left\{e^{i2.3x}\right\}$, $\boldsymbol{V}_{\sqrt{2}} = \operatorname{span}\left\{e^{i\sqrt{2}x}\right\}$, etc. A **transform** is loosely defined as a function that maps a family of functions into an analysis. A very useful transform (a "*Fourier transform*") for Fourier Analysis is

$$\left[\tilde{\mathbf{F}}\mathsf{f}\right](\omega) \triangleq \frac{1}{\sqrt{2\pi}} \int_{\mathbb{R}} \mathsf{f}(x) e^{-i\omega x}\, \mathrm{d}x$$

---

101 [109], [119], page 240

102 The word *analysis* comes from the Greek word ἀνάλυσις, meaning "dissolution" ( [140], page 23, ⟨entry 359⟩), which in turn means "the resolution or separation into component parts" ( [21], http://dictionary.reference.com/browse/dissolution)

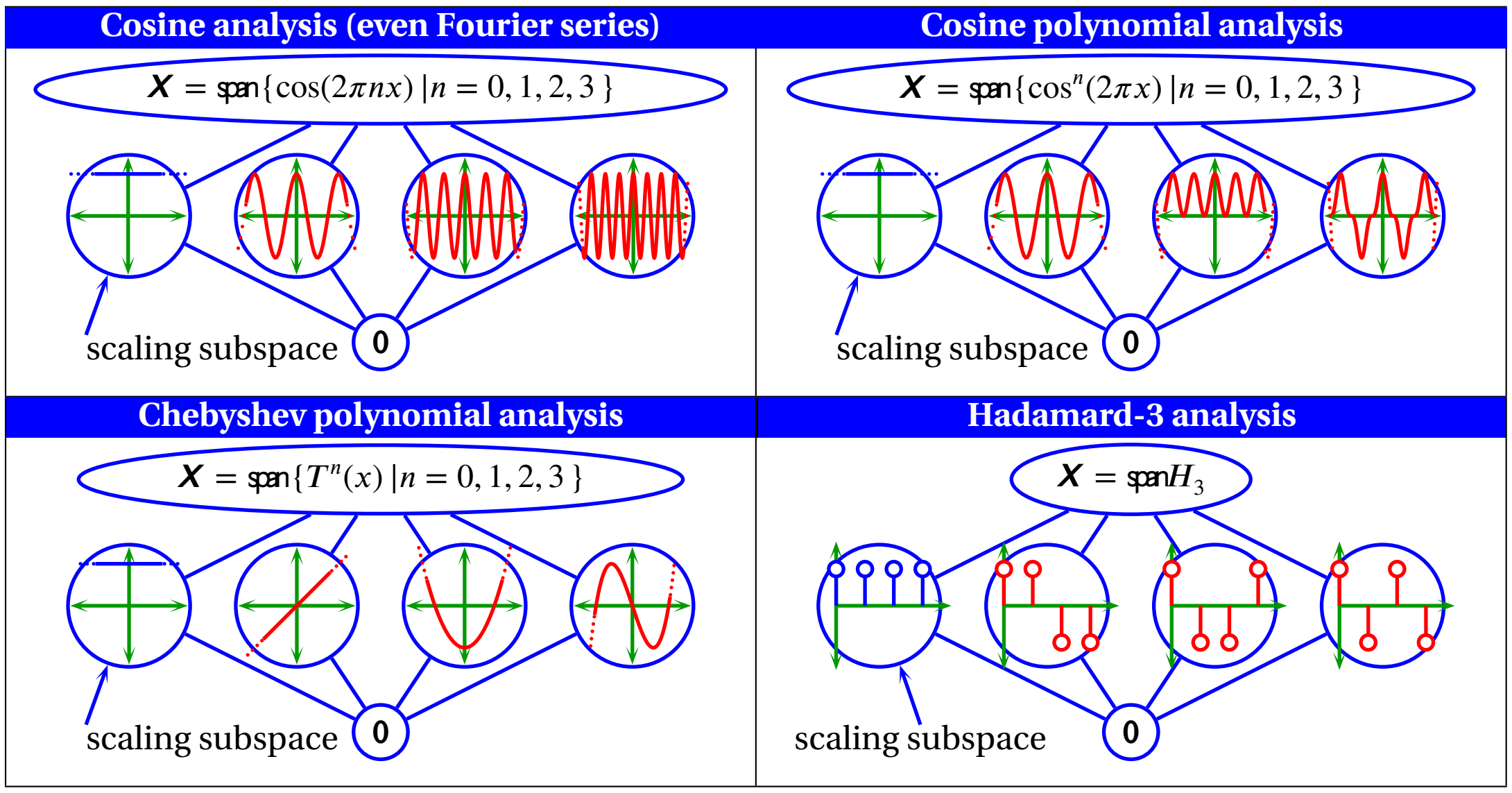

Figure 8: Examples of order structures for selected analyses (Example 4.11 page 42)

An analysis can be partially characterized by its order structure with respect to an order relation such as the set inclusion relation $\subseteq$. Most transforms have a very simple M-$n$ order structure, as illustrated to the right.[103] The M-$n$ lattices for $n \geq 3$ are *modular* (Lemma 1.56 page 16) but not *distributive* (Theorem 1.57 page 16). Analyses typically have one subspace that is a *scaling* subspace; and this subspace is often simply a family of constants (as is the case with Fourier Analysis).

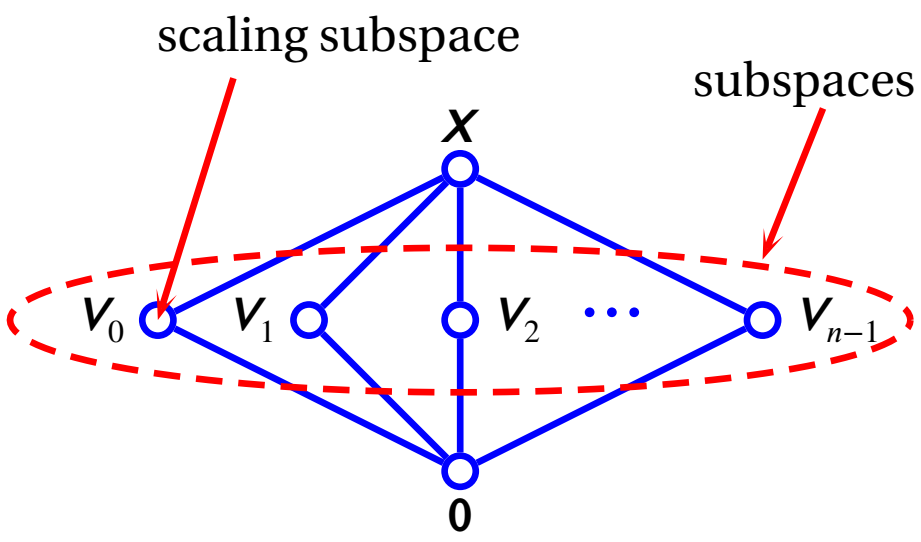

An analysis can be represented using three different structures:

① sequence of subspaces
② sequence of basis vectors
③ sequence of basis coefficients

These structures are isomorphic to each other, and can therefore be used interchangeably.

**Example 4.11** [104] Some examples of the order structures of some analyses are illustrated in Figure 8 (page 42).

[103] [75], page 29, ⟨§2.2⟩
[104] [75], pages 30–31

## 4.3 Multiresolution analysis

A multiresolution analysis provides "coarse" approximations of a function in a linear space $L^2_{\mathbb{R}}$ at multiple "scales" or "resolutions". Key to this process is a sequence of *scaling functions*. Most traditional transforms feature a single *scaling function* $\phi(x)$ set equal to one ($\phi(x) = 1$). This allows for convenient representation of the most basic functions, such as constants.[105] A multiresolution system, on the other hand, uses a generalized form of the scaling concept:[106]

(1) Instead of the scaling function simply being set *equal to unity* ($\phi(x) = 1$), a multiresolution analysis (Definition 4.12 page 43) is often constructed in such a way that the scaling function $\phi(x)$ forms a *partition of unity* such that $\sum_{n\in\mathbb{Z}} \mathbf{T}^n\phi(x) = 1$.

(2) Instead of there being *just one* scaling function, there is an entire sequence of scaling functions $\left(\!\left(\mathbf{D}^j\phi(x)\right)\!\right)_{j\in\mathbb{Z}}$, each corresponding to a different "*resolution*".

**Definition 4.12** [107] Let $\left(\!\left(V_j\right)\!\right)_{j\in\mathbb{Z}}$ be a sequence of subspaces on $L^2_{\mathbb{R}}$. Let $A^-$ be the *closure* of a set $A$. The sequence $\left(\!\left(V_j\right)\!\right)_{j\in\mathbb{Z}}$ is a **multiresolution analysis** on $L^2_{\mathbb{R}}$ if

1. $V_j = V_j^-$ $\quad \forall j\in\mathbb{Z}$ (*closed*) and
2. $V_j \subset V_{j+1}$ $\quad \forall j\in\mathbb{Z}$ (*linearly ordered*) and
3. $\left(\bigcup_{j\in\mathbb{Z}} V_j\right)^- = L^2_{\mathbb{R}}$ (*dense* in $L^2_{\mathbb{R}}$) and
4. $f \in V_j \iff \mathbf{D}f \in V_{j+1}$ $\quad \forall j\in\mathbb{Z}, f\in L^2_{\mathbb{R}}$ (*self-similar*) and
5. $\exists\phi$ such that $\left\{\mathbf{T}^n\phi|_{n\in\mathbb{Z}}\right\}$ is a *Riesz basis* for $V_0$.

A *multiresolution analysis* is also called an **MRA**. An element $V_j$ of $\left(\!\left(V_j\right)\!\right)_{j\in\mathbb{Z}}$ is a **scaling subspace** of the space $L^2_{\mathbb{R}}$. The pair $\left(L^2_{\mathbb{R}}, \left(\!\left(V_j\right)\!\right)\right)$ is a **multiresolution analysis space**, or **MRA space**. The function $\phi$ is the **scaling function** of the *MRA space*.

The traditional definition of the *MRA* also includes the following:

6. $f \in V_j \iff \mathbf{T}^n f \in V_j$ $\quad \forall n,j\in\mathbb{Z}, f\in L^2_{\mathbb{R}}$ (*translation invariant*)
7. $\bigcap_{j\in\mathbb{Z}} V_j = \{\mathbb{0}\}$ (*greatest lower bound* is **0**)

However, these follow from the *MRA* as defined in Definition 4.12 (Proposition 4.13 page 44, Proposition 4.14 page 44).

---

[105] [95], page 8

[106] The concept of a scaling space was perhaps first introduced by Taizo Iijima in 1959 in Japan, and later as the *Gaussian Pyramid* by Burt and Adelson in the 1980s in the West. [118], page 70, [92], [24], [4], [111], [6], [80], [166], ⟨historical survey⟩

[107] [85], page 44, [119], page 221, ⟨Definition 7.1⟩ , [118], page 70, [122], page 21, ⟨Definition 2.2.1⟩, [27], page 284, ⟨Definition 13.1.1⟩, [8], pages 451–452, ⟨Definition 7.7.6⟩, [163], pages 300–301, ⟨Definition 10.16⟩, [35], pages 129–140, ⟨Riesz basis: page 139⟩

**Proposition 4.13** [108] *Let* MRA *be defined as in Definition 4.12 page 43.*

$$\left\{ \left(\!\left( \boldsymbol{V}_j \right)\!\right)_{j\in\mathbb{Z}} \textit{ is an } \text{MRA} \right\} \implies \underbrace{\left\{ \mathsf{f} \in \boldsymbol{V}_j \iff \mathbf{T}^n \mathsf{f} \in \boldsymbol{V}_j \quad \forall n,j\in\mathbb{Z}, \mathsf{f}\in\boldsymbol{L}^2_{\mathbb{R}} \right\}}_{\text{TRANSLATION INVARIANT}}$$

**Proposition 4.14** [109] *Let* MRA *be defined as in Definition 4.12 page 43.*

$$\left\{ \left(\!\left( \boldsymbol{V}_j \right)\!\right)_{j\in\mathbb{Z}} \textit{ is an } \text{MRA} \right\} \implies \left\{ \bigcap_{j\in\mathbb{Z}} \boldsymbol{V}_j = \{\mathbb{0}\} \qquad (\text{GREATEST LOWER BOUND } \textit{is } \mathbf{0}) \right\}$$

The MRA (Definition 4.12 page 43) is more than just an interesting mathematical toy. Under some very "reasonable" conditions (next proposition), as $j \to \infty$, the *scaling subspace* $\boldsymbol{V}_j$ is *dense* in $\boldsymbol{L}^2_{\mathbb{R}}$ ...meaning that with the MRA we can represent any "reasonable" function to within an arbitrary accuracy.

**Proposition 4.15** [110]

$$\left\{ \begin{array}{lll} (1). & \left(\!\left(\mathbf{T}^n \phi\right)\!\right) \textit{ is a } \text{RIESZ SEQUENCE} & \textit{and} \\ (2). & \tilde{\phi}(\omega) \textit{ is } \text{CONTINUOUS } \textit{at } 0 & \textit{and} \\ (3). & \tilde{\phi}(0) \neq 0 & \end{array} \right\} \implies \left\{ \left( \bigcup_{j\in\mathbb{Z}} \boldsymbol{V}_j \right)^{-} = \boldsymbol{L}^2_{\mathbb{R}} \qquad (\text{DENSE } \textit{in } \boldsymbol{L}^2_{\mathbb{R}}) \right\}$$

A *multiresolution analysis* (Definition 4.12 page 43) together with the set inclusion relation $\subseteq$ form the *linearly ordered set* (Definition 1.4 page 4) $\left( \left(\!\left( \boldsymbol{V}_j \right)\!\right), \subseteq \right)$, illustrated to the right by a *Hasse diagram* (Definition 1.6 page 4). Subspaces $\boldsymbol{V}_j$ increase in "size" with increasing $j$. That is, they contain more and more vectors (functions) for larger and larger $j$—with the upper limit of this sequence being $\boldsymbol{L}^2_{\mathbb{R}}$. Alternatively, we can say that approximation within a subspace $\boldsymbol{V}_j$ yields greater "*resolution*" for increasing $j$.[111]

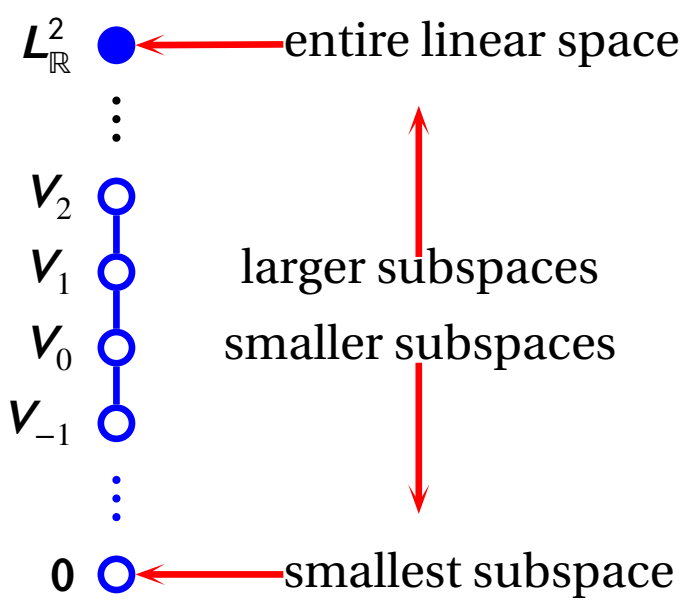

**Remark 4.16** [112]Note that the *greatest lower bound* (*g.l.b.*) of the linearly ordered set $\left( \left(\!\left( \boldsymbol{V}_j \right)\!\right), \subseteq \right)$ is $\mathbf{0}$ (Proposition 4.14 page 44): All linear subspaces contain the zero vector. So the intersection of any two subspaces must at least contain $\mathbb{0}$. If the intersection of any two linear subspaces $\boldsymbol{X}$ and $\boldsymbol{Y}$ is exactly $\{\mathbb{0}\}$, then for any vector in the sum of those subspaces ($\boldsymbol{u} \in \boldsymbol{X} \dotplus \boldsymbol{Y}$) there are **unique** vectors $\mathsf{f} \in \boldsymbol{X}$ and $\mathsf{g} \in \boldsymbol{Y}$ such that $\boldsymbol{u} = \mathsf{f} + \mathsf{g}$. This is *not* necessarily true if the intersection contains more than just $\{\mathbb{0}\}$ .

---

[108] [85], page 45, ⟨Theorem 1.6⟩, [75], pages 32–33, ⟨Proposition 2.1⟩

[109] [168], pages 19–28, ⟨Proposition 2.14⟩, [85], page 45, ⟨Theorem 1.6⟩, [141], pages 313–314, ⟨Lemma 6.4.28⟩, [75], pages 33–35, ⟨Proposition 2.2⟩

[110] [168], pages 28–31, ⟨Proposition 2.15⟩, [75], pages 35–37, ⟨Proposition 2.3⟩

[111] [123], page 83, ⟨Theorem 3.2.12⟩, [106], page 67, ⟨Theorem 2.14⟩, [76], ⟨Theorem 7.1⟩

[112] [75], page 38, ⟨§2.3.2 Order structure⟩

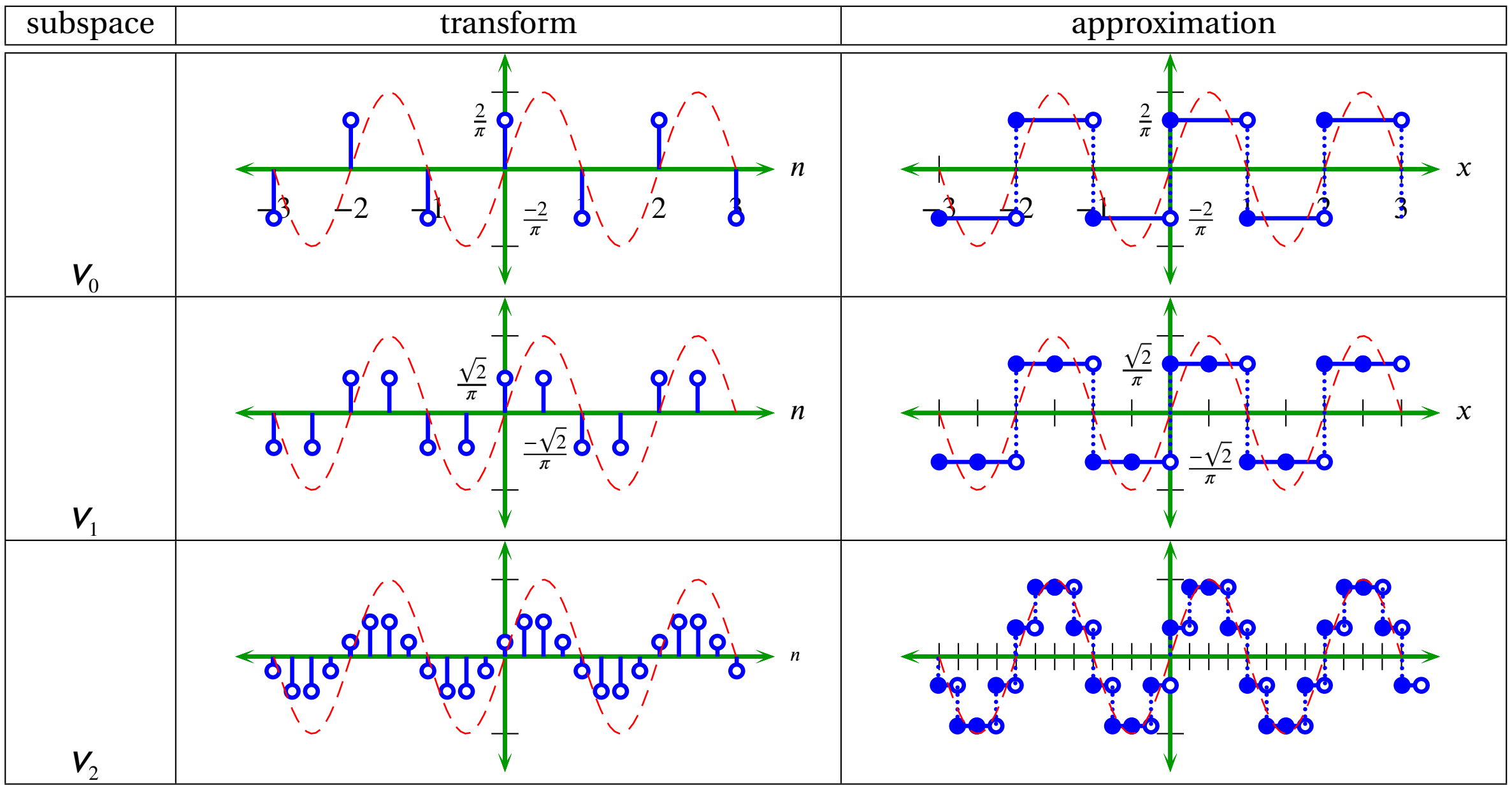

Figure 9: Example approximations of $\sin(\pi x)$ in 3 Haar scaling subspaces (see Example 4.17 page 45)

**Example 4.17**

In the *Haar* MRA, the scaling function $\phi(x)$ is the *pulse function*

$$\phi(x) = \begin{cases} 1 & \text{for } x \in [0,\ 1) \\ 0 & \text{otherwise.} \end{cases}$$

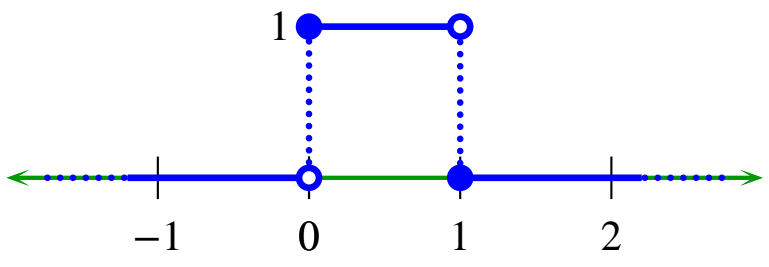

In the subspace $\boldsymbol{V}_j$ ($j \in \mathbb{Z}$) the scaling functions are

$$\mathbf{D}^j\phi(x) = \begin{cases} (2)^{j/2} & \text{for } x \in \left[0,\ \left(2^{-j}\right)\right) \\ 0 & \text{otherwise.} \end{cases}$$

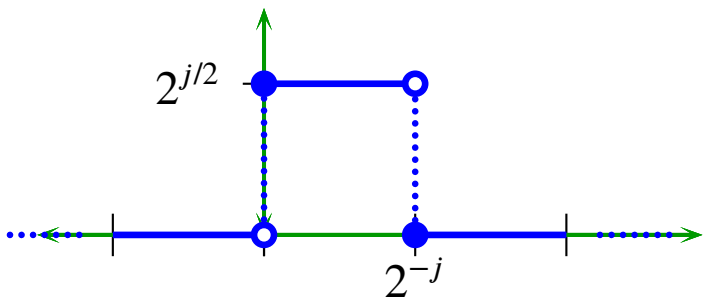

The scaling subspace $\boldsymbol{V}_0$ is the span $\boldsymbol{V}_0 \triangleq \operatorname{span}\left\{\mathbf{T}^n\phi|_{n\in\mathbb{Z}}\right\}$. The scaling subspace $\boldsymbol{V}_j$ is the span $\boldsymbol{V}_j \triangleq \operatorname{span}\left\{\mathbf{D}^j\mathbf{T}^n\phi\,|n \in \mathbb{Z}\right\}$. Note that $\left\|\mathbf{D}^j\mathbf{T}^n\phi\right\|$ for each resolution $j$ and shift $n$ is unity:

$$\begin{aligned}
\left\|\mathbf{D}^j\mathbf{T}^n\phi\right\|^2 &= \|\phi\|^2 \\
&= \int_0^1 |1|^2\,\mathrm{d}x && \text{by definition of } \|\cdot\| \text{ on } \boldsymbol{L}^2_{\mathbb{R}} \\
&= 1
\end{aligned}$$

Let $\mathsf{f}(x) = \sin(\pi x)$. Suppose we want to project $\mathsf{f}(x)$ onto the subspaces $\boldsymbol{V}_0$, $\boldsymbol{V}_1$, $\boldsymbol{V}_2$, ....

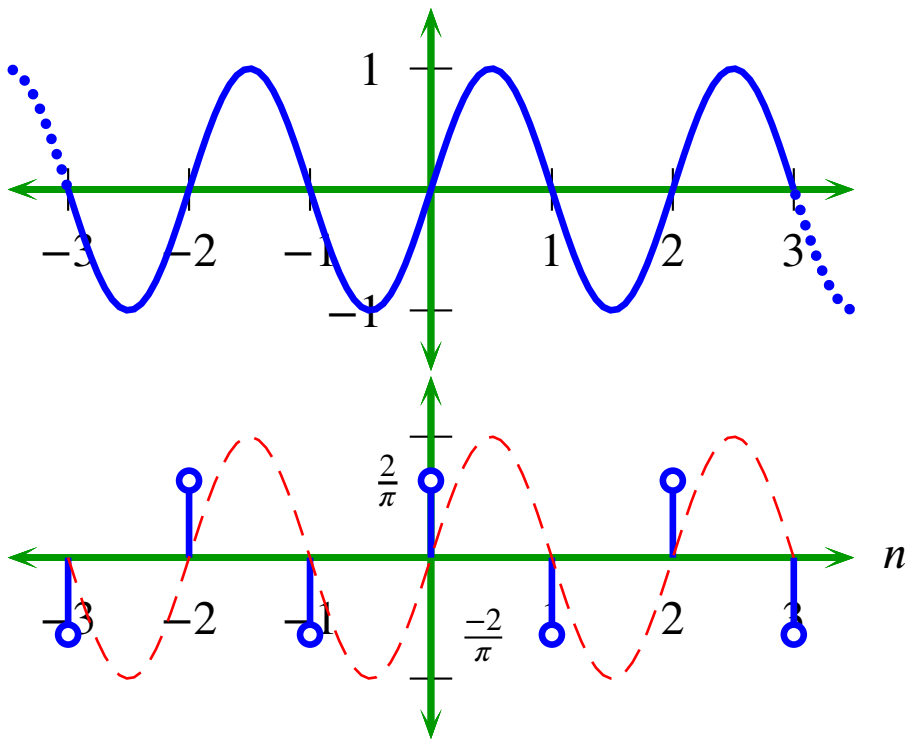

The values of the transform coefficients for the subspace $\boldsymbol{V}_j$ are given by

$$
\begin{aligned}
\left[\mathbf{R}_j \mathsf{f}(x)\right](n) &= \frac{1}{\|\mathbf{D}^j\mathbf{T}^n\phi\|^2}\left\langle \mathsf{f}(x) \mid \mathbf{D}^j\mathbf{T}^n\phi\right\rangle \\
&= \frac{1}{\|\phi\|^2}\left\langle \mathsf{f}(x) \mid 2^{j/2}\phi\left(2^j x - n\right)\right\rangle && \text{by Proposition 4.4 page 39} \\
&= 2^{j/2}\left\langle \mathsf{f}(x) \mid \phi\left(2^j x - n\right)\right\rangle \\
&= 2^{j/2}\int_{2^{-j}n}^{2^{-j}(n+1)} \mathsf{f}(x)\,\mathrm{d}x \\
&= 2^{j/2}\int_{2^{-j}n}^{2^{-j}(n+1)} \sin(\pi x)\,\mathrm{d}x \\
&= 2^{j/2}\left(-\frac{1}{\pi}\right)\cos\left(\pi x\right)\Big|_{2^{-j}n}^{2^{-j}(n+1)} \\
&= \frac{2^{j/2}}{\pi}\left[\cos\left(2^{-j}n\pi\right) - \cos\left(2^{-j}(n+1)\pi\right)\right]
\end{aligned}
$$

And the projection $\mathbf{A}_n\mathsf{f}(x)$ of the function $\mathsf{f}(x)$ onto the subspace $\boldsymbol{V}_j$ is

$$
\begin{aligned}
\mathbf{A}_j\mathsf{f}(x) &= \sum_{n\in\mathbb{Z}}\left\langle \mathsf{f}(x) \mid \mathbf{D}^j\mathbf{T}^n\phi\right\rangle \mathbf{D}^j\mathbf{T}^n\phi \\
&= \frac{2^{j/2}}{\pi}\sum_{n\in\mathbb{Z}}\left[\cos\left(2^{-j}n\pi\right) - \cos\left(2^{-j}(n+1)\pi\right)\right]2^{j/2}\phi\left(2^j x - n\right) \\
&= \frac{2^{j}}{\pi}\sum_{n\in\mathbb{Z}}\left[\cos\left(2^{-j}n\pi\right) - \cos\left(2^{-j}(n+1)\pi\right)\right]\phi\left(2^j x - n\right)
\end{aligned}
$$

The transforms into the subspaces $\boldsymbol{V}_0$, $\boldsymbol{V}_1$, and $\boldsymbol{V}_2$, as well as the approximations in those subspaces are as illustrated in Figure 9 (page 45).

## 4.4 Wavelet analysis

The term "wavelet" comes from the French word "*ondelette*", meaning "small wave". And in essence, wavelets are "small waves" (as opposed to the "long waves" of Fourier analysis) that form a basis for the Hilbert space $\boldsymbol{L}^2_{\mathbb{R}}$.[113]

**Definition 4.18** [114] Let $\mathbf{T}$ and $\mathbf{D}$ be as defined in Definition 4.1 page 38. A function $\psi(x)$ in $\boldsymbol{L}^2_{\mathbb{R}}$ is a **wavelet function** for $\boldsymbol{L}^2_{\mathbb{R}}$ if

$$\left\{\mathbf{D}^j\mathbf{T}^n\psi\,|_{j,n\in\mathbb{Z}}\right\} \text{ is a \textit{Riesz basis} for } \boldsymbol{L}^2_{\mathbb{R}}.$$

In this case, $\psi$ is also called the **mother wavelet** of the basis $\left\{\mathbf{D}^j\mathbf{T}^n\psi\,|_{j,n\in\mathbb{Z}}\right\}$. The sequence of subspaces $\left(\!\left(\boldsymbol{W}_j\right)\!\right)_{j\in\mathbb{Z}}$ is the **wavelet analysis** induced by $\psi$, where each subspace $\boldsymbol{W}_j$ is defined as

$$\boldsymbol{W}_j \triangleq \operatorname{span}\left\{\mathbf{D}^j\mathbf{T}^n\psi\,\middle|\,n\in\mathbb{Z}\right\}.$$

A *wavelet analysis* $\left(\!\left(\boldsymbol{W}_j\right)\!\right)$ is often constructed from a *multiresolution anaysis* (Definition 4.12 page 43) $\left(\!\left(\boldsymbol{V}_j\right)\!\right)$ under the relationship

$$\boldsymbol{V}_{j+1} = \boldsymbol{V}_j \dotplus \boldsymbol{W}_j, \qquad \text{where } \dotplus \text{ is subspace addition (\textit{Minkowski addition}).}$$

By this relationship alone, $\left(\!\left(\boldsymbol{W}_j\right)\!\right)$ is in no way uniquely defined in terms of a multiresolution analysis $\left(\!\left(\boldsymbol{V}_j\right)\!\right)$. In general there are many possible complements of a subspace $\boldsymbol{V}_j$. To uniquely define such a wavelet subspace, one or more additional constraints are required. One of the most common additional constraints is *orthogonality*, such that $\boldsymbol{V}_j$ and $\boldsymbol{W}_j$ are orthogonal to each other.

**Definition 4.19** Let $\left(\boldsymbol{L}^2_{\mathbb{R}}, \left(\!\left(\boldsymbol{V}_j\right)\!\right), \phi, \left(\!\left(h_n\right)\!\right)\right)$ be a multiresolution system (Definition 4.12 page 43) and $\left(\!\left(\boldsymbol{W}_j\right)\!\right)_{j\in\mathbb{Z}}$ a wavelet analysis (Definition 4.18 page 47) with respect to $\left(\!\left(\boldsymbol{V}_j\right)\!\right)_{j\in\mathbb{Z}}$. Let $\left(\!\left(g_n\right)\!\right)_{n\in\mathbb{Z}}$ be a sequence of coefficients such that $\psi = \sum_{n\in\mathbb{Z}} g_n\mathbf{D}\mathbf{T}^n\phi$.

A **wavelet system** is the tuple

$$\left(\boldsymbol{L}^2_{\mathbb{R}}, \left(\!\left(\boldsymbol{V}_j\right)\!\right), \left(\!\left(\boldsymbol{W}_j\right)\!\right), \phi, \psi, \left(\!\left(h_n\right)\!\right), \left(\!\left(g_n\right)\!\right)\right)$$

and the sequence $\left(\!\left(g_n\right)\!\right)_{n\in\mathbb{Z}}$ is the **wavelet coefficient sequence**.

**Theorem 4.20** [115] *Let* $\left(\boldsymbol{L}^2_{\mathbb{R}}, \left(\!\left(\boldsymbol{V}_j\right)\!\right), \left(\!\left(\boldsymbol{W}_j\right)\!\right), \phi, \psi, \left(\!\left(h_n\right)\!\right), \left(\!\left(g_n\right)\!\right)\right)$ *be a* WAVELET SYSTEM *(Definition 4.19 page 47). Let* $\boldsymbol{V}_1 \dotplus \boldsymbol{V}_2$ *represent* MINKOWSKI ADDITION *of two subspaces* $\boldsymbol{V}_1$ *and* $\boldsymbol{V}_2$ *of a Hilbert space* $\boldsymbol{H}$.

$$\begin{aligned}
\boldsymbol{L}^2_{\mathbb{R}} &= \lim_{j\to\infty} \boldsymbol{V}_j && (\boldsymbol{L}^2_{\mathbb{R}} \text{ is equivalent to one very large scaling subspace})\\
&= \boldsymbol{V}_j \dotplus \boldsymbol{W}_j \dotplus \boldsymbol{W}_{j+1} \dotplus \boldsymbol{W}_{j+2} \dotplus \cdots && \left(\begin{array}{l}\boldsymbol{L}^2_{\mathbb{R}} \text{ is equivalent to one scaling space}\\ \text{and a sequence of wavelet subspaces}\end{array}\right)\\
&= \cdots \dotplus \boldsymbol{W}_{-2} \dotplus \boldsymbol{W}_{-1} \dotplus \boldsymbol{W}_0 \dotplus \boldsymbol{W}_1 \dotplus \boldsymbol{W}_2 \dotplus \cdots && (\boldsymbol{L}^2_{\mathbb{R}} \text{ is equivalent to a sequence of wavelet subspaces})
\end{aligned}$$

113 [158], page ix, [7], page 191
114 [168], page 17, ⟨Definition 2.1⟩, [75], page 50, ⟨Definition 2.4⟩
115 [75], page 53, ⟨Theorem 2.8⟩

✎PROOF:

(1) Proof for (1):

$$\mathbf{L}^2_{\mathbb{R}} = \lim_{j\to\infty} \mathbf{V}_j \qquad \text{by Definition 4.12 page 43}$$

(2) Proof for (2):

$$\begin{aligned}
\underbrace{\mathbf{V}_j \dotplus \mathbf{W}_j}_{\mathbf{V}_{j+1}} \dotplus \mathbf{W}_{j+1} \dotplus \mathbf{W}_{j+2} \dotplus \cdots &= \underbrace{\mathbf{V}_{j+1} \dotplus \mathbf{W}_{j+1}}_{\mathbf{V}_{j+2}} \dotplus \mathbf{W}_{j+2} \dotplus \mathbf{W}_{j+3} \dotplus \cdots \\
&= \underbrace{\mathbf{V}_{j+2} \dotplus \mathbf{W}_{j+2}}_{\mathbf{V}_{j+3}} \dotplus \mathbf{W}_{j+3} \dotplus \mathbf{W}_{j+4} \dotplus \cdots \\
&= \underbrace{\mathbf{V}_{j+3} \dotplus \mathbf{W}_{j+3}}_{\mathbf{V}_{j+4}} \dotplus \mathbf{W}_{j+4} \dotplus \mathbf{W}_{j+5} \dotplus \cdots \\
&= \underbrace{\mathbf{V}_{j+5} \dotplus \mathbf{W}_{j+5}}_{\mathbf{V}_{j+5}} \dotplus \mathbf{W}_{j+6} \dotplus \mathbf{W}_{j+6} \dotplus \cdots \\
&= \lim_{j\to\infty} \mathbf{V}_{j+5} \dotplus \mathbf{W}_{j+5} \dotplus \mathbf{W}_{j+6} \dotplus \mathbf{W}_{j+6} \dotplus \cdots \\
&= \mathbf{L}^2_{\mathbb{R}}
\end{aligned}$$

(3) Proof for (3):

$$\begin{aligned}
\mathbf{L}^2_{\mathbb{R}} &= \underbrace{\mathbf{V}_0}_{\mathbf{V}_{-1} \dotplus \mathbf{W}_{-1}} \dotplus \mathbf{W}_0 \dotplus \mathbf{W}_1 \dotplus \mathbf{W}_2 \dotplus \mathbf{W}_3 \dotplus \cdots && \text{by (2)} \\
&= \underbrace{\mathbf{V}_{-1}}_{\mathbf{V}_{-2} \dotplus \mathbf{W}_{-2}} \mathbf{W}_{-1} \dotplus \mathbf{W}_0 \dotplus \mathbf{W}_1 \dotplus \mathbf{W}_2 \dotplus \mathbf{W}_3 \dotplus \cdots \\
&= \underbrace{\mathbf{V}_{-2}}_{\mathbf{V}_{-3} \dotplus \mathbf{W}_{-3}} \mathbf{W}_{-2} \dotplus \mathbf{W}_{-1} \dotplus \mathbf{W}_0 \dotplus \mathbf{W}_1 \dotplus \mathbf{W}_2 \dotplus \mathbf{W}_3 \dotplus \cdots \\
&= \underbrace{\mathbf{V}_{-3}}_{\mathbf{V}_{-4} \dotplus \mathbf{W}_{-4}} \mathbf{W}_{-3} \dotplus \mathbf{W}_{-2} \dotplus \mathbf{W}_{-1} \dotplus \mathbf{W}_0 \dotplus \mathbf{W}_1 \dotplus \mathbf{W}_2 \dotplus \mathbf{W}_3 \dotplus \cdots \\
&\vdots \\
&= \cdots \dotplus \mathbf{W}_{-3} \dotplus \mathbf{W}_{-2} \dotplus \mathbf{W}_{-1} \dotplus \mathbf{W}_0 \dotplus \mathbf{W}_1 \dotplus \mathbf{W}_2 \dotplus \mathbf{W}_3 \dotplus \cdots
\end{aligned}$$

✏

**Remark 4.21** In the special case that two subspaces $\mathbf{W}_1$ and $\mathbf{W}_2$ are *orthogonal* to each other, then the *subspace addition* operation $\mathbf{W}_1 \dotplus \mathbf{W}_2$ is frequently expressed as $\mathbf{W}_1 \oplus \mathbf{W}_2$. In the case of an *orthonormal wavelet system*, the expressions in Theorem 4.20 (page 47)

could be expressed as

$$\begin{aligned} L^2_{\mathbb{R}} &= \lim_{j\to\infty} V_j \\ &= V_j \oplus W_j \oplus W_{j+1} \oplus W_{j+2} \oplus \cdots \\ &= \cdots \oplus W_{-2} \oplus W_{-1} \oplus W_0 \oplus W_1 \oplus W_2 \oplus \cdots . \end{aligned}$$

.

## 4.5 Fast Wavelet Transform (FWT)

*Filter bank*s can be used to implement a "*Fast Wavelet Transform*" (*FWT*). This is illustrated in Figure 10 page 49.[116]

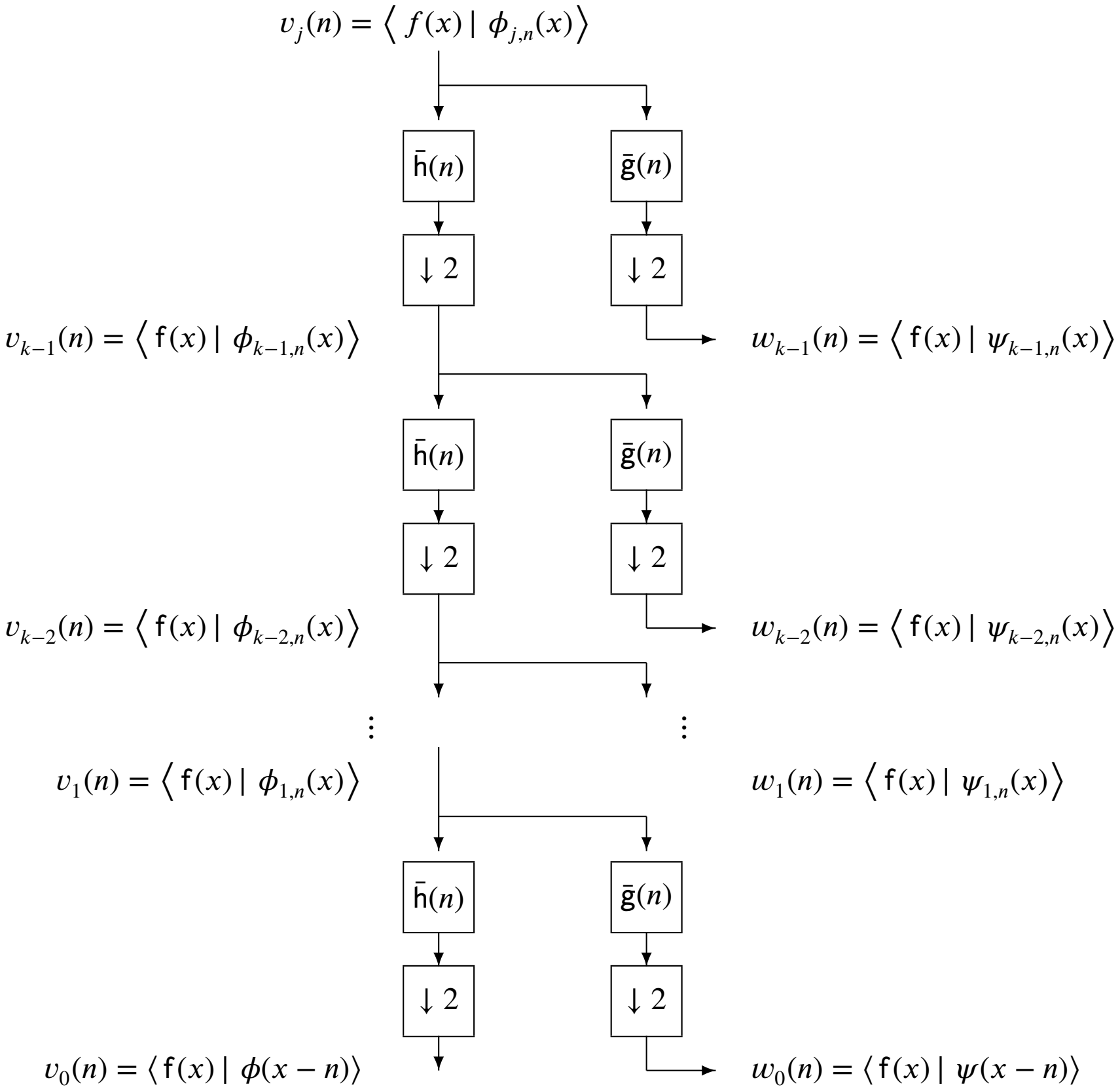

Figure 10: $k$-Stage Fast Wavelet Transform (FWT)

116 [119], page 257, ⟨Figure 7.12⟩, [75], pages 371–372, ⟨Figure L.1⟩

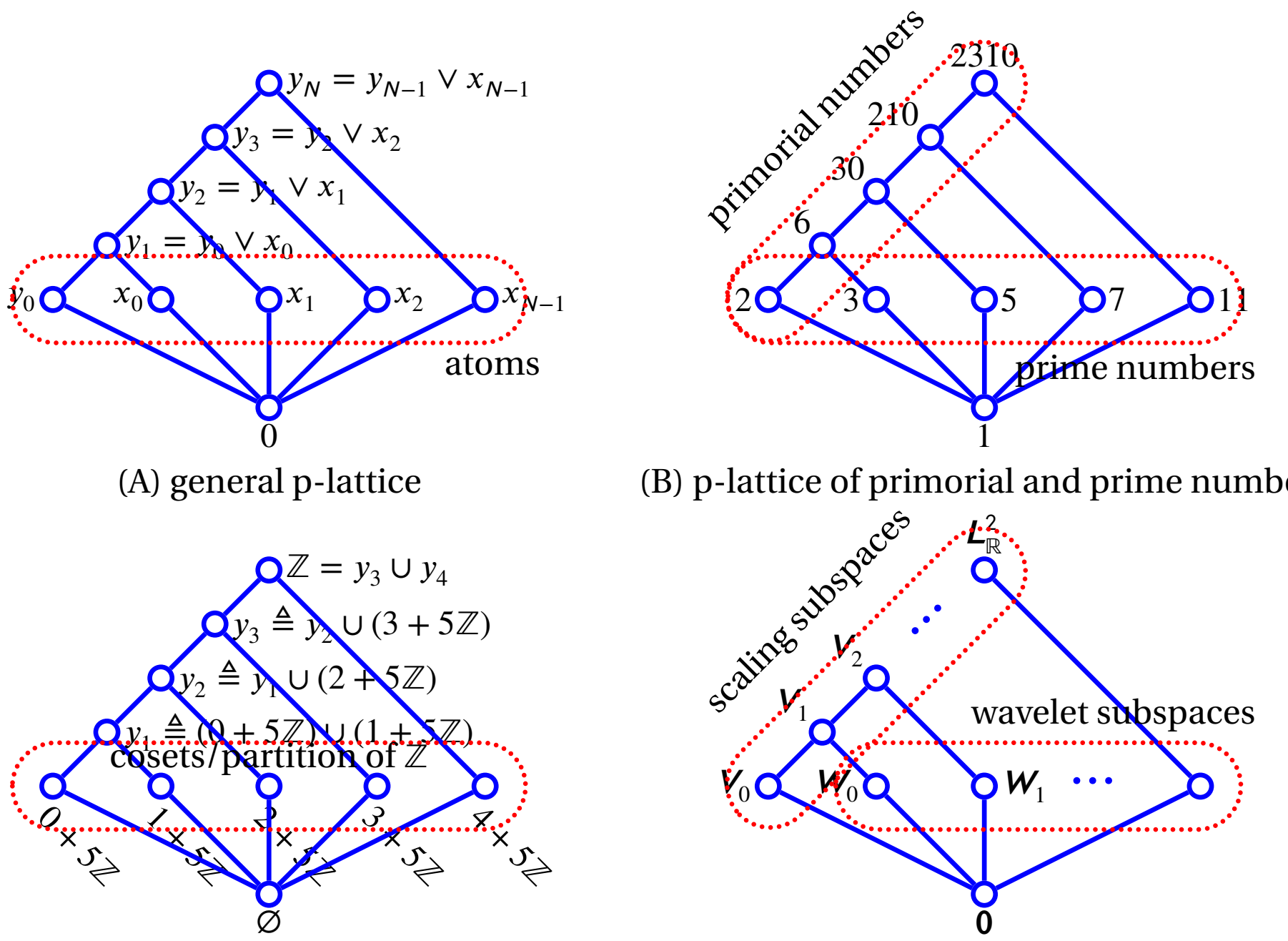

(A) general p-lattice (B) p-lattice of primorial and prime numbers

(C) p-lattice of cosets (D) p-lattice of MRA and wavelet subspaces

Figure 11: Some selected *primorial lattices* (see Example 5.2 page 50–Example 5.5 page 51)

# 5 Main Results

## 5.1 Primorial Lattices

**Definition 5.1** Let $X \triangleq \{0, x_0, x_1, \ldots, x_N, y_0, y_1, \ldots, y_N\}$ be a set.
A *lattice* $\boldsymbol{L} \triangleq (X, \vee, \wedge; \leq)$ is **primorial** if

1. 0 is the *least element* of $\boldsymbol{L}$ and
2. $\boldsymbol{L}$ is *atomic* (Definition 1.44 page 13) and $\{y_0, x_0, x_1, \ldots, x_N\}$ are *atoms* of $\boldsymbol{L}$ and
3. $y_{n+1} = y_n \vee x_n$ .

A lattice that is *primorial* is a **primorial lattice**, or simply a **p-lattice**.

**Example 5.2** A general *primorial lattice* is illustrated to in Figure 11 page 50 (A).

**Example 5.3** [117] The set of *primorial numbers* and *prime numbers* ordered by the *divides* ("|") relation forms a *primorial lattice*, as illustrated in Figure 11 page 50 (B).

117 [75], page 30, [2] ⟨http://oeis.org/A002110⟩

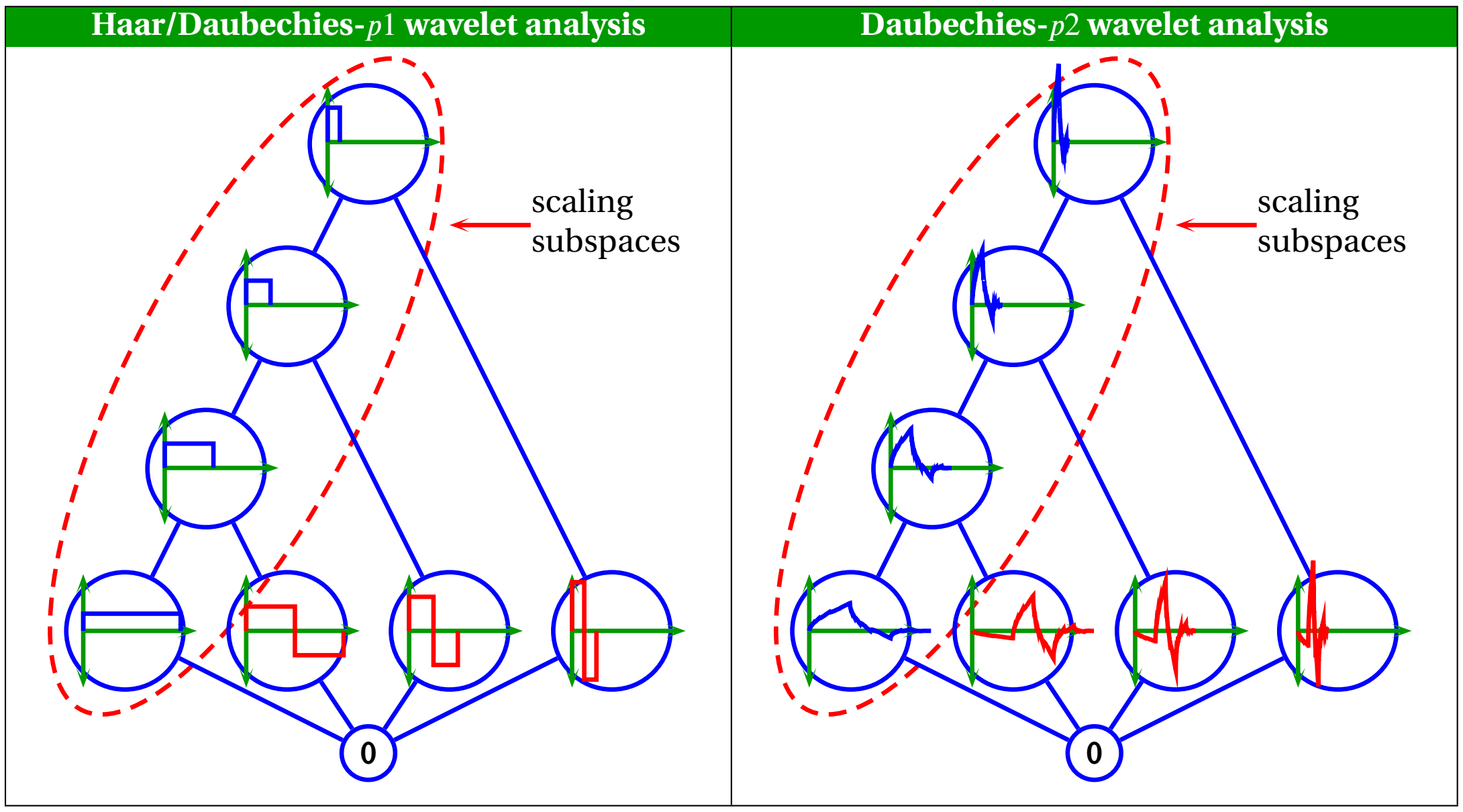

Figure 12: some MRA-wavelet systems

**Example 5.4** Any *partition*, along with successive unions of the partition elements, generates a *primorial lattice*. One example of this is the *cosets* of $\mathbb{Z}$, which generate a *finite* primorial lattice, as illustrated in Figure 11 page 50 (C).

**Example 5.5** A special characteristic of MRA-wavelet analysis is that it's order structure with respect to the $\subseteq$ relation is not a simple $M_n$ *lattice* (as is with the case of Fourier and several other analyses). Rather, it is a *primorial lattice* as illustrated in Figure 11 page 50 (D) and in Figure 12 page 51.

**Proposition 5.6** [118] *Let* $\boldsymbol{L} \triangleq (X, \vee, \wedge; \leq)$ *be a* LATTICE.

$$\left\{ \begin{array}{l} \boldsymbol{L} \text{ is} \\ \textbf{primorial} \end{array} \right\} \implies \left\{ \begin{array}{llll} 1. & \boldsymbol{L} \text{ is NONDISTRIBUTIVE} & \text{(Definition 1.53 page 15)} & \text{and} \\ 2. & \boldsymbol{L} \text{ is NONMODULAR} & \text{(Definition 1.47 page 14)} & \text{and} \\ 3. & \boldsymbol{L} \text{ is COMPLEMENTED} \iff \boldsymbol{L} \text{ is FINITE} & \text{(Definition 1.63 page 17)} & \text{and} \\ 4. & \boldsymbol{L} \text{ is NOT UNIQUELY COMPLEMENTED} & \text{(Definition 1.63 page 17)} & \text{and} \\ 5. & \boldsymbol{L} \text{ is NONORTHOCOMPLEMENTED} & \text{(Definition 1.72 page 20)} & \text{and} \\ 6. & \boldsymbol{L} \text{ is NONBOOLEAN} & \text{(Definition 1.69 page 18)} & . \end{array} \right\}$$

117 [75], page 72, ⟨Section 2.4.3 Order structure⟩
118 [75], page 52, ⟨Proposition 2.6⟩

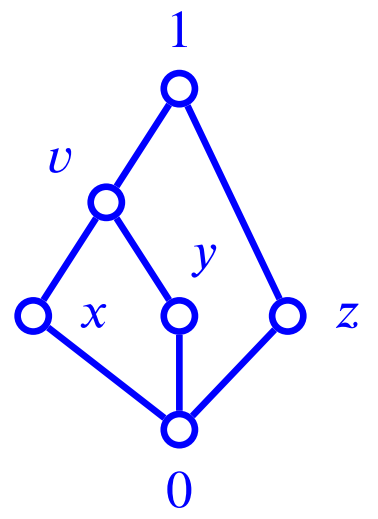

✎PROOF:

(1) Proof that $\mathbf{L}$ is *nondistributive*:
   (a) $\mathbf{L}$ contains the *N5 lattice* (Definition 1.49 page 14).
   (b) Because $\mathbf{L}$ contains the $N5$ lattice, $\mathbf{L}$ is *nondistributive* (Theorem 1.57 page 16).
(2) Proof that $\mathbf{L}$ is *nonmodular* and *nondistributive*:
   (a) $\mathbf{L}$ contains the $N5$ lattice (Definition 1.49 page 14).
   (b) Because $\mathbf{L}$ contains the $N5$ lattice, $\mathbf{L}$ is *nonmodular* (Theorem 1.50 page 14).
(3) Proof that $\mathbf{L}$ is *noncomplemented*:

$$\begin{aligned} x' &= y' = v' = z \\ z' &= \{x, y, v\} \\ x'' &= (x')' \\ &= z' \\ &= \{x, y, v\} \\ &\neq x \end{aligned}$$

(4) Proof that $\mathbf{L}$ is *nonBoolean*:
   (a) $\mathbf{L}$ is *nondistributive* (item 1 page 52).
   (b) Because $\mathbf{L}$ is *nondistributive*, it is *nonBoolean* (Definition 1.69 page 18).

✎

## 5.2 Reduction operator on boolean lattices

**Definition 5.7** Let $\mathbb{B}$ be the set of all *bounded lattices* (Definition 1.39 page 12). Let $\mathbf{L}_2^N \triangleq (X, \vee, \wedge, 0, 1; \leq)$ be a *Boolean lattice* (Definition 1.69 page 18) with $2^N$ elements and $N \in \mathbb{N}$ ($N$ is a positive integer). The operator $\mathbf{R}$ is the **lattice reduction operator** of $\mathbf{L}_2^N$ and $\mathbf{R}\mathbf{L}_2^N$ is the **reduction of** $\mathbf{L}_2^N$ if

$$\mathbf{R}\mathbf{L}_2^N \triangleq \left\{ \mathbf{L} \in \mathbb{B} \;\middle|\; \begin{array}{lll} 1. & \mathbf{L} \text{ is a } 2^{N-1} \text{ element } \textit{Boolean lattice} & \text{and} \\ 2. & \mathbf{L} \subseteq \mathbf{L}_2^N & \text{and} \\ 3. & \{0,1\} \in \mathbf{L} & \text{and} \\ 4. & \{x,y\} \text{ is an } \textit{orthocomplemented pair} \text{ in } \mathbf{L} \implies & \\ & \{x,y\} \text{ is an } \textit{orthocomplemented pair} \text{ in } \mathbf{L}_2^N & \end{array} \right\}$$

Note that in Definition 5.7, the *order relation* $\le$ is the same for both $\boldsymbol{L}_2^N$ and any $\boldsymbol{L}$ in $\mathbf{R}\boldsymbol{L}_2^N$. That is, if $x \le y$ in $\boldsymbol{L}_2^N$, then $x \le y$ in $\boldsymbol{L}$ as well.

**Example 5.8** Let $\boldsymbol{L}_2^2$ be a *Boolean lattice* (Definition 1.69 page 18) of order 2. Let $\mathbf{R}$ be the *lattice reduction operator* $\mathbf{R}$ and $\mathbf{R}\boldsymbol{L}_2^2$ be the *reduction of* $\boldsymbol{L}_2^2$ (Definition 5.7 page 52). Then $\mathbf{R}\boldsymbol{L}_2^2$ yields a set of exactly one $2^{2-1}$ value Boolean lattice, as illustrated next:

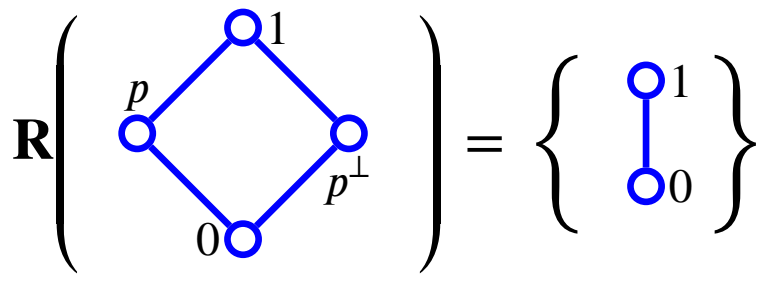

**Example 5.9** Let $\boldsymbol{L}_2^3$ be a *Boolean lattice* (Definition 1.69 page 18) of order 3. Let $\mathbf{R}$ be the **lattice reduction operator R** and $\mathbf{R}\boldsymbol{L}_2^3$ be the **reduction of** $\boldsymbol{L}_2^3$ (Definition 5.7 page 52). The operation $\mathbf{R}\boldsymbol{L}_2^3$ yields a set of three $2^2$ value Boolean lattices, as illustrated next:

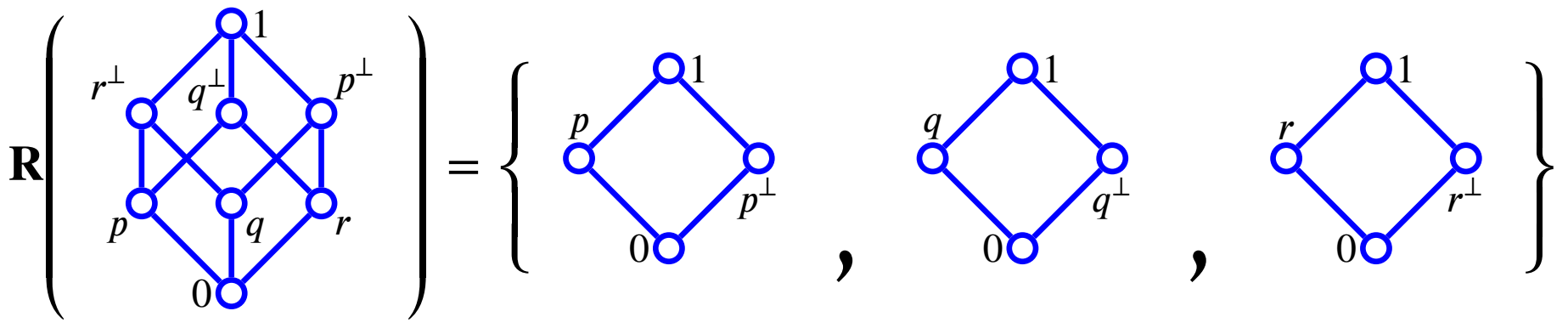

**Example 5.10** Let $\boldsymbol{L}_2^4$ be a *Boolean lattice* (Definition 1.69 page 18) of order 4. Let $\mathbf{R}$ be the **lattice reduction operator R** and $\mathbf{R}\boldsymbol{L}_2^4$ be the **reduction of** $\boldsymbol{L}_2^4$ (Definition 5.7 page 52). The operation $\mathbf{R}\boldsymbol{L}_2^4$ yields a set of ten $2^3$ value Boolean lattices, as illustrated in Figure 13 (page 54).

**Remark 5.11** In a *boolean lattice* $\boldsymbol{L}_2^N$ (Definition 1.69 page 18), besides the pair $\{0,1\}$, there are a total of $2^{N-1} - 1$ *orthocomplemented* (Definition 1.72 page 20) pairs of elements. But note that any arbitrary $2^{N-1} - 2$ pairs of orthocomplemented pairs does not in general generate a *boolean lattice*. The lattice $\boldsymbol{L}_2^4$, for example, has $2^{4-1} - 1 = 7$ orthocomplemented pairs besides $\{0,1\}$. To generate an $\boldsymbol{L}_2^3$ lattice, we need 3 orthocomplemented pairs. There are $\binom{7}{3} = \frac{7!}{3!4!} = 35$ ways of selecting 3 pairs from $\boldsymbol{L}_2^4$, but only 10 of these ways generate a *boolean lattice* (Example 5.10 page 53). All other ways fail.

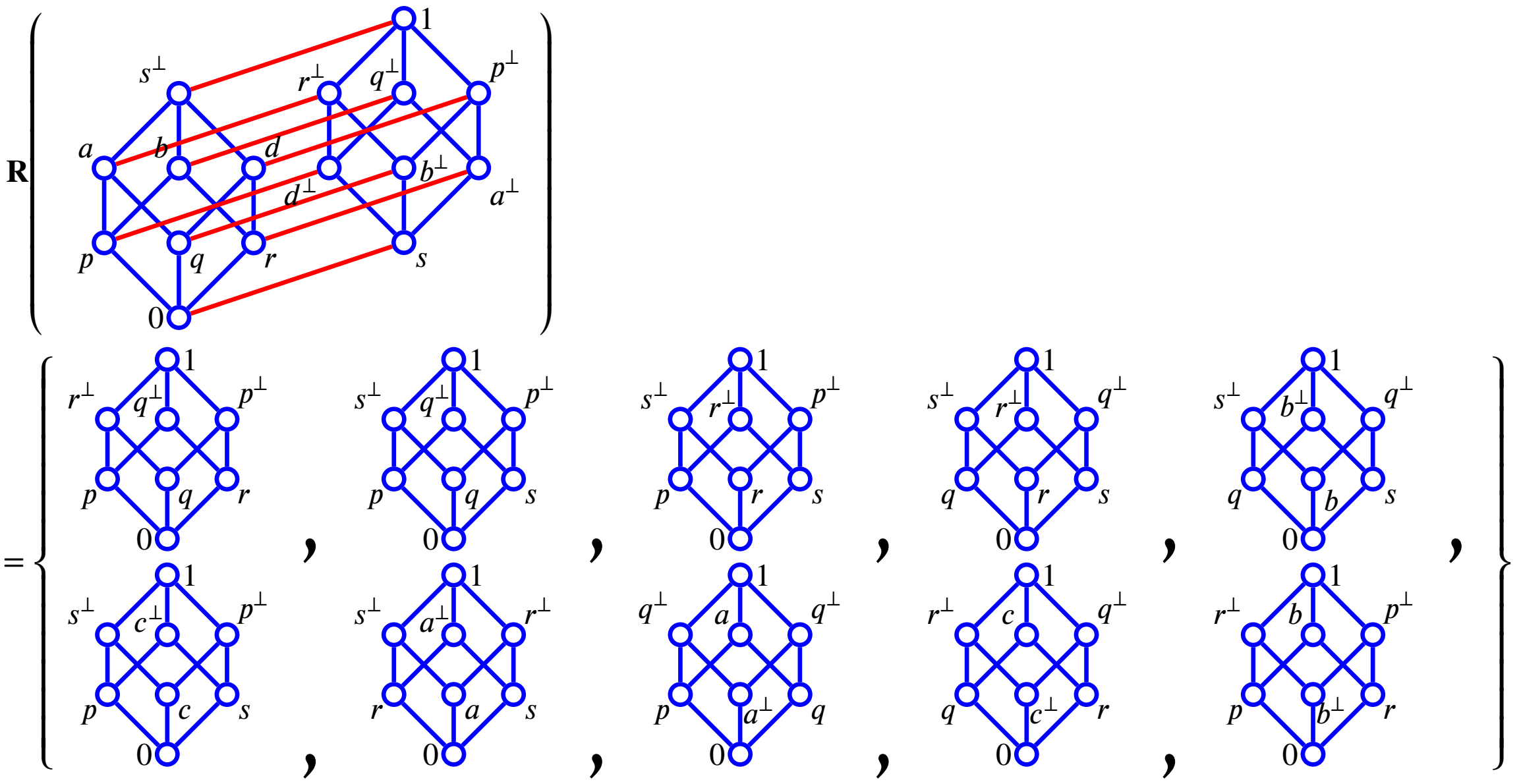

Figure 13: **reduction of $\boldsymbol{L}_2^4$** (Example 5.10 page 53)

For example, if we were to select the pairs $\{0, w, w^\perp, a, a^\perp, b, b^\perp, 1\}$, we would get the *orthocomplemented*, but **non-boolean** (Definition 1.69 page 18) lattice illustrated to the right; In particular, it is *complemented*, but *non-distributive*. For example, $w^\perp \wedge (a \vee b) = w^\perp \neq 0 = 0 \vee 0 = (w^\perp \wedge a) \vee (w^\perp \wedge b)$. Alternatively, note that the set $\{1, a, w, 0, b^\perp, w^\perp\}$ together with the ordering relation $\leq$ form an $O_6$ *sublattice* (Definition 1.73 page 20), which contains an $N_5$ *sublattice*, which implies that the lattice to the right is *non-distributive* (by the *Birkhoff distributivity criterion* Theorem 1.57 page 16).

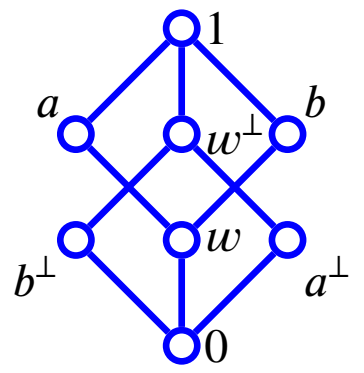

**Example 5.12** Let $\boldsymbol{L}_2^5$ be a *Boolean lattice* (Definition 1.69 page 18) of order 5. Let $\mathbf{R}$ be the **lattice reduction operator R** and $\mathbf{R}\boldsymbol{L}_2^5$ be the **reduction of $\boldsymbol{L}_2^5$** (Definition 5.7 page 52). The result of the operation $\mathbf{R}\boldsymbol{L}_2^5$ is partially illustrated in Figure 14 (page 55).

## 5.3 Difference operator on bounded lattices

**Definition 5.13** Let $X \setminus Y$ be the standard *set difference* of a set $X$ and a set $Y$. Let $\boldsymbol{L}_x \triangleq (X, \vee, \wedge, 0, 1; \leq)$ and $\boldsymbol{L}_y \triangleq (Y, \vee, \wedge, 0, 1; \leq)$ be *bounded lattices* (Definition 1.39 page 12).
The **bounded lattice difference $\boldsymbol{L}_x \oslash \boldsymbol{L}_y$** of $\boldsymbol{L}_x$ and $\boldsymbol{L}_y$ is the lattice $\boldsymbol{L}$ such that

$$\boldsymbol{L} \triangleq ((X \setminus Y) \cup \{0, 1\}, \vee, \wedge, 0, 1; \leq)$$

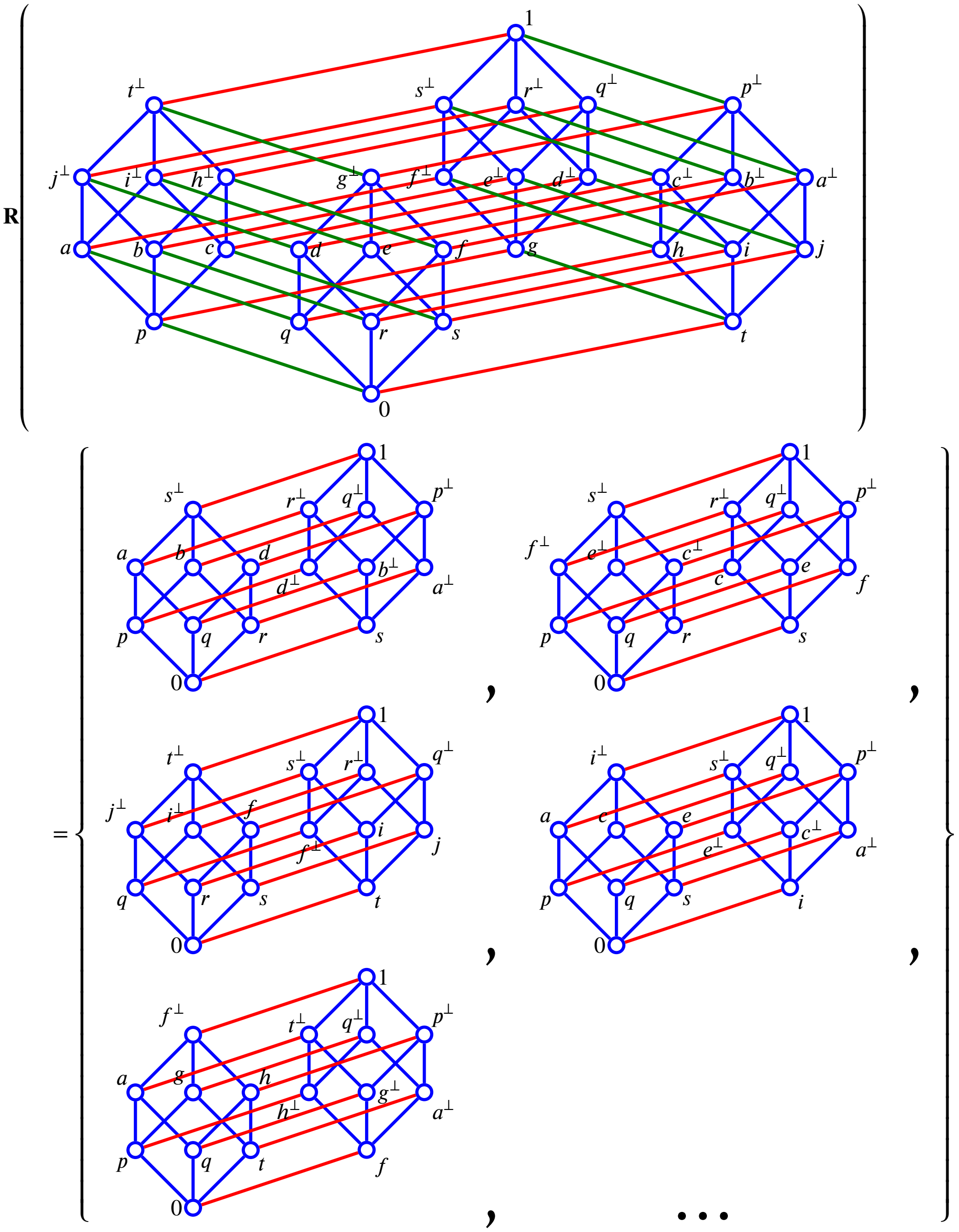

Figure 14: **reduction of $L_2^5$** (Example 5.12 page 54)

**Example 5.14** Let $\oslash$ be the *bounded lattice difference* operator (Definition 5.13 page 54).

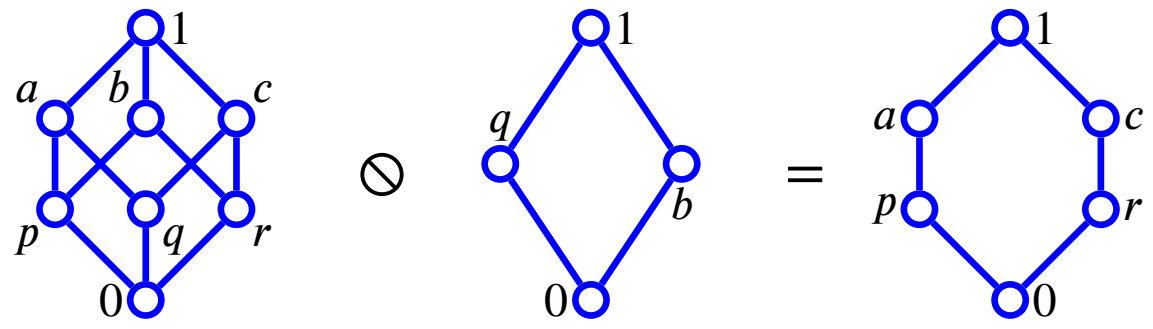

**Proposition 5.15** *Let $\mathbb{B}$ be the set of all* BOUNDED LATTICES *(Definition 1.39 page 12). Let $\oslash$ be the* BOUNDED LATTICE DIFFERENCE *operator (Definition 5.13 page 54).*
$(\mathbb{B}, \oslash, \subseteq)$ *is a* D-POSET *(Definition 3.18 page 37).*

**Theorem 5.16** *Let $\boldsymbol{L} \triangleq \boldsymbol{L}_2^N \oslash \boldsymbol{L}_2^{N-1}$ be the* BOUNDED LATTICE DIFFERENCE *(Definition 5.13 page 54) of a* BOOLEAN LATTICE *$\boldsymbol{L}_2^N$ (Definition 1.69 page 18) and a* BOOLEAN LATTICE *$\boldsymbol{L}_2^{N-1}$ selected from the set $\mathbf{R}\boldsymbol{L}_2^N$ (Definition 5.7 page 52). Let $X \triangleq \left\{ \boldsymbol{L}_2^n \,|\, n = 1, 2, \ldots \right\} \cup \left\{ \boldsymbol{L}_2^n \oslash \boldsymbol{L}_2^{n-1} \,|\, n = 2, 3, \ldots \right\}$.*

1. *$\boldsymbol{L}_2^N \oslash \boldsymbol{L}_2^{N-1}$ is an **orthocomplemented lattice** (Definition 1.72 page 20)* *and*
2. *The structure $\mathbb{P} \triangleq (X, \vee, \wedge; \subseteq)$ is a **primorial lattice** (Definition 5.1 page 50).*

✎PROOF:

(1) Proof that $\boldsymbol{L}_2^N \oslash \boldsymbol{L}_2^{N-1}$ is an **orthocomplemented lattice**:

- (a) $\boldsymbol{L}_2^N$ is a *Boolean lattice* by definition.
- (b) $\boldsymbol{L}_2^{N-1}$ is also a *Boolean lattice* (Definition 5.7 page 52).
- (c) Every lattice that is *Boolean* is also *orthocomplemented* (Proposition 1.80 page 23).
- (d) By definition of $\boldsymbol{L}_2^N \oslash \boldsymbol{L}_2^{N-1}$, *orthocomplemented pairs* are removed from $\boldsymbol{L}_2^N$ and the orthocomplemented pair $\{0, 1\}$ is put back in.
- (e) What remains in $\boldsymbol{L}_2^N \oslash \boldsymbol{L}_2^{N-1}$ is a set of *orthocomplemented pair*s, ordered with the same ordering relation $\leq$ that orders $\boldsymbol{L}_2^N$.
- (f) All remaining *orthocomplemented pair*s are still *involutory*: $x = x^{\perp\perp} \quad \forall x \in X$
- (g) All remaining *orthocomplemented pair*s are still *antitone* because the *ordering relation* $\leq$ in $\boldsymbol{L}_2^N$ and $\boldsymbol{L}_2^N \oslash \boldsymbol{L}_2^{N-1}$ is the same.
- (h) All remaining *orthocomplemented pair*s still have the *non-contradiction* property because suppose that in $\boldsymbol{L}_2^N \oslash \boldsymbol{L}_2^{N-1}$, there is an element $x$ such that $x \wedge x^{\perp} = m \neq 0$. Then in $\boldsymbol{L}_2^N$, it would also be true that $x \wedge x^{\perp} \neq 0$. This cannot be true (is a contradiction); so therefore for all $x$ in $\boldsymbol{L}_2^N \oslash \boldsymbol{L}_2^{N-1}$, $x \wedge x^{\perp} = 0$ (*non-contradiction* property).
- (i) So $\boldsymbol{L}_2^N \oslash \boldsymbol{L}_2^{N-1}$ is an *orthocomplemented lattice* (Definition 1.72 page 20).

(2) Proof that $\left( X \triangleq \left\{ \boldsymbol{L}_2^n \,|\, n = 1, 2, \ldots \right\} \cup \left\{ \boldsymbol{L}_2^n \oslash \boldsymbol{L}_2^{n-1} \,|\, n = 2, 3, \ldots \right\}, \subseteq \right)$ is a **primorial lattice**: This follows directly from the construction of the *bounded lattice difference* (Definition 5.13 page 54) and the definition of *primorial lattice*s (Definition 5.1 page 50).

✏

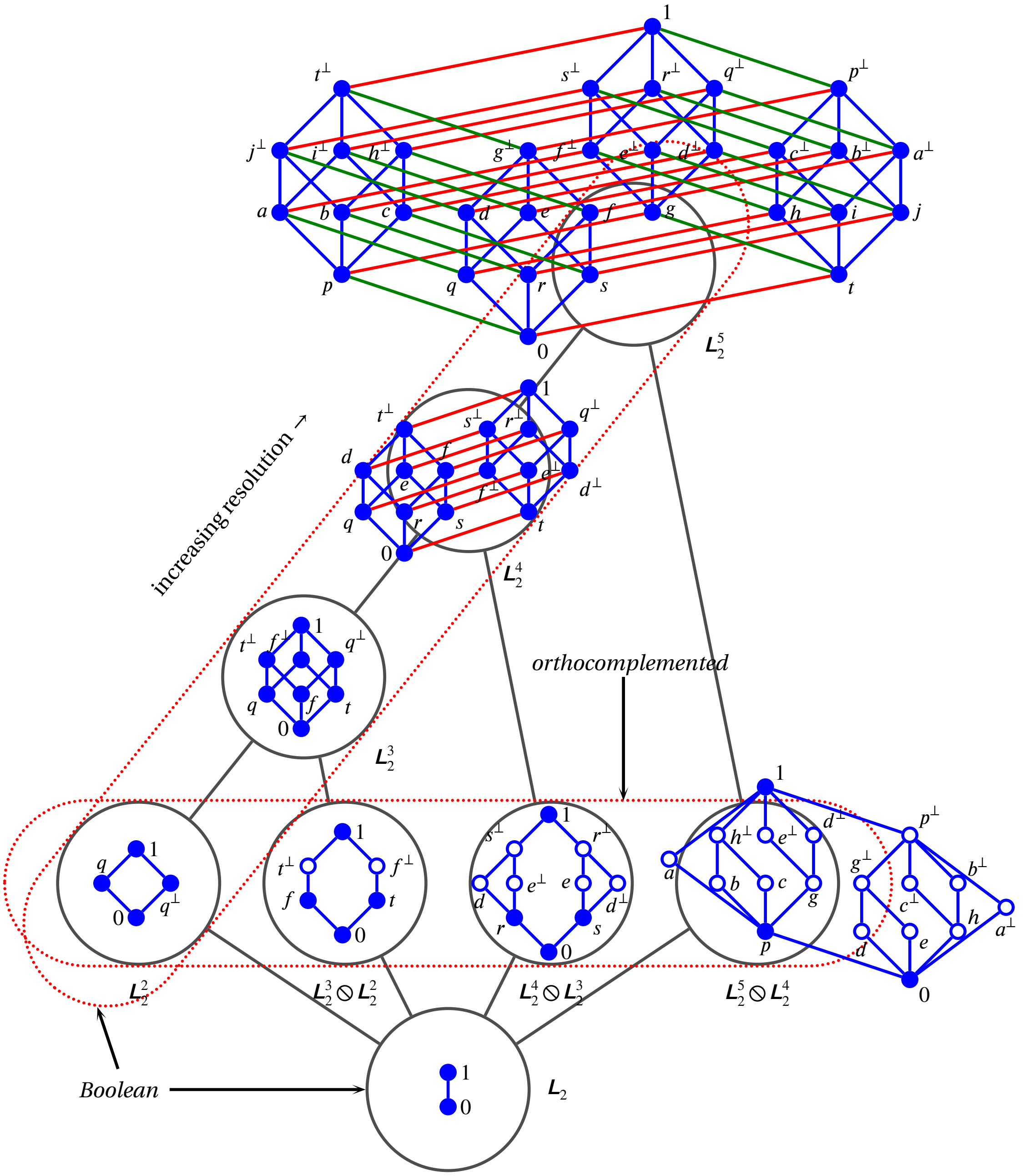

Figure 15: *a primorial lattice generated by* $\boldsymbol{L}_2^5$

**Definition 5.17** Let $\boldsymbol{L}_2^N$ be a $2^N$ element *Boolean lattice* (Definition 1.69 page 18). The lattice $\mathbb{P}$ as described in Theorem 5.16 is a **primorial lattice generated by $\boldsymbol{L}_2^N$**.

**Example 5.18** Figure 15 (page 57) illustrates a *primorial lattice generated by* $\boldsymbol{L}_2^5$.

## 5.4 Projections on primorial lattices

This section introduces three lattice projections. When performing analysis in a *primorial lattice* (Definition 5.1 page 50), it is necessary to project a point that exists in a lattice of "high resolution" onto a lattice $\boldsymbol{L}$ of lower resolution that may or may not contain this point. The three projections introduced here are the

1. *zero primorial projection* (Definition 5.19 page 58) which assigns to 0 any point that does not exist in $\boldsymbol{L}$
2. *Sasaki primorial projection* (Definition 5.20 page 58) which assigns a projection value using the *Sasaki projection* (Definition 2.22 page 31)
3. *metric primorial projection* (Definition 5.22 page 59) which assigns a projection value based on a *lattice metric* (Definition 2.7 page 27).

**Definition 5.19** Let $\mathbb{P}$ be a *primorial lattice* (Definition 5.17 page 58) generated by a *Boolean lattice* $\boldsymbol{L}_2^N$ (Definition 1.69 page 18). Let $\boldsymbol{L} \triangleq (Y, \vee, \wedge, 0, 1\,;\, \leq)$ be a lattice in $\mathbb{P}$. Let $\mathbb{x} \triangleq (\!(x_n)\!)$ be a *sequence* over the set $X$. The **zero primorial projection $\Phi_{\boldsymbol{L}}(x)$** of $x$ onto $\boldsymbol{L}$ is defined as

$$\Phi^z_{\boldsymbol{L}}(x) \triangleq \bigvee_{\boldsymbol{L}} \left[\{x, 0\} \cap Y\right] \qquad \forall x \in X$$

The **zero primorial projection $\Phi^z_{\boldsymbol{L}}(\mathbb{x})$** of $\mathbb{x}$ onto $\boldsymbol{L}$ is defined as

$$\Phi^z_{\boldsymbol{L}}(\mathbb{x}) \triangleq (\!(y_n)\!) \text{ where } y_n \triangleq \Phi^z_{\boldsymbol{L}}(x_n) \qquad \forall x_n \in (\!(x_n)\!),\, y_n \in (\!(y_n)\!).$$

**Definition 5.20** Let $\mathbb{P}$ and $\mathbb{x}$ be defined as in Definition 5.19 (page 58). Let $\mathbb{P}$ be a *primorial lattice* (Definition 5.17 page 58) generated by a *Boolean lattice* $\boldsymbol{L}_2^N$ (Definition 1.69 page 18). Let $\boldsymbol{L} \triangleq (Y, \vee, \wedge, 0, 1\,;\, \leq)$ be a lattice in $\mathbb{P}$. Let $\mathbb{x} \triangleq (\!(x_n)\!)$ be a *sequence* over the set $X$. The **Sasaki primorial projection $\Phi^s_{\boldsymbol{L}}(x)$** of $x$ onto $\boldsymbol{L}$ is defined as

$$\Phi^s_{\boldsymbol{L}}(x) \triangleq \bigvee_{\boldsymbol{L}} \left[\left\{\phi_y(x) \,|\, y \in Y\right\} \cap Y\right] \qquad \forall x \in \boldsymbol{L}$$

where $\phi_y(x)$ is the *Sasaki projection* of $x$ onto $y$ (Definition 2.22 page 31) in the smallest *Boolean lattice* $\boldsymbol{L}_2^M$ that contains both $x$ and $\boldsymbol{L}$. The **Sasaki primorial projection $\Phi^s_{\boldsymbol{L}}(\mathbb{x})$** of $\mathbb{x}$ onto $\boldsymbol{L}$ is defined as

$$\Phi^s_{\boldsymbol{L}}(\mathbb{x}) \triangleq (\!(y_n)\!) \text{ where } y_n \triangleq \Phi^s_{\boldsymbol{L}}(x_n) \qquad \forall x_n \in (\!(x_n)\!).$$

The *Sasaki primorial projection* yields a kind of *maxmini* (Theorem 1.35 page 11) result:

**Proposition 5.21** *Let $\Phi_{L}(x)$ be the* SASAKI PRIMORIAL PROJECTION *of $x$ onto $L$ in a* PRIMORIAL LATTICE $\mathbb{P}$.

$$\Phi^s_{L}(x) = \bigvee_{L} [\{x \wedge y \,|\, y \in Y\} \cap Y] \qquad \forall x \in X$$

✎PROOF:

$$\begin{aligned}
\Phi^s_{L}(x) &\triangleq \bigvee \left[\left\{\phi_y(x) \,|\, y \in Y\right\} \cap Y\right] && \text{by def. of } \textit{Sasaki primorial projection} \text{ (Definition 5.20 page 58)}\\
&\triangleq \bigvee \left[\left\{(x \vee y^{\perp}) \wedge y \,|\, y \in Y\right\} \cap Y\right] && \text{by definition of } \textit{Sasaki projection} \text{ (Definition 2.22 page 31)}\\
&= \bigvee \left[\left\{(x \wedge y) \vee (y^{\perp} \wedge y) \,|\, y \in Y\right\} \cap Y\right] && \text{by } \textit{distributive} \text{ prop. (Theorem 1.70 page 19)}\\
&= \bigvee [\{(x \wedge y) \vee (0) \,|\, y \in Y\} \cap Y] && \text{by } \textit{noncontradiction} \text{ property (Theorem 1.70 page 19)}\\
&= \bigvee [\{x \wedge y \,|\, y \in Y\} \cap Y] && \text{by } \textit{bounded} \text{ property (Theorem 1.70 page 19)}
\end{aligned}$$

✎

**Definition 5.22** Let $\mathbb{P}$ and $\mathbb{x}$ be defined as in Definition 5.19.
The *metric primorial projection* $\Phi^m_{L}(x)$ of $x$ onto $L$ is defined as

$$\Phi^m_{L}(x) \triangleq \bigwedge_{L} \left[\overline{\mathsf{B}}(x, r) \cap Y\right] \qquad \text{where}$$

1. $\overline{\mathsf{B}}(x, r)$ is the *closed ball* in $\left(L^M_2, \mathsf{d}\right)$ with the smallest radius $r$ that contains $x$ *and*
2. $\left(L^M_2, \mathsf{d}\right)$ is a *metric lattice* (Definition 2.7 page 27) *and*
3. $L^M_2$ is the smallest *Boolean lattice* (Definition 1.69 page 18) containing $x$ *and*
4. the *valuation* function defining $\mathsf{d}$ is the *height* function on $L^M_2$.

The *metric primorial projection* $\Phi_{L}(\mathbb{x})$ of $\mathbb{x}$ onto $L$ is defined as

$$\Phi_{L}(\mathbb{x}) \triangleq \langle\!\langle y_n \rangle\!\rangle \text{ such that } y_n \triangleq \Phi_{L}(x_n).$$

**Example 5.23** Here are examples of the *primorial projections* $\Phi^z_{O_6}(x)$ (Definition 5.19 page 58), $\Phi^s_{O_6}(x)$ (Definition 5.20 page 58), and $\Phi^m_{O_6}(x)$ (Definition 5.22 page 59) in the *primorial lattice* (Definition 5.1 page 50) generated by the *Boolean lattice* (Definition 1.69 page 18) $L^5_2 \triangleq (X, \vee, \wedge, 0, 1; \le)$ as illustrated in Figure 15 page 57 onto the lattice $O_6 \triangleq L^3_2 \oslash L^2_2 \triangleq (Y, \vee, \wedge, 0, 1; \le)$.

| projection | $x$ in $O_6 \triangleq L^3_2 \oslash L^2_2$ | | | | | | $x$ in $L^3_2$ | | $x$ in $L^4_2$ | | | | $x$ in $L^5_2$ | | | | | |
|---|---|---|---|---|---|---|---|---|---|---|---|---|---|---|---|---|---|---|
| $x =$ | $0$ | $f$ | $t$ | $t^{\perp}$ | $f^{\perp}$ | $1$ | $q$ | $q^{\perp}$ | $r$ | $r^{\perp}$ | $s$ | $s^{\perp}$ | $g$ | $g^{\perp}$ | $p$ | $p^{\perp}$ | $d$ | $d^{\perp}$ |
| $\Phi^z_{O_6}(x) =$ | $0$ | $f$ | $t$ | $t^{\perp}$ | $f^{\perp}$ | $1$ | $0$ | $0$ | $0$ | $0$ | $0$ | $0$ | $0$ | $0$ | $0$ | $0$ | $0$ | $0$ |
| $\Phi^s_{O_6}(x) =$ | $0$ | $f$ | $t$ | $t^{\perp}$ | $f^{\perp}$ | $1$ | $0$ | $1$ | $0$ | $f^{\perp}$ | $0$ | $f^{\perp}$ | $t$ | $f$ | $0$ | $1$ | $0$ | $t$ |
| $\Phi^m_{O_6}(x) =$ | $0$ | $f$ | $t$ | $t^{\perp}$ | $f^{\perp}$ | $1$ | $0$ | $0$ | $0$ | $f^{\perp}$ | $0$ | $f^{\perp}$ | $t$ | $f$ | $0$ | $1$ | $0$ | $t$ |

✎PROOF:

(1) Proof for zero primorial projection values:

$$
\begin{aligned}
\Phi^z_{O_6}(0) &= \bigvee\left[(\{0\}\cup\{0\})\cap\left\{0,f,t,t^\perp,f^\perp,1\right\}\right] &&= \bigvee[\{0\}] &&= 0\\
\Phi^z_{O_6}(f) &= \bigvee\left[(\{f\}\cup\{0\})\cap\left\{0,f,t,t^\perp,f^\perp,1\right\}\right] &&= \bigvee[\{0,f\}] &&= f\\
\Phi^z_{O_6}(t) &= \bigvee\left[(\{t\}\cup\{0\})\cap\left\{0,f,t,t^\perp,f^\perp,1\right\}\right] &&= \bigvee[\{0,t\}] &&= t\\
\Phi^z_{O_6}(t^\perp) &= \bigvee\left[\left(\left\{t^\perp\right\}\cup\{0\}\right)\cap\left\{0,f,t,t^\perp,f^\perp,1\right\}\right] &&= \bigvee\left[\left\{0,t^\perp\right\}\right] &&= t^\perp\\
\Phi^z_{O_6}(f^\perp) &= \bigvee\left[\left(\left\{f^\perp\right\}\cup\{0\}\right)\cap\left\{0,f,t,t^\perp,f^\perp,1\right\}\right] &&= \bigvee\left[\left\{0,f^\perp\right\}\right] &&= f^\perp\\
\Phi^z_{O_6}(1) &= \bigvee\left[(\{1\}\cup\{0\})\cap\left\{0,f,t,t^\perp,f^\perp,1\right\}\right] &&= \bigvee[\{1,0\}] &&= 1\\
\Phi^z_{O_6}(q) &= \bigvee\left[(\{q\}\cup\{0\})\cap\left\{0,f,t,t^\perp,f^\perp,1\right\}\right] &&= \bigvee[\{0\}] &&= 0\\
\Phi^z_{O_6}(q^\perp) &= \bigvee\left[\left(\left\{q^\perp\right\}\cup\{0\}\right)\cap\left\{0,f,t,t^\perp,q^\perp,1\right\}\right] &&= \bigvee\left[\left\{0,q^\perp\right\}\right] &&= 0\\
\Phi^z_{O_6}(r) &= \bigvee\left[(\{r\}\cup\{0\})\cap\left\{0,f,t,t^\perp,f^\perp,1\right\}\right] &&= \bigvee[\{0\}] &&= 0\\
\Phi^z_{O_6}(r^\perp) &= \bigvee\left[\left(\left\{r^\perp\right\}\cup\{0\}\right)\cap\left\{0,f,t,t^\perp,r^\perp,1\right\}\right] &&= \bigvee\left[\left\{0,r^\perp\right\}\right] &&= 0\\
\Phi^z_{O_6}(s) &= \bigvee\left[(\{s\}\cup\{0\})\cap\left\{0,f,t,t^\perp,f^\perp,1\right\}\right] &&= \bigvee[\{0\}] &&= 0\\
\Phi^z_{O_6}(s^\perp) &= \bigvee\left[\left(\left\{s^\perp\right\}\cup\{0\}\right)\cap\left\{0,f,t,t^\perp,r^\perp,1\right\}\right] &&= \bigvee[\{0\}] &&= 0\\
\Phi^z_{O_6}(g) &= \bigvee\left[(\{g\}\cup\{0\})\cap\left\{0,f,t,t^\perp,f^\perp,1\right\}\right] &&= \bigvee[\{0\}] &&= 0\\
\Phi^z_{O_6}(g^\perp) &= \bigvee\left[\left(\left\{g^\perp\right\}\cup\{0\}\right)\cap\left\{0,f,t,t^\perp,r^\perp,1\right\}\right] &&= \bigvee[\{0\}] &&= 0\\
\Phi^z_{O_6}(p) &= \bigvee\left[(\{p\}\cup\{0\})\cap\left\{0,f,t,t^\perp,f^\perp,1\right\}\right] &&= \bigvee[\{0\}] &&= 0\\
\Phi^z_{O_6}(p^\perp) &= \bigvee\left[\left(\left\{p^\perp\right\}\cup\{0\}\right)\cap\left\{0,f,t,t^\perp,r^\perp,1\right\}\right] &&= \bigvee[\{0\}] &&= 0\\
\Phi^z_{O_6}(d) &= \bigvee\left[(\{d\}\cup\{0\})\cap\left\{0,f,t,t^\perp,f^\perp,1\right\}\right] &&= \bigvee[\{0\}] &&= 0\\
\Phi^z_{O_6}(d^\perp) &= \bigvee\left[\left(\left\{d^\perp\right\}\cup\{0\}\right)\cap\left\{0,f,t,t^\perp,r^\perp,1\right\}\right] &&= \bigvee[\{0\}] &&= 0
\end{aligned}
$$

(2) Proof for Sasaki primorial projection (Definition 5.20 page 58):

$$
\begin{aligned}
\Phi^s_{O_6}(0) &= \bigvee[\{0\wedge y\,|\,y\in Y\}\cap Y] &&= \bigvee[\{0,0,0,0,0,0\}\cap Y] &&= \bigvee\{0\} &&= 0\\
\Phi^s_{O_6}(f) &= \bigvee[\{f\wedge y\,|\,y\in Y\}\cap Y] &&= \bigvee[\{0,f,0,f,0,f\}\cap Y] &&= \bigvee\{0,f\} &&= f\\
\Phi^s_{O_6}(t) &= \bigvee[\{t\wedge y\,|\,y\in Y\}\cap Y] &&= \bigvee[\{0,0,t,0,t,t\}\cap Y] &&= \bigvee\{0,t\} &&= t\\
\Phi^s_{O_6}(t^\perp) &= \bigvee\left[\left\{t^\perp\wedge y\,|\,y\in Y\right\}\cap Y\right] &&= \bigvee\left[\left\{0,f,0,t^\perp,q,t^\perp\right\}\cap Y\right] &&= \bigvee\left\{0,f,t^\perp\right\} &&= t^\perp\\
\Phi^s_{O_6}(f^\perp) &= \bigvee\left[\left\{f^\perp\wedge y\,|\,y\in Y\right\}\cap Y\right] &&= \bigvee\left[\left\{0,0,t,q,f^\perp,f^\perp\right\}\cap Y\right] &&= \bigvee\left\{0,t,f^\perp\right\} &&= f^\perp\\
\Phi^s_{O_6}(1) &= \bigvee[\{1\wedge y\,|\,y\in Y\}\cap Y] &&= \bigvee\left[\left\{0,f,t,f^\perp,t^\perp,1\right\}\cap Y\right] &&= \bigvee Y &&= 1\\
\Phi^s_{O_6}(q) &= \bigvee[\{q\wedge y\,|\,y\in Y\}\cap Y] &&= \bigvee[\{0,0,0,q,0,q\}\cap Y] &&= \bigvee\{0\} &&= 0\\
\Phi^s_{O_6}(q^\perp) &= \bigvee\left[\left\{q^\perp\wedge y\,|\,y\in Y\right\}\cap Y\right] &&= \bigvee\left[\left\{0,f,t,f,t,q^\perp\right\}\cap Y\right] &&= \bigvee\{0,f,t\} &&= 1
\end{aligned}
$$

$$\begin{aligned}
\Phi^s_{O_6}(r) &= \bigvee\left[\{r\wedge y \mid y\in Y\}\cap Y\right] &&= \bigvee\left[\{0,r,0,r,0,r\}\cap Y\right] &&= \bigvee\{0\} &&= 0\\
\Phi^s_{O_6}(r^\perp) &= \bigvee\left[\left\{r^\perp\wedge y \mid y\in Y\right\}\cap Y\right] &&= \bigvee\left[\left\{0,s,t,e,f^\perp,r^\perp\right\}\cap Y\right] &&= \bigvee\left\{0,t,f^\perp\right\} &&= f^\perp\\
\Phi^s_{O_6}(s) &= \bigvee\left[\{s\wedge y \mid y\in Y\}\cap Y\right] &&= \bigvee\left[\{0,s,0,s,0,s\}\cap Y\right] &&= \bigvee\{0\} &&= 0\\
\Phi^s_{O_6}(s^\perp) &= \bigvee\left[\left\{s^\perp\wedge y \mid y\in Y\right\}\cap Y\right] &&= \bigvee\left[\left\{0,0,t,d,f^\perp,s^\perp\right\}\cap Y\right] &&= \bigvee\left\{0,t,f^\perp\right\} &&= f^\perp\\
\Phi^s_{O_6}(g) &= \bigvee\left[\{g\wedge y \mid y\in Y\}\cap Y\right] &&= \bigvee\left[\{0,0,t,p,g,g\}\cap Y\right] &&= \bigvee\{0,t\} &&= t\\
\Phi^s_{O_6}(g^\perp) &= \bigvee\left[\left\{g^\perp\wedge y \mid y\in Y\right\}\cap Y\right] &&= \bigvee\left[\left\{0,f,0,g^\perp,0,g^\perp\right\}\cap Y\right] &&= \bigvee\{0,f\} &&= f\\
\Phi^s_{O_6}(p) &= \bigvee\left[\{p\wedge y \mid y\in Y\}\cap Y\right] &&= \bigvee\left[\{0,0,0,p,p,p\}\cap Y\right] &&= \bigvee\{0\} &&= 0\\
\Phi^s_{O_6}(p^\perp) &= \bigvee\left[\left\{p^\perp\wedge y \mid y\in Y\right\}\cap Y\right] &&= \bigvee\left[\left\{0,f,t,g^\perp,t,p^\perp\right\}\cap Y\right] &&= \bigvee\{0,f,t\} &&= 1\\
\Phi^s_{O_6}(d) &= \bigvee\left[\{d\wedge y \mid y\in Y\}\cap Y\right] &&= \bigvee\left[\{0,r,0,d,0,d\}\cap Y\right] &&= \bigvee\{0\} &&= 0\\
\Phi^s_{O_6}(d^\perp) &= \bigvee\left[\left\{d^\perp\wedge y \mid y\in Y\right\}\cap Y\right] &&= \bigvee\left[\left\{0,s,t,0,g,d^\perp\right\}\cap Y\right] &&= \bigvee\{0,t\} &&= t
\end{aligned}$$

(3) Proof for metric primorial projection (Definition 5.22 page 59):

$$\begin{aligned}
\Phi^m_{O_6}(0) &= \bigwedge\left[\overline{\mathsf{B}}(0,0)\cap Y\right] &&= \bigwedge\left[\{0\}\cap\left\{0,f,t,t^\perp,f^\perp,1\right\}\right] &&= \bigwedge\{0\} &&= 0\\
\Phi^m_{O_6}(f) &= \bigwedge\left[\overline{\mathsf{B}}(f,0)\cap Y\right] &&= \bigwedge\left[\{f\}\cap\left\{0,f,t,t^\perp,f^\perp,1\right\}\right] &&= \bigwedge\{f\} &&= f\\
\Phi^m_{O_6}(t) &= \bigwedge\left[\overline{\mathsf{B}}(t,0)\cap Y\right] &&= \bigwedge\left[\{t\}\cap\left\{0,f,t,t^\perp,f^\perp,1\right\}\right] &&= \bigwedge\{t\} &&= t\\
\Phi^m_{O_6}(t^\perp) &= \bigwedge\left[\overline{\mathsf{B}}\left(t^\perp,0\right)\cap Y\right] &&= \bigwedge\left[\left\{t^\perp\right\}\cap\left\{0,f,t,t^\perp,f^\perp,1\right\}\right] &&= \bigwedge\left\{t^\perp\right\} &&= t^\perp\\
\Phi^m_{O_6}(f^\perp) &= \bigwedge\left[\overline{\mathsf{B}}\left(f^\perp,0\right)\cap Y\right] &&= \bigwedge\left[\left\{f^\perp\right\}\cap\left\{0,f,t,t^\perp,f^\perp,1\right\}\right] &&= \bigwedge\left\{f^\perp\right\} &&= f^\perp\\
\Phi^m_{O_6}(1) &= \bigwedge\left[\overline{\mathsf{B}}(1,0)\cap Y\right] &&= \bigwedge\left[\{1\}\cap\left\{0,f,t,t^\perp,f^\perp,1\right\}\right] &&= \bigwedge\{1\} &&= 1\\
\Phi^m_{O_6}(q) &= \bigwedge\left[\overline{\mathsf{B}}(q,1)\cap Y\right] &&= \bigwedge\left[\left\{q,0,t^\perp\right\}\cap Y\right] &&= \bigwedge\left\{0,t^\perp\right\} &&= 0\\
\Phi^m_{O_6}(q^\perp) &= \bigwedge\left[\overline{\mathsf{B}}\left(q^\perp,1\right)\cap Y\right] &&= \bigwedge\left[\left\{q^\perp,t,1\right\}\cap Y\right] &&= \bigwedge\{t,1\} &&= t\\
\Phi^m_{O_6}(r) &= \bigwedge\left[\overline{\mathsf{B}}(r,1)\cap Y\right] &&= \bigwedge\left[\{r,0,d,f\}\cap Y\right] &&= \bigwedge\{0,f\} &&= 0\\
\Phi^m_{O_6}(r^\perp) &= \bigwedge\left[\overline{\mathsf{B}}\left(r^\perp,1\right)\cap Y\right] &&= \bigwedge\left[\left\{r^\perp,d^\perp,f^\perp,1\right\}\cap Y\right] &&= \bigwedge\left\{f^\perp,1\right\} &&= f^\perp\\
\Phi^m_{O_6}(s) &= \bigwedge\left[\overline{\mathsf{B}}(s,1)\cap Y\right] &&= \bigwedge\left[\left\{s,0,e,f,e^\perp\right\}\cap Y\right] &&= \bigwedge\{0,f\} &&= 0\\
\Phi^m_{O_6}(s^\perp) &= \bigwedge\left[\overline{\mathsf{B}}\left(s^\perp,1\right)\cap Y\right] &&= \bigwedge\left[\left\{s^\perp,e^\perp,f^\perp,d,1\right\}\cap Y\right] &&= \bigwedge\left\{f^\perp,1\right\} &&= f^\perp\\
\Phi^m_{O_6}(g) &= \bigwedge\left[\overline{\mathsf{B}}(g,1)\cap Y\right] &&= \bigwedge\left[\left\{g,p,f^\perp,e^\perp,d^\perp,t\right\}\cap Y\right] &&= \bigwedge\left\{f^\perp,t\right\} &&= t\\
\Phi^m_{O_6}(g^\perp) &= \bigwedge\left[\overline{\mathsf{B}}\left(g^\perp,1\right)\cap Y\right] &&= \bigwedge\left[\left\{g^\perp,d,e,f,p^\perp,t^\perp\right\}\cap Y\right] &&= \bigwedge\left\{f,t^\perp\right\} &&= f\\
\Phi^m_{O_6}(p) &= \bigwedge\left[\overline{\mathsf{B}}(p,1)\cap Y\right] &&= \bigwedge\left[\{p,0,p,a,b,c,g\}\cap Y\right] &&= \bigwedge\{0\} &&= 0
\end{aligned}$$

$$\Phi^m_{O_6}(p^\perp) = \bigwedge \left[\overline{\mathsf{B}}\left(p^\perp, 1\right) \cap Y\right] \quad = \bigwedge \left[\left\{p^\perp, a^\perp, b^\perp, c^\perp, g^\perp, 1\right\} \cap Y\right] \quad = \bigwedge \{1\} \qquad = 1$$

$$\begin{aligned}
\Phi^m_{O_6}(d) &= \bigwedge \left[\overline{\mathsf{B}}(d,2) \cap \left\{0, f, t, t^\perp, f^\perp, 1\right\}\right] \\
&= \bigwedge \left[\left\{0, a, b, d, e, f, h, i, q, r, c^\perp, g^\perp, j^\perp, p^\perp, s^\perp, t^\perp\right\} \cap \left\{0, f, t, t^\perp, f^\perp, 1\right\}\right] \\
&= \bigwedge \left\{0, f, t^\perp\right\} \\
&= 0 \\
\Phi^m_{O_6}(d^\perp) &= \bigwedge \left[\overline{\mathsf{B}}\left(d^\perp,2\right) \cap \left\{0, f, t, t^\perp, f^\perp, 1\right\}\right] \\
&= \bigwedge \left[\left\{c, g, j, p, s, t, a^\perp, b^\perp, d^\perp, e^\perp, f^\perp, h^\perp, i^\perp, q^\perp, r^\perp, 1\right\} \cap \left\{0, f, t, t^\perp, f^\perp, 1\right\}\right] \\
&= \bigwedge \left\{t, f^\perp, 1\right\} \\
&= t
\end{aligned}$$

✏

## 5.5 A generalized probability function

This paper introduces a new definition for a lattice-valued probability function (next).

**Definition 5.24** Let $\boldsymbol{L} \triangleq (X, \vee, \wedge, \neg, 0, 1; \leq)$ be a *lattice with negation* (Definition 2.16 page 30). Let Ⓓ be the *distributivity* relation (Definition 1.52 page 15). A function $\mathsf{p}$ in $\mathbb{R}^{\boldsymbol{X}}$ is a **probability** on $\boldsymbol{L}$ if

1. $\mathsf{p}(0) = 0$ (*nondegenerate*) and
2. $\mathsf{p}(1) = 1$ (*normalized*) and
3. $x \leq y \implies \mathsf{p}(x) \leq \mathsf{p}(y) \quad \forall x,y \in X$ (*monotone*) and
4. $\left\{\begin{array}{ll} x \wedge y = 0 & \text{and} \\ (z, x, y) \in Ⓓ & \forall z \in X \end{array}\right\} \implies \mathsf{p}(x \vee y) = \mathsf{p}(x) + \mathsf{p}(y) \quad \forall x,y \in X$ (*additive*).

If $\mathsf{p}$ is a *probability* on a *lattice with negation* $\boldsymbol{L}$, then $(\boldsymbol{L}, \mathsf{p})$ is a **probability space**.

**Remark 5.25** Definition 5.24 page 62 (previous) is not any standard definition of the *probability function*. On a *Boolean lattice*, the **measure-theoretic probability** function, due to A. N. Kolmogorov, is defined as[119]

(1). $\mathsf{p}(1) = 1$ (*normalized*) and
(2). $\mathsf{p}(x) \geq 0 \quad \forall x \in X$ (*nonnegative*) and
(3). $\displaystyle\bigwedge_{n=1}^{\infty} x_n = 0 \implies \mathsf{p}\left(\bigvee_{n=1}^{\infty} x_n\right) = \sum_{n=1}^{\infty} \mathsf{p}(x_n) \quad \forall x_n \in X$ (*σ-additive*) .

[119] [13], pages 22–23, ⟨Probability Measures⟩, [103], [102], page 16, ⟨*field of probability*⟩, [137], pages 8–9, ⟨Definition 2.3(13)⟩, [99], page 27

The advantage of this definition is that p is a *measure*, and hence all the power of measure theory is subsequently at one's disposal in using p. However, it has often been argued that the requirement of *σ-additivity* is unnecessary for a probability function. Even as early as 1930, de Finetti argued against it, in what became a kind of polite running debate with Fréchet.[120] In fact, Kolmogorov himself provided some argument against *σ-additivity* when referring to the closely related *Axiom of Continuity* saying, "Since the new axiom is essential for infinite fields of probability only, it is almost impossible to elucidate its empirical meaning…For, in describing any observable random process we can obtain only finite fields of probability.…" But in its support he added, "This limitation has been found expedient in researches of the most diverse sort."[121]

There are several other definitions of probability that only require *additivity* rather than *σ-additivity*. On a *Boolean lattice*, the **traditional probability** function is defined as[122]

(1). $\mathsf{p}(1) = 1$ (*normalized*) and

(2). $\mathsf{p}(x) \geq 0 \quad \forall x \in X$ (*nonnegative*) and

(3). $x \wedge y = 0 \implies \mathsf{p}(x \vee y) = \mathsf{p}(x) + \mathsf{p}(y) \quad \forall x,y \in X$ (*additive*) .

This definition implies (on a *Boolean lattice*) that

(a). $\mathsf{p}(0) = 0$ (*nondegenerate*) and

(b). $\mathsf{p}(x) \leq 1 \quad \forall x \in X$ (*upper bounded*) and

(c). $\mathsf{p}(x) = 1 - \mathsf{p}(\neg x) \quad \forall x \in X$ and

(d). $\mathsf{p}(x \vee y) \leq \mathsf{p}(x) + \mathsf{p}(y) \quad \forall x,y \in X$ (*subadditive*) and

(e). $\mathsf{p}(x \vee y) = \mathsf{p}(x) + \mathsf{p}(y) - \mathsf{p}(x \wedge y) \quad \forall x,y \in X$ and

(f). $x \leq y \implies \mathsf{p}(x) \leq \mathsf{p}(y) \quad \forall x,y \in X$ (*monotone*) .

On a *distributive pseudocomplemented lattice*, the **generalized probability** function has been defined as[123]

(1). $\mathsf{p}(0) = 0$ (*nondegenerate*) and

(2). $\mathsf{p}(1) = 1$ (*normalized*) and

(3). $0 \leq \mathsf{p}(1) \leq 1$ and

(4). $\mathsf{p}(x \vee y) = \mathsf{p}(x) + \mathsf{p}(y) - \mathsf{p}(x \wedge y) \quad \forall x,y \in X$ .

On an *orthomodular lattice*, or a *finite modular lattice*, the **quantum probability** function is defined as[124]

(1). $\mathsf{p}(0) = 0$ (*nondegenerate*) and

(2). $\mathsf{p}(1) = 1$ (*normalized*) and

(3). $x \perp y \implies \mathsf{p}(x \vee y) = \mathsf{p}(x) + \mathsf{p}(y) \quad \forall x,y \in X$ (*additive*) .

However, for lattices that are not *distributive*, *modular*, or *orthomodular*, none of these definitions work out so well. Take for example the $O_6$ *lattice* with the "very reasonable" probability function given in Example 5.31 (page 66). This probability space $(O_6, \mathsf{p})$ fails to be any of the 4 probability functions defined in this Remark. It fails to be a *measure-theoretic* or *traditional probability* function because

$$a \wedge b = 0 \qquad \text{but} \qquad \mathsf{p}(a \vee b) = \mathsf{p}(1) = 1 \neq \tfrac{1}{3} + \tfrac{1}{2} = \mathsf{p}(a) + \mathsf{p}(b) \; .$$

---

[120] [60], [65], [59], [66], [58], [28], pages 258–260
[121] [102], page 15
[122] [138], pages 21–22, [102], page 2, ⟨§1. Axioms I–V⟩
[123] [129], page 118, [128]
[124] [74], page 126, ⟨DEFINITIONS⟩, [129], page 118

It fails to be a *generalized probability* function because

$$\mathsf{p}(a \vee b) = \mathsf{p}(1) = 1 \neq \tfrac{1}{3} + \tfrac{1}{2} - 0 = \mathsf{p}(a) + \mathsf{p}(b) - \mathsf{p}(0) = \mathsf{p}(a) + \mathsf{p}(b) - \mathsf{p}(a \wedge b) \; .$$

It fails to be an *quantum probability* function because

$$a \perp b = 0 \qquad \text{but} \qquad \mathsf{p}(a \vee b) = \mathsf{p}(1) = 1 \neq \tfrac{1}{3} + \tfrac{1}{2} = \mathsf{p}(a) + \mathsf{p}(b) \; .$$

In each of these cases, the function $\mathsf{p}$ fails to be *additive*. The solution of Definition 5.24 (page 62) is simply to "switch off" *additivity* when the lattice is not *distributive*. This method is a little "crude", but at least it allows us to define probability on a very wide class of lattices, while retaining compatibility with the *Boolean* case (Proposition 5.26 page 64, Proposition 5.27 page 64, Proposition 5.28 page 65).

**Proposition 5.26** [125] *Let* $(\boldsymbol{L}, \mathsf{p})$ *be a* PROBABILITY SPACE *(Definition 5.24 page 62).*

$$0 \;\leq\; \mathsf{p}(x) \;\leq\; 1 \qquad \forall x \in X$$

✎PROOF:

$$\begin{aligned} 0 &= \mathsf{p}(0) && \text{by previous result} \\ &\leq \mathsf{p}(x) && \text{because } 0 \leq x \text{ and } \textit{monotone} \text{ property (Definition 5.24 page 62)} \\ \mathsf{p}(x) &\leq \mathsf{p}(1) && \text{because } x \leq 1 \text{ and } \textit{monotone} \text{ property (Definition 5.24 page 62)} \\ &= 1 && \text{by property of } \mathsf{p} \text{ (Definition 5.24 page 62)} \end{aligned}$$

✏

**Proposition 5.27** [126] *Let* $(\boldsymbol{L}, \mathsf{p})$ *be a* PROBABILITY SPACE *(Definition 5.24 page 62).*

$$\left\{ \begin{array}{l} \boldsymbol{L} \textit{ is} \\ \text{ORTHOCOMPLEMENTED} \end{array} \right\} \quad \Longrightarrow \quad \left\{ \mathsf{p}(x) \;=\; 1 - \mathsf{p}(\neg x) \quad \forall x \in X \right\}$$

✎PROOF:

$$\begin{aligned} 1 - \mathsf{p}(\neg x) &= \mathsf{p}(1) - \mathsf{p}(\neg x) && \text{by Definition 5.24 page 62} \\ &= \mathsf{p}(x \vee \neg x) - \mathsf{p}(\neg x) && \text{by } \textit{excluded middle} \text{ property of } \textit{ortho negation} \text{ (Definition 2.14 page 29)} \\ &= \mathsf{p}(x) + \mathsf{p}(\neg x) - \mathsf{p}(\neg x) && \text{because } (x)(\neg x) = 0 \text{ and } \textit{additive} \text{ property (Definition 5.24 page 62)} \\ &= \mathsf{p}(x) \end{aligned}$$

✏

---

125 📘 [138], page 21, ⟨(2-11)⟩
126 📘 [138], page 21, ⟨(2-12)⟩

**Proposition 5.28** [127] *Let* $(\boldsymbol{L}, \mathsf{p})$ *be a* PROBABILITY SPACE *(Definition 5.24 page 62)*.

$$\left\{\begin{array}{l}\boldsymbol{L} \text{ is} \\ \text{BOOLEAN}\end{array}\right\} \implies \left\{\begin{array}{llllll} 1. & \mathsf{p}(x \vee y) & = & \mathsf{p}(x) + \mathsf{p}(y) - \mathsf{p}(x \wedge y) & \forall x,y \in X & \text{and} \\ 2. & \mathsf{p}(x \vee y) & \leq & \mathsf{p}(x) + \mathsf{p}(y) & \forall x,y \in X & \text{(BOOLE'S INEQUALITY)} \end{array}\right.$$

✎PROOF:

(1) lemma: Proof that $\mathsf{p}((\neg x) \wedge y) = \mathsf{p}(y) - \mathsf{p}(x \wedge y)$:

$$\begin{aligned}
\mathsf{p}(y) - \mathsf{p}(xy) &= \mathsf{p}(1 \wedge y) - \mathsf{p}(xy) && \text{by definition of 1 and } \wedge \text{ (Definition 1.28 page 9)} \\
&= \mathsf{p}[(x \vee \neg x)y] - \mathsf{p}(xy) && \text{by } \textit{excluded middle} \text{ property of } \textit{Boolean lattices} \\
&= \mathsf{p}(xy \vee \neg xy) - \mathsf{p}(xy) && \text{by } \textit{distributive} \text{ property of } \textit{Boolean lattices} \\
&= \mathsf{p}(xy) + \mathsf{p}(\neg xy) - \mathsf{p}(xy) && \text{because } (xy)(\neg xy) = 0 \text{ and by } \textit{additive} \text{ property} \\
&= \mathsf{p}(\neg xy)
\end{aligned}$$

(2) Proof that $\mathsf{p}(x \vee y) = \mathsf{p}(x) + \mathsf{p}(y) - \mathsf{p}(x \wedge y)$:

$$\begin{aligned}
\mathsf{p}(x \vee y) &= \mathsf{p}(x \vee \neg xy) && \text{by property of } \textit{Boolean lattices} \\
&= \mathsf{p}(x) + \mathsf{p}(\neg xy) && \text{because } (x)(\neg xy) = 0 \text{ and by } \textit{additive} \text{ property} \\
&= \mathsf{p}(x) + \mathsf{p}(y) - \mathsf{p}(x \wedge y) && \text{by item 1 (page 65)}
\end{aligned}$$

✏

**Example 5.29** The function $\neg$ on the lattice $\boldsymbol{L}$ as illustrated to the right is a *Kleene negation* (Definition 2.14 page 29). Together with the probability function $\mathsf{p}$, also illustrated to the right, the pair $(\boldsymbol{L}, \mathsf{p})$ is a *probability space* (Definition 5.24 page 62).

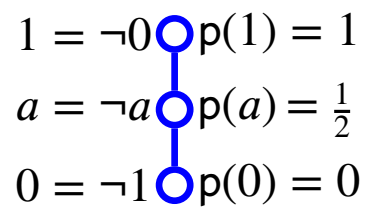

**Example 5.30** The *lattice with negation* $\boldsymbol{L}$ (Definition 2.16 page 30) illustrated to the right is a *Boolean lattice*. Together with the probability function $\mathsf{p}$, also illustrated to the right, the pair $(\boldsymbol{L}, \mathsf{p})$ is a *probability space* (Definition 5.24 page 62).

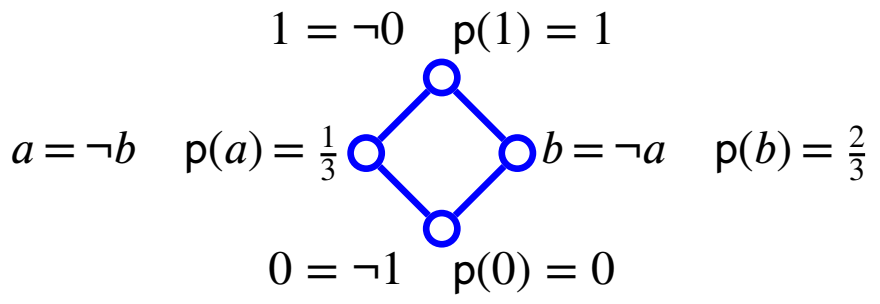

[127] [138], page 21, ⟨(2-13)⟩, [57], pages 22–23, ⟨(7.4),(7.6)⟩

**Example 5.31** The *lattice with negation* $\boldsymbol{L}$ (Definition 2.16 page 30) illustrated to the right is an *orthocomplemented* $O_6$ *lattice* (Definition 1.73 page 20). Together with the probability function $\mathsf{p}$, also illustrated to the right, the pair $(\boldsymbol{L}, \mathsf{p})$ is a *probability space* (Definition 5.24 page 62).

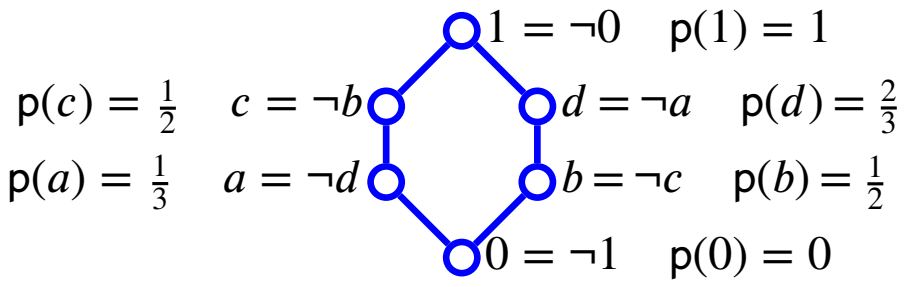

## 5.6 Applications

This section discusses some possible applications of *primorial lattice*s.

### 5.6.1 Logic analysis

Let $\boldsymbol{L}_2^N$ be a $2^N$-valued *Boolean logic* (Definition 2.27 page 33). Let $\mathbb{P}$ be the *primorial lattice generated by* $\boldsymbol{L}_2^N$ (Definition 5.17 page 58). The *sequence* of lattices $\left(\!\left(\boldsymbol{L}_2^N, \boldsymbol{L}_2^{N-1}, \ldots, \boldsymbol{L}_2^2, \boldsymbol{L}_2\right)\!\right)$ in $\mathbb{P}$ are *Boolean logic*s with decreasing "resolution" (higher values of $n$ in $\boldsymbol{L}_2^n$ correspond to greater resolution). Thus, we can reduce a very complex logic in $\boldsymbol{L}_2^N$ to a simpler lower resolution logic.

Moreover, the sequence of *ortho logics* (Definition 2.27 page 33) in $\mathbb{P}$

$$\left(\!\left(\boldsymbol{L}_2^N \oslash \boldsymbol{L}_2^{N-1},\; \boldsymbol{L}_2^{N-1} \oslash \boldsymbol{L}_2^{N-2},\; \ldots,\; \boldsymbol{L}_2^3 \oslash \boldsymbol{L}_2^2,\; \boldsymbol{L}_2\right)\!\right)$$

represents the *Boolean logic* $\boldsymbol{L}_2^N$ at $N-1$ progressively lower "frequencies". Alternatively, we could say that the *Boolean logic* at resolution $N$ is "decomposed" into (or *analyzed* by) $N-1$ *ortho logic*s. Moreover, a proposition $p$ in a higher resolution space can be projected into a lower resolution space (including the two-value classic logic space) by a *projection operator* (Section 5.4 page 58).

### 5.6.2 Fuzzy logic analysis

*Fuzzy logics* (Definition 2.27 page 33) can be constructed on *Boolean* and *orthocomplemented* lattices[128] such that together with the subset ordering relation $\subseteq$, form of a *primorial lattice* $\mathbb{P}$ (Definition 5.1 page 50). A Boolean fuzzy logic $\boldsymbol{L}_2^N$ can then be rendered at $N-1$ different "resolutions" using the Boolean lattices of $\mathbb{P}$ and analyzed at $N-1$ "frequencies" using the orthocomplemented lattices of $\mathbb{P}$, as described in Section 5.6.1 (page 66).

[128] [77], ⟨§2.2⟩

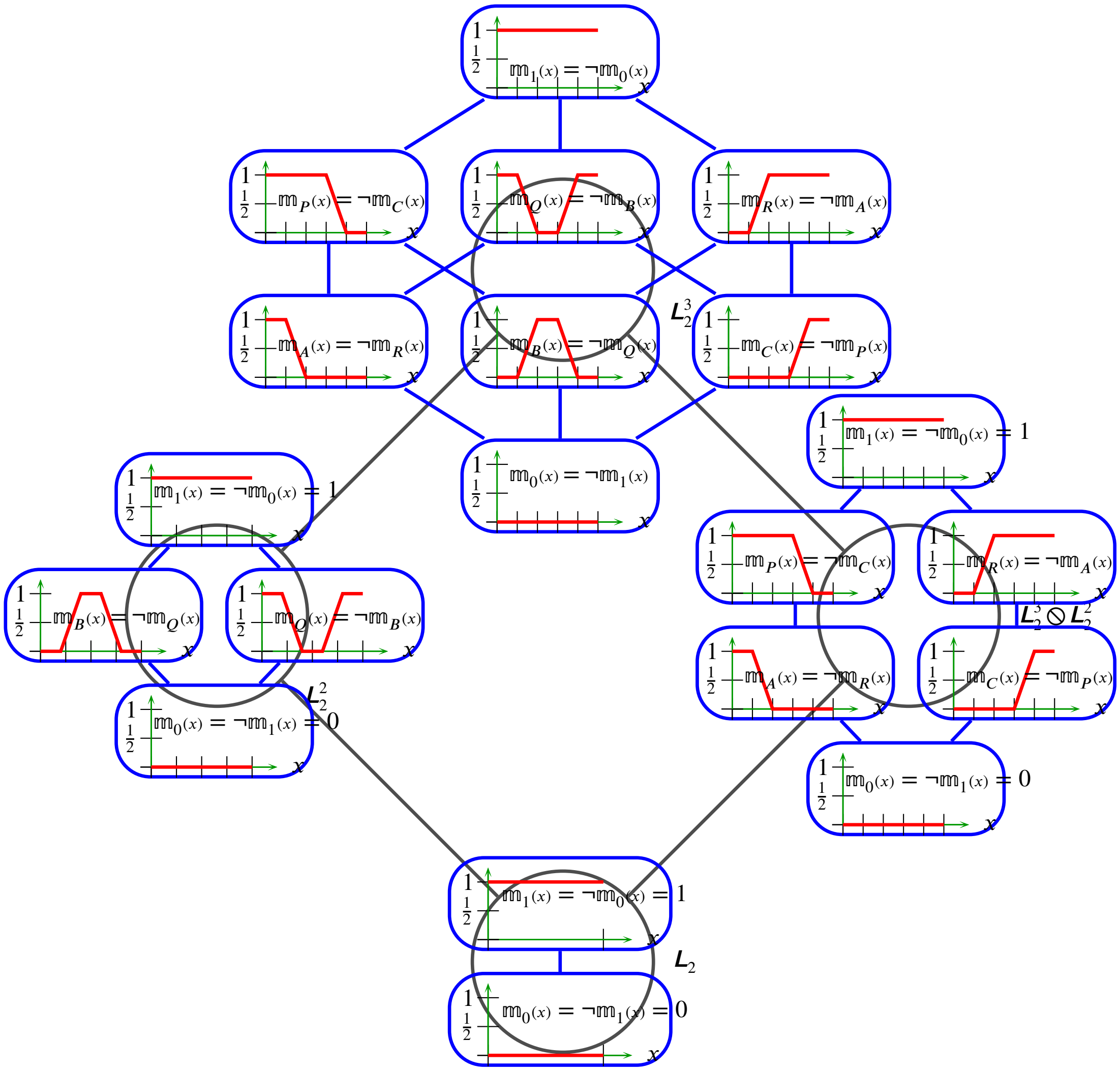

Figure 16: primorial lattice for *fuzzy subset logic* (Example 5.32 page 68)

**Example 5.32** Figure 16 (page 67) illustrates a *fuzzy subset logic*[129] on a primorial lattice. The lattice $\mathbf{L}_2^3$ contains both *monotonic* and *non-monotonic membership functions*. These are separated into lower resolution spaces $\mathbf{L}_2^2$ containing the *non-monotonic* membership functions (neglecting 1 and 0), $\mathbf{L}_2^3 \oslash \mathbf{L}_2^2$ containing the *monotonic* membership functions, and $\mathbf{L}_2$ containing crisp set logic. A *projection operator* (Section 5.4 page 58) can be used to project a membership function onto any of these spaces as perhaps called for by a given application.

### 5.6.3 Probability analysis

A *logic* is a *lattice with negation* (Definition 2.16 page 30) and with an *implication* function defined on it. A *probability* is a *lattice with negation* and with a *probability* function (Definition 5.24 page 62) defined on it.

Let $\mathbf{L}_2^N$ be the $2^N$-element Boolean lattice generated by an $N$-event *Boolean probability space* (Definition 5.24 page 62). Let $\mathbb{P}$ be the *primorial lattice* (Definition 5.1 page 50) generated by $\mathbf{L}_2^N$. Then in $\mathbb{P}$, the probability space can be rendered at progressively lower resolutions using the Boolean lattices of $\mathbb{P}$, and can be analyzed at assorted "frequencies" using the ortho-complemented lattices of $\mathbb{P}$.

**Example 5.33** A *primorial lattice* with a probability function is illustrated in Figure 17 (page 69).

### 5.6.4 Symbolic sequence analysis

**Definitions.** Finding some properties of a sequence $\mathbb{x}$ that is constructed over a field $\mathbb{F}$ may be referred to as *sequence analysis* or *discrete-time signal analysis*. If we somehow mathematically alter $\mathbb{x}$ with an operator $\mathbf{A}$ to produce a new sequence $\mathbb{y} \triangleq \mathbf{A}\mathbb{x}$, then this may be referred to as *sequence processing*, or more commonly as *discrete-time signal processing* or *digital signal processing* (*DSP*).

**Basis theory.** Sequence analysis and sequence processing typically make use of basis theory. In basis theory in general (of which Fourier analysis and wavelet analysis are special cases), we represent some point $\mathbb{x}$ ($\mathbb{x}$ is a sequence) in a Banach space (a complete normed linear space) by a linear combination of a basis sequence $(\!(x_n)\!)$ such that

$$\mathbb{x} \stackrel{\star}{=} \sum_{n\in Z} a_n x_n$$

[129] [77], ⟨§3.2⟩

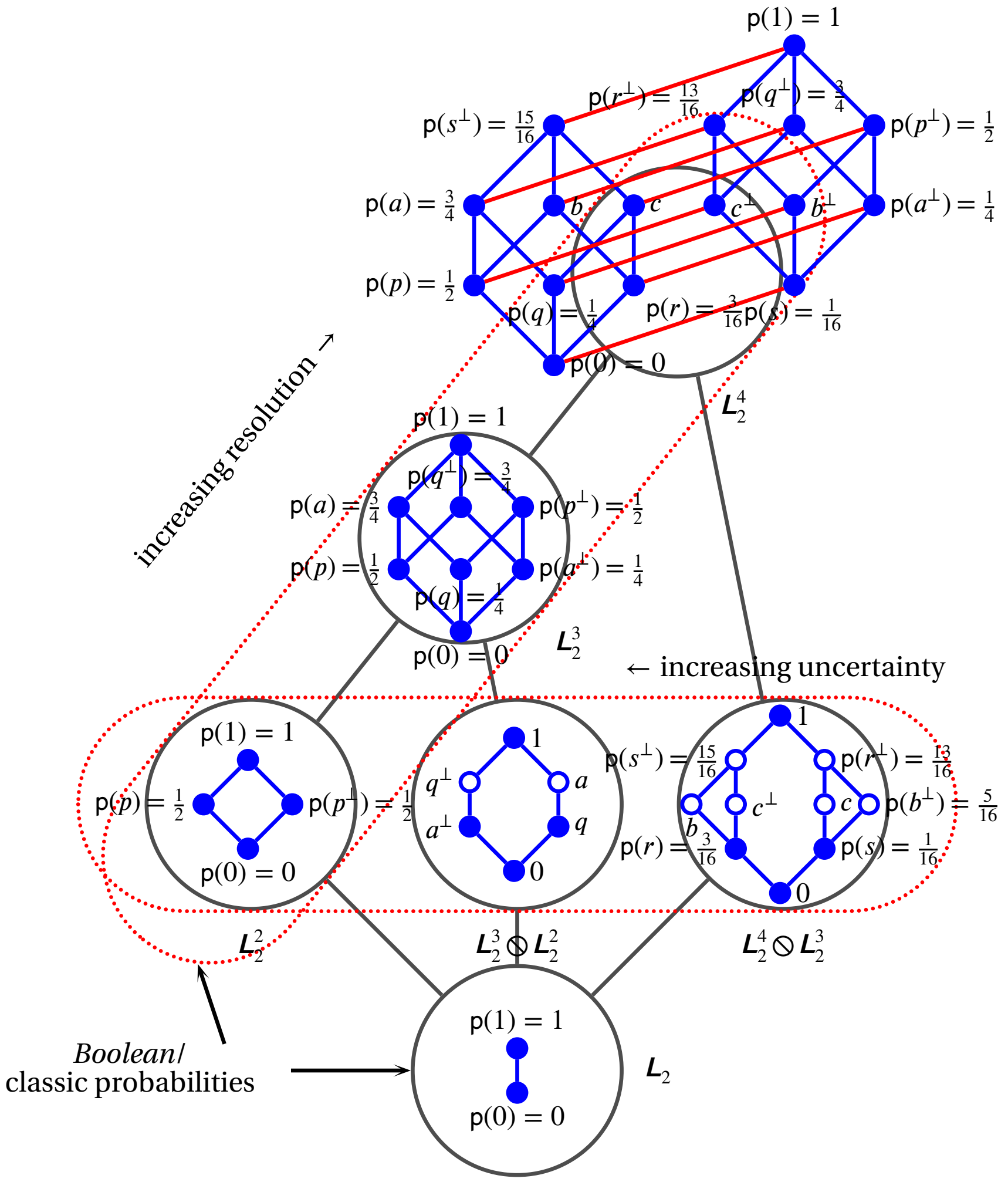

Figure 17: primorial lattice with probability function (Example 5.33 page 68)

where $\stackrel{\star}{=}$ represents strong convergence with respect to the norm $\|\cdot\|$ of the Banach space. Each element $a_n$ is a member of the field $\mathbb{F}$ of the Banach space and the sequence $(\!(a_n)\!)$ is often referred to as a "transform" (Fourier transform, discrete-time Fourier transform, wavelet transform, etc.)

In order to be able to successfully compute any transform (such as a Fourier transform or wavelet transform) in a Banach space or even a finite linear space, the sequence $\mathbb{x}$ needs to be somehow related to the field $\mathbb{F}$ over which the Banach space is constructed.

**The problem.** Let $\tilde{\mathbf{F}}$ be the discrete-time Fourier transform operator and $\mathbf{W}$ be a discrete-time wavelet transform. Suppose we want to compute $\tilde{\mathbf{F}}\mathbb{x}$ or $\mathbf{W}\mathbb{x}$. This is a problem in *symbolic sequence analysis* and *symbolic signal processing* in general because of the following reasons:

1. The symbols in $\mathbb{x}$ have no field structure; so we can't even add them.
2. The symbols in $\mathbb{x}$ have no order structure; so if $A$, $B$, and $C$ are symbols, we can't say, for example, $A < B$ or $B < C$, etc.
3. The symbols in $\mathbb{x}$ have no topology except for some arguably trivial topologies;[130] so we can't say, for example, that $A$ is "closer" to $B$ than it is to $C$, etc.

In fact, *symbol sequence analysis* does not just cause problems for Fourier or wavelet analysis only—it causes problems for basis theory in general because a basis is constructed in a Banach space, and symbolic sequences are in general not constructed in Banach spaces.

A kind of "hack" solution may be to map the symbols to points $(\!(p_1, p_2, \dots, p_N)\!)$ in the *complex plane* $\mathbb{C}$. If these points are chosen such that they are distinct, not on either the real or imaginary axes, and $|p_1| = |p_2| = \dots = |p_N|$, then that would seem to be a good start, because now the mapped symbols have a field structure, and they are arguably unordered (arguably we can't say any one of them is greater or less than any other, just as in the original symbol sequence).

But we still have the topology problem. If we map, say, 4 symbols to 4 points in $\mathbb{C}$ as $p_1 = 1$, $p_2 = -1$, $p_3 = i$, and $p_4 = -i$, then "$p_1$" is closer (with respect to the metric induced by the norm $|\cdot|$) to "$p_3$" then it is to "$p_2$":

$$\mathsf{d}(p_1, p_3) = |p_1 - p_3| = \left(p_1^2 - p_3^2\right)^{1/2} = \left(1^2 - i^2\right)^{1/2} = \sqrt{2} \lneq 2 = \left(2^2 - 0^2\right)^{1/2} = \mathsf{d}(p_1, p_2)$$

This unwanted topological property is introduced by the mapping, will affect the transform, but yet is not a property of the original symbolic sequence.

---

[130] These topologies include the *indiscrete topology* $\{\varnothing, X\}$ where $X \triangleq \{A, B, C\}$, *discrete topology* $2^X$ (references: [126], page 77, [107], page 107, ⟨Example 3.J⟩, [156], pages 42–43, ⟨II.4⟩, [44], page 18 ), and the topology induced by the *discrete metric* $\mathsf{d}(x, y) \triangleq \{1 \text{ for } x \neq y, 0 \text{ for } x = y\}$ (references: [67], page 13, [31], page 24, [101], page 19, ⟨Example 2.1⟩ ).

"Frequency" properties may be useful in *symbolic sequence analysis* and *symbolic sequence processing*. But the point here is that any kind of basis theory technique (including Fourier or wavelet techniques) may result in a kind of imperfect "hack" solution.

**Proposed solution.** The solution proposed here is to perform symbolic sequence analysis using primorial lattices. Suppose we have a sequence $\mathbb{x}$ over a set of $N$ symbols (each element in the sequence can be any one of $N$ different symbols). Let $\mathbb{P}$ be the primorial lattice generated by $\boldsymbol{L}_2^N$. The orthogonal $N$ atoms of $\boldsymbol{L}_2^N$ represent the $N$ symbols. The element $A \vee B$ in $\boldsymbol{L}_2^N$, where $A$ and $B$ are 2 symbols, represents the event of a particular position in the sequence being $A$ OR $B$ (it is not possible for a particular position to be both $A$ AND $B$).

Any symbol in $\boldsymbol{L}_2^N$ can be projected onto any other Boolean or orthocomplemented lattice in $\mathbb{P}$ by use of a *lattice projection* (Section 5.4 page 58). The result of projecting an entire sequence onto a lattice in $\mathbb{P}$ is another sequence (Definition 5.19 page 58). So after projection, a sequence on $\boldsymbol{L}_2^N$ results in $N-1$ sequences of lower resolution and $N-1$ sequences of assorted frequencies. This is similar in form to the *Fast Wavelet Transform*, as illustrated in Figure 10 (page 49).

### 5.6.5 Symbolic sequence processing (SSP)

**Introduction.** The previous section discusses symbolic sequence analysis—meaning we are not trying to change the properties of the sequence, we are only trying to understand its properties. This section discusses *symbolic sequence processing* (or *symbolic signal processing*)—meaning we *are* trying to change the properties of the sequence.

*Digital signal processing* (*DSP*) or *discrete-time signal processing* operates on a sequence constructed over a field $\mathbb{F}$, where $\mathbb{F}$ is typically either $\mathbb{R}$ or $\mathbb{C}$. Often by use of simple multiplication and addition operations on elements of the sequence, one can change the properties of the sequence. Often when the properties are related to Fourier analysis, the DSP operations are called "filtering".

**The problem.** Multiplication and addition operations commonly used in DSP require field properties. In symbolic sequence processing, we don't in general have a field.

**Proposed solution.** Sequence processing of, or "filtering" on, a symbolic sequence $\mathbb{x}$ can be performed by judicious selection and/or rejection of the various projections onto the logics in the primorial lattice $\mathbb{P}$.

For example, if one wants $\mathbb{x}$ at a lower "resolution"s, then simply select the sequence from a projection onto the *Boolean logic* at resolution lower than $N$. If one wants to "filter out" the "high frequency" components of $\mathbb{x}$, then simply discard the projections onto the higher frequency orthocomplemented lattices before synthesizing a new sequence from the "low frequency" component sequences.

Synthesis of two projection sequences $\mathbb{y}$ and $\mathbb{z}$ into a new sequence $\mathbb{x}'$ can be performed, for example, by pointwise join such that

$$\begin{aligned}\mathbb{y} \oplus \mathbb{z} &\triangleq (\!( y_n )\!)_{n\in\mathbb{Z}} \bigvee (\!( z_n )\!)_{n\in\mathbb{Z}} \\ &\triangleq (\!( y_n \vee z_n )\!)_{n\in\mathbb{Z}} \\ &\triangleq (\!( x_n )\!)_{n\in\mathbb{Z}} \\ &\triangleq \mathbb{x}\end{aligned}$$

### 5.6.6 Genomic Signal Processing (GSP)

*Genomic Signal Processing* (*GSP*) is simply a special case of *Symbolic Sequence Processing* with $N = 4$. In GSP, the 4 symbols are commonly referred to as $A$, $C$, $T$, and $G$, each of which corresponds to a nucleobase (adenine, thymine, cytosine, and guanine, respectively).[131] The sequence itself is called a *genome*. A typical genome sequence contains a large number of symbols (about 3 billion for humans, 29751 for the SARS virus).[132]

**Example 5.34** Traditionally in GSP, the symbols $(A \vee T)$ and $(C \vee G)$ are of special interest. Portions of a genome sequence high in $(A \vee T)$ content separate at lower temperatures than do those with high $(C \vee G)$ content.[133] Therefore, one could construct a primorial lattice induced by $\boldsymbol{L}_2^4$ that allows for convenient analysis of $A \vee T$ and/or $C \vee G$ in some lower resolution space. An example is illustrated in Figure 18 (page 73).

**Example 5.35** In some cases, genomic sequences with more than 4 symbols ($N > 4$) have been studied.[134] Figure 19 (page 74) illustrates a primorial lattice with an extra symbol $X$ in the higher resolution $L_2^5$ Boolean lattice, but with only the symbols $A$, $C$, $G$, and $T$ in the lower resolution $L_2^4$ Boolean lattice. The symbol $X$ can be projected onto any of the lower resolution spaces using a *projection operator* (Section 5.4 page 58).

[131] [121], ⟨Mendel (1853): gene coding uses discrete symbols⟩, [165], page 737, ⟨Watson and Crick (1953): gene coding symbols are adenine, thymine, cytosine, and guanine⟩, [164], page 965, [142], page 52

[132] [1], ⟨http://www.ncbi.nlm.nih.gov/genome/guide/human/⟩, ⟨Homo sapiens, NC_000001–NC_000022 (22 chromosome pairs), NC_000023 (X chromosome), NC_000024 (Y chromosome), NC_012920 (mitochondria)⟩, [1], ⟨http://www.ncbi.nlm.nih.gov/nuccore/30271926⟩, ⟨SARS coronavirus, NC_004718.3⟩ [150], ⟨homo sapien chromosome 1⟩, [149], ⟨SARS coronavirus⟩

[133] [32], page 13, ⟨Remark 1.2⟩

[134] [30], [53]

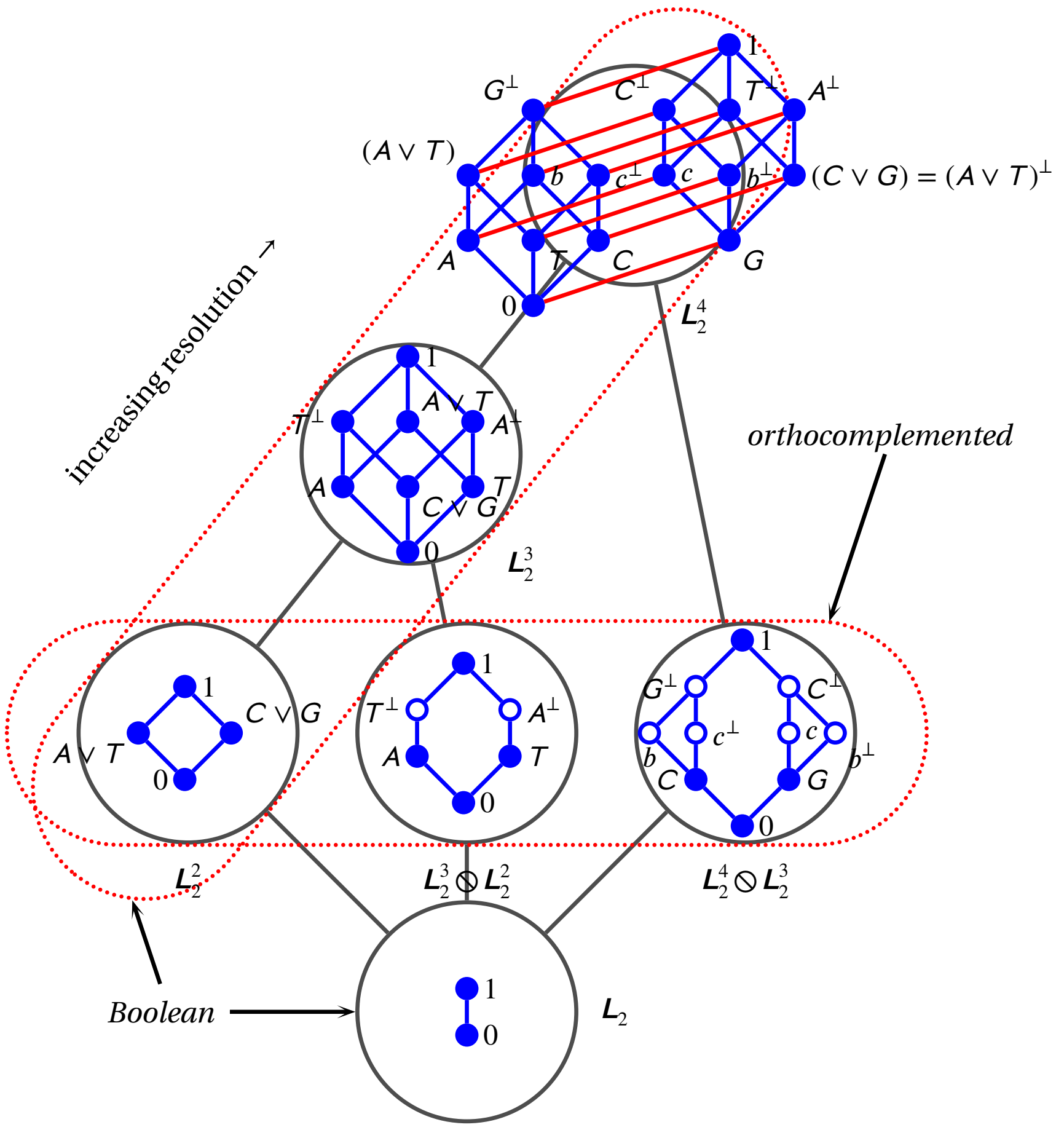

Figure 18: primorial lattice for genomic signal processing (GSP) with $A \vee T$ and $C \vee G$ analysis features (Example 5.34 page 72)

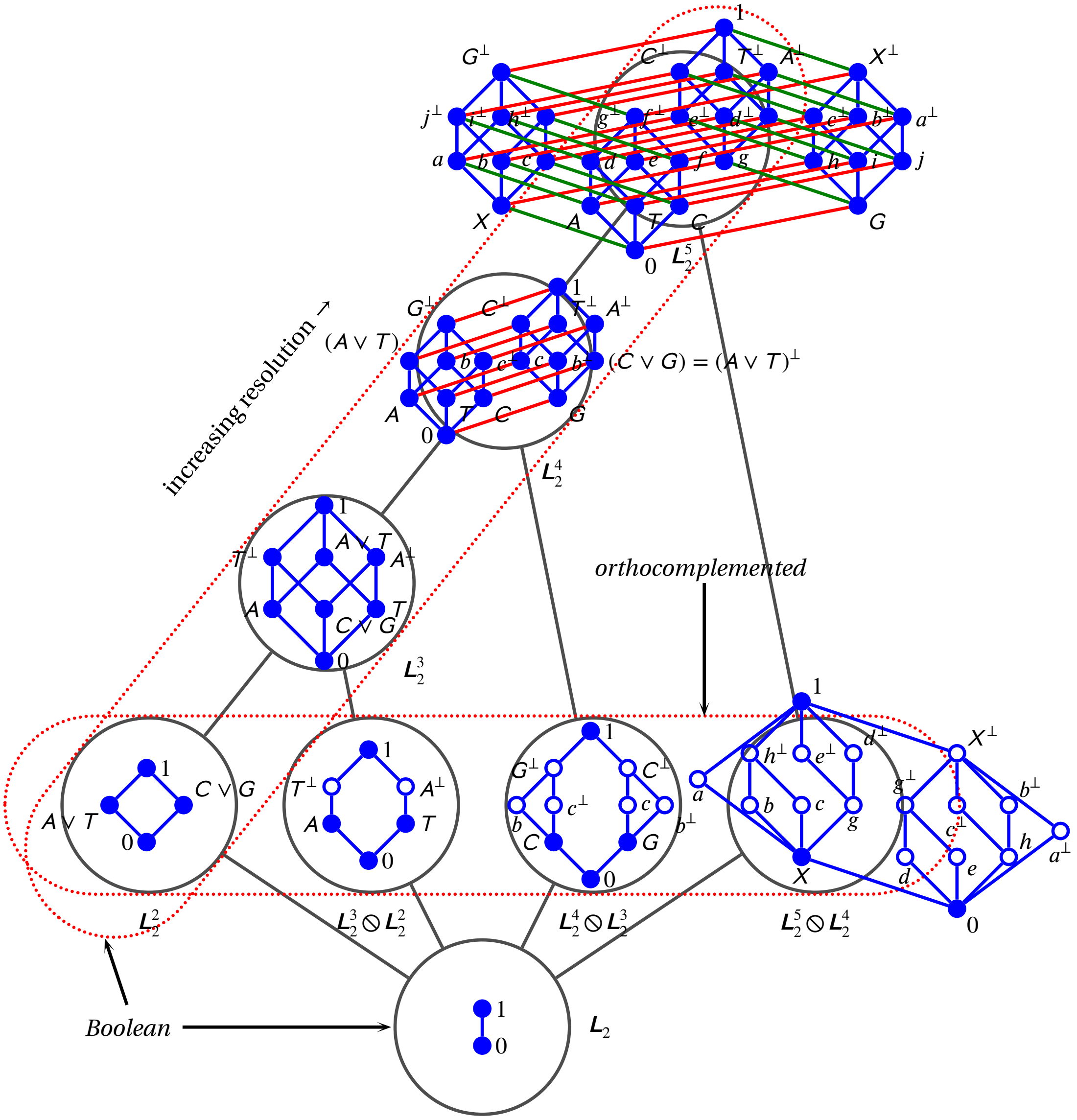

Figure 19: primorial lattice for genomic signal processing (GSP) with extra symbol $X$ (Example 5.35 page 72)

# References

[1] *GenBank* (2014) Available at `http://www.ncbi.nlm.nih.gov`

[2] *On-Line Encyclopedia of Integer Sequences*, World Wide Web (2014) Available at `http://oeis.org/`

[3] **M E Adams**, *Uniquely Complemented Lattices*, from: "The Dilworth theorems: selected papers of Robert P. Dilworth'', (K P Bogart, R S Freese, J P Kung, editors), Birkhäuser, Boston (1990) 79–84

[4] **E H Adelson**, **P J Burt**, *Image Data Compression with the Laplacian Pyramid*, from: "Proceedings of the Pattern Recognition and Information Processing Conference'', IEEE Computer Society Press, Dallas Texas (1981) 218–223

[5] **C D Aliprantis**, **O Burkinshaw**, *Principles of Real Analysis*, 3 edition, Acedemic Press, London (1998)

[6] **L Alvarez**, **F Guichard**, **P-L Lions**, **J M Morel**, *Axioms and fundamental equations of image processing*, Archive for Rational Mechanics and Analysis 123 (1993) 199–257

[7] **K E Atkinson**, **W Han**, *Theoretical Numerical Analysis: A Functional Analysis Framework*, volume 39 of *Texts in Applied Mathematics*, 3 edition, Springer (2009)

[8] **G Bachman**, **L Narici**, **E Beckenstein**, *Fourier and Wavelet Analysis*, Universitext Series, Springer (2000)

[9] **R Balbes**, **P Dwinger**, *Distributive Lattices*, University of Missouri Press, Columbia (1975)2011 reprint edition available (ISBN 9780983801108)

[10] **R Bellman**, **M Giertz**, *On the analytic formalism of the theory of fuzzy sets*, Information Sciences 5 (1973) 149–156

[11] *A Prelude to Sampling, Wavelets, and Tomography*, from: "Sampling, Wavelets, and Tomography'', (J Benedetto, A I Zayed, editors), Applied and Numerical Harmonic Analysis, Springer (2004) 1–32

[12] **L Beran**, *Orthomodular Lattices: Algebraic Approach*, Mathematics and Its Applications (East European Series), D. Reidel Publishing Company, Dordrecht (1985)

[13] **P Billingsley**, *Probability And Measure*, 3 edition, Wiley series in probability and mathematical statistics (1995)

[14] **G Birkhoff**, *On the Combination of Subalgebras*, Mathematical Proceedings of the Cambridge Philosophical Society 29 (1933) 441–464

[15] **G Birkhoff**, *The Logic of Quantum Mechanics*, Annals of Mathematics 37 (1936) 823–843

[16] **G Birkhoff**, *Lattices and their applications*, Bulletin of the American Mathematical Society 44 (1938) 1:793–800

[17] **G Birkhoff**, *Lattice Theory*, 2 edition, American Mathematical Society, New York (1948)

[18] **G Birkhoff**, *Lattice Theory*, volume 25 of *Colloquium Publications*, 3 edition, American Mathematical Society, Providence (1967)

[19] **G Birkhoff**, **P Hall**, *Applications of lattice algebra*, Mathematical Proceedings of the Cambridge Philosophical Society 30 (1934) 115–122

[20] **G Birkhoff**, **J V Neumann**, *The Logic of Quantum Mechanics*, The Annals of Mathematics 37 (1936) 823–843

[21] **D Black**, **I Brooks**, **e a Robert Groves** (editors), *Collins English Dictionary—Complete and Unabridged*, 10 edition (2009)

[22] **T S Blyth**, *Lattices and ordered algebraic structures*, Springer, London (2005)

[23] **S Burris**, **H P Sankappanavar**, *A Course in Universal Algebra*, 1 edition, Graduate texts in mathematics 78, Springer-Verlag, New York (1981)2000 edition available for free online

[24] **P J Burt**, **E H Adelson**, *The Laplacian Pyramid As A Compact Image Code*, IEEE Transactions On Communications COM-3L (1983) 532–540

[25] **P G Casazza**, **M C Lammers**, *Bracket Products for Weyl-Heisenberg Frames*, from: “Gabor Analysis and Algorithms: Theory and Applications'', (H G Feichtinger, T Strohmer, editors), Applied and Numerical Harmonic Analysis, Birkhäuser (1998) 71–98

[26] **G Cattaneo**, **D Ciucci**, *Lattices with Interior and Closure Operators and Abstract Approximation Spaces*, from: “Transactions on Rough Sets X'', (J F Peters, A Skowron, editors), Lecture notes in computer science 5656, Springer (2009) 67–116

[27] **O Christensen**, *An Introduction to Frames and Riesz Bases*, Applied and Numerical Harmonic Analysis, Birkhäuser, Boston/Basel/Berlin (2003)

[28] **D M Cifarelli**, **E Regazzini**, *De Finetti's Contribution to Probability and Statistics*, Statistical Science 11 (1996) 253–282

[29] **D W Cohen**, *An Introduction to Hilbert Space and Quantum Logic*, Problem Books in Mathematics, Springer-Verlag, New York (1989)

[30] **A Comnish-Bowden**, *Nomenclature for incompletely specified bases in nucleic acid sequences: recommendations 1984*, Nucleic Acids Research 13 (1985) 3021–3030

[31] **E T Copson**, *Metric Spaces*, Cambridge tracts in mathematics and mathematical physics 57, Cambridge University Press, London (1968)

[32] **N Cristianini**, **M W Hahn**, *Introduction to Computational Genomics: A Case Studies Approach*, Cambridge University Press (2007)

[33] **X Dai**, **D R Larson**, *Wandering vectors for unitary systems and orthogonal wavelets*, volume 134 of *Memoirs of the American Mathematical Society*, American Mathematical Society, Providence R.I. (1998)

[34] **X Dai**, **S Lu**, *Wavelets in subspaces*, Michigan Math. J. 43 (1996) 81–98

[35] **I Daubechies**, *Ten Lectures on Wavelets*, Society for Industrial and Applied Mathematics, Philadelphia (1992)

[36] **B A Davey**, **H A Priestley**, *Introduction to Lattices and Order*, 2 edition, Cambridge mathematical text books, Cambridge University Press, Cambridge (2002)

[37] **A C Davis**, *A characterization of complete lattices*, Pacific Journal of Mathematics 5 (1955) 311–319

[38] **R Dedekind**, *Ueber die von drei Moduln erzeugte Dualgruppe*, Mathematische Annalen 53 (1900) 371–403Regarding the Dual Group Generated by Three Modules

[39] **D Devidi**, *Negation: Philosophical Aspects*, from: "Encyclopedia of Language & Linguistics'', (K Brown, editor), 2 edition, Elsevier (2006) 567–570

[40] **D Devidi**, *Negation: Philosophical Aspects*, from: "Concise Encyclopedia of Philosophy of Language and Linguistics'', (A Barber, R J Stainton, editors), Elsevier (2010) 510–513

[41] **E Deza**, **M-M Deza**, *Dictionary of Distances*, Elsevier Science, Amsterdam (2006)

[42] **M-M Deza**, **E Deza**, *Encyclopedia of Distances*, Springer (2009)

[43] **M M Deza**, **M Laurent**, *Geometry of Cuts and Metrics*, volume 15 of *Algorithms and Combinatorics*, Springer, Berlin/Heidelberg/New York (1997)

[44] **E DiBenedetto**, *Real Analysis*, Birkhäuser Advanced Texts, Birkhäuser, Boston (2002)

[45] **J A Dieudonné**, *Foundations of Modern Analysis*, Academic Press, New York (1969)

[46] **R Dilworth**, *Lattices With Unique Complements*, Transactions of the American Mathematical Society 57 (1945) 123–154

[47] **R Dilworth**, *A Decomposition Theorem for Partially Ordered Sets*, Annals of Mathematics 51 (1950) 161–166

[48] **R Dilworth**, *A Decomposition Theorem for Partially Ordered Sets*, from: "The Dilworth theorems: selected papers of Robert P. Dilworth'', (K P Bogart, R S Freese, J P Kung, editors), Birkhäuser (1990), Boston (1950) ?

[49] **R Dilworth**, *Aspects of distributivity*, Algebra Universalis 18 (1984) 4–17

[50] **J Doner**, **A Tarski**, *An extended arithmetic of ordinal numbers*, Fundamenta Mathematicae 65 (1969) 95–127

[51] **J M Dunn**, *Generalized Ortho Negation*, from: "Negation: A Notion in Focus'', (H Wansing, editor), Perspektiven der Analytischen Philosophie / Perspectives in Analytical Philosophy 7, De Gruyter (1996) 3–26

[52] **J M Dunn**, *A Comparative Study of Various Model-theoretic Treatments of Negation: A History of Formal Negation*, from: "What is Negation?'', (D M Gabbay, H Wansing, editors), Applied Logic Series 13, De Gruyter (1999) 23–52

[53] **L Elnitski**, **R C Hardison**, **J Li**, **S Yang**, **D Kolbe**, **P Eswara**, **M J O'Connor**, **S Schwartz**, **W Miller**, **F Chiaromonte**, *Distinguishing Regulatory DNA From Neutral Sites*, Genome Research 13 (2003) 64–72

[54] **M Erné**, **J Heitzig**, **J Reinhold**, *On the number of distributive lattices*, The Electronic Journal of Combinatorics 9 (2002)

[55] **Euclid**, *Elements* (circa 300BC)

[56] **J D Farley**, *Chain Decomposition Theorems For Ordered Sets And Other Musings* (1997) 1–12

[57] **W Feller**, *An Introduction to Probability Theory and its Applications Volume I*, 3, revised edition, Wiley series in probability and mathematical statistics, John Wiley & Sons (1970)

[58] **B de Finetti**, *Ancora sull'estensione alle classi numerabili del teorema delle probabilità totali*, Rendiconti del Reale Istituto Lombardo di Scienze e Lettere 63 (1930) 1063–1069

[59] **B de Finetti**, *A proposito dell'estensione del teorema delle probabilita totali alle classi numerabili*, Rendiconti del Reale Istituto Lombardo di Scienze e Lettere 63 (1930) 901–905

[60] **B de Finetti**, *Sui passaggi al limite nel calcolo delle probabilità*, Rendiconti del Reale Istituto Lombardo di Scienze e Lettere 63 (1930) 155–166

[61] **B Forster**, **P Massopust** (editors), *Four Short Courses on Harmonic Analysis: Wavelets, Frames, Time-Frequency Methods, and Applications to Signal and Image Analysis*, Applied and Numerical Harmonic Analysis, Springer (2009)

[62] **D J Foulis**, *A note on orthomodular lattices*, Portugaliae Mathematica 21 (1962) 65–72

[63] **M R Fréchet**, *Sur quelques points du calcul fonctionnel (On some points of functional calculation)*, Rendiconti del Circolo Matematico di Palermo 22 (1906) 1–74Rendiconti del Circolo Matematico di Palermo (Statements of the Mathematical Circle of Palermo)

[64] **M R Fréchet**, *Les Espaces abstraits et leur théorie considérée comme introduction a l'analyse générale*, Borel series, Gauthier-Villars, Paris (1928)Abstract spaces and their theory regarded as an introduction to general analysis

[65] **M Fréchet**, *Sur l'extension du théorème des probabilités totales au cas d'une suite infinie d'evénéments*, Ren- diconti del Reale Istituto Lombardo di Scienze e Lettere 63 (1930) 899–900

[66] **M Fréchet**, *Sur l'extension du théorème des probabilités totales au cas d'une suite infinie d'evénéments*, Ren- diconti del Reale Istituto Lombardo di Scienze e Lettere 63 (1930) 1059–1062

[67] **J R Giles**, *Introduction to the Analysis of Metric Spaces*, Australian Mathematical Society lecture series 3, Cambridge University Press, Cambridge (1987)

[68] **S Givant**, **P Halmos**, *Introduction to Boolean Algebras*, Undergraduate Texts in Mathematics, Springer (2009)

[69] **T N T Goodman**, **S L Lee**, **W S Tang**, *Wavelets in Wandering Subspaces*, Transactions of the A.M.S. 338 (1993) 639–654Transactions of the American Mathematical Society

[70] **S Gottwald**, *Many-Valued Logic And Fuzzy Set Theory*, from: “Mathematics of Fuzzy Sets: Logic, Topology, and Measure Theory'', (U Höhle, S Rodabaugh, editors), The Handbooks of Fuzzy Sets 3, Kluwer Academic Publishers (1999) 5–90

[71] **G A Grätzer**, *Lattice Theory; first concepts and distributive lattices*, A Series of books in mathematics, W. H. Freeman & Company, San Francisco (1971)

[72] **G A Grätzer**, *General Lattice Theory*, 2 edition, Birkhäuser Verlag, Basel (2003)

[73] **G A Grätzer**, *Two Problems That Shaped a Century of Lattice Theory*, Notices of the American Mathematical Society 54 (2007) 696–707

[74] **R Greechie**, *Orthomodular lattices admitting no states*, Journal of Combinatorial Theory, Series A 10 (1971) 119–132

[75] **D J Greenhoe**, *Wavelet Structure and Design*, volume 3 of *Mathematical Structure and Design series*, Abstract Space Publishing (2013)

[76] **D J Greenhoe**, *Analysis Structure and Properties*, volume 2 of *Mathematical Structure and Design series*, Abstract Space Publishing (2014)To be published, Lord willing, in 2014 or 2015

[77] **D J Greenhoe**, *Lattice compatible operators for fuzzy logic*, Journal of Logic and Analysis (2014)In review as of 2014 August 21

[78] **D J Greenhoe**, *Sets Relations and Order Structures*, volume 1 of *Mathematical Structure and Design series*, Abstract Space Publishing (2014)To be published, Lord willing, in 2014 or 2015

[79] **S Gudder**, *Quantum Probability*, Probability and Mathematical Statistics, Academic Press (1988)

[80] **F Guichard**, **J-M Morel**, **R Ryan**, *Contrast invariant image analysis and PDE's* (2012) Available at `http://dev.ipol.im/~morel/JMMBook2012.pdf`

[81] **P R Halmos**, *Naive Set Theory*, The University Series in Undergraduate Mathematics, D. Van Nostrand Company, Inc., Princeton, New Jersey (1960)

[82] **F Hausdorff**, *Set Theory*, 3 edition, Chelsea Publishing Company, New York (1937)1957 translation of the 1937 German *Grundzüge der Mengenlehre*

[83] **C Heil**, *A Basis Theory Primer*, expanded edition edition, Applied and Numerical Harmonic Analysis, Birkhäuser, Boston (2011)

[84] **J Heitzig**, **J Reinhold**, *Counting Finite Lattices*, Journal Algebra Universalis 48 (2002) 43–53

[85] **E Hernández**, **G Weiss**, *A First Course on Wavelets*, CRC Press, New York (1996)

[86] **J R Higgins**, *Sampling Theory in Fourier and Signal Analysis: Foundations*, Oxford Science Publications, Oxford University Press (1996)

[87] **S S Holland, Jr**, *A Radon-Nikodym Theorem in Dimension Lattices*, Transactions of the American Mathematical Society 108 (1963) 66–87

[88] **S S Holland, Jr**, *The Current Interest in Orthomodular Lattices*, from: “Trends in Lattice Theory'', (J C Abbott, editor), Van Nostrand-Reinhold, New York (1970) 41–126From Preface: “The present volume contains written versions of four talks on lattice theory delivered to a symposium on Trends in Lattice Theory held at the United States Naval Academy in May of 1966.”

[89] **E V Huntington**, *Sets of Independent Postulates for the Algebra of Logic*, Transactions of the American Mathematical Society 5 (1904) 288–309

[90] **E V Huntington**, *New Sets of Independent Postulates for the Algebra of Logic, With Special Reference to Whitehead and Russell's Principia Mathematica*, Transactions of the American Mathematical Society 35 (1933) 274–304

[91] **K Husimi**, *Studies on the foundations of quantum mechanics I*, Proceedings of the Physico-Mathematical Society of Japan 19 (1937) 766–789

[92] **T Iijima**, *Basic theory of pattern observation*, Papers of Technical Group on Automata and Automatic Control (1959)See Weickert 1999 for historical information

[93] **V I Istrăţescu**, *Inner Product Structures: Theory and Applications*, Mathematics and Its Applications, D. Reidel Publishing Company (1987)

[94] **L Iturrioz**, *Ordered structures in the description of quantum systems: mathematical progress*, from: "Methods and applications of mathematical logic: proceedings of the VII Latin American Symposium on Mathematical Logic held July 29-August 2, 1985'', volume 69, Sociedade Brasileira de Lógica, Sociedade Brasileira de Matemática, and the Association for Symbolic Logic, AMS Bookstore (1988), Providence Rhode Island (1985) 55–75

[95] **B Jawerth**, **W Sweldens**, *An Overview of Wavelet Based Multiresolutional Analysis*, SIAM Review 36 (1994) 377–412

[96] **S Jenei**, *Structure of Girard Monoids on [0,1]*, from: "Topological and Algebraic Structures in Fuzzy Sets: A Handbook of Recent Developments in the Mathematics of Fuzzy Sets'', (S E Rodabaugh, E P Klement, editors), Trends in Logic 20, Springer (2003) 277–308

[97] **W S Jevons**, *Pure Logic or the Logic of Quality Apart from Quantity; with Remarks on Boole's System and the Relation of Logic and Mathematics*, Edward Stanford, London (1864)

[98] **G Kalmbach**, *Orthomodular Lattices*, Academic Press, London, New York (1983)

[99] **G Kalmbach**, *Measures and Hilbert Lattices*, World Scientific, Singapore (1986)

[100] **D W Kammler**, *A First Course in Fourier Analysis*, 2 edition, Cambridge University Press (2008)

[101] **M A Khamsi**, **W Kirk**, *An Introduction to Metric Spaces and Fixed Point Theory*, John Wiley, New York (2001)

[102] **A N Kolmogorov**, *Foundations of the theory of probability*, 2 edition, Chelsea Publishing Company, New Yourk (1933)1956 2nd edition English translation of A. N. Kolmogorov's 1933 "Grundbegriffe der Wahrscheinlichkeitsrechnung"

[103] **A N Kolmogorov**, *Grundbegriffe der Wahrscheinlichkeitsrechnung*, Springer, Berlin (1933)

[104] **F Kôpka**, **F Chovanec**, *D-Posets*, Mathematica Slovaca 44 (1994) 21–34Communicated by Anatolij Dvurečenskij

[105] **A Korselt**, *Bemerkung zur Algebra der Logik*, Mathematische Annalen 44 (1894) 156–157Referenced by Birkhoff(1948)p.133

[106] **C S Kubrusly**, *The Elements of Operator Theory*, 1 edition, Springer (2001)

[107] **C S Kubrusly**, *The Elements of Operator Theory*, 2 edition, Springer (2011)

[108] **R Larson**, **S Andima**, *The lattice of topologies: a survey*, Rocky Mountain Journal of Mathematics 5 (1975) 177–198

[109] **P G Lemarié** (editor), *Les Ondelettes en 1989*, volume 1438 of *Lecture Notes in Mathematics*, Springer-Verlag, Berlin (1990)

[110] **R Lidl**, **G Pilz**, *Applied Abstract Algebra*, Undergraduate texts in mathematics, Springer, New York (1998)

[111] **T Lindeberg**, *Scale-Space Theory in Computer Vision*, The Springer International Series in Engineering and Computer Science, Springer (1993)

[112] **L H Loomis**, *The Lattice Theoretic Background of the Dimension Theory of Operator Algebras*, volume 18 of *Memoirs of the American Mathematical Society*, American Mathematical Society, Providence RI (1955)

[113] **S MacLane**, **G Birkhoff**, *Algebra*, 3 edition, AMS Chelsea Publishing, Providence (1999)

[114] **M D MacLaren**, *Atomic orthocomplemented lattices*, Pacific Journal of Mathematics 14 (1964) 597–612

[115] **F Maeda**, *Kontinuierliche Geometrien*, volume 95, Springer-Verlag, Berlin (1958)

[116] **F Maeda**, **S Maeda**, *Theory of Symmetric lattices*, volume 173 of *Die Grundlehren der mathematischen Wissenschaften in Einzeldarstellungen*, Springer-Verlag, Berlin/New York (1970)

[117] **S Maeda**, *On Conditions for the Orthomodularity*, Proceedings of the Japan Academy 42 (1966) 247–251

[118] **S G Mallat**, *Multiresolution Approximations and Wavelet Orthonormal Bases of $L^2(\mathbb{R})$*, Transactions of the American Mathematical Society 315 (1989) 69–87

[119] **S G Mallat**, *A Wavelet Tour of Signal Processing*, 2 edition, Elsevier (1999)

[120] **R N McKenzie**, *Equational Bases for Lattice Theories*, Mathematica Scandinavica 27 (1970) 24–38

[121] **G Mendel**, *Experiments In Plant Hybridization* (1853) Available at `http://old.esp.org/foundations/genetics/classical/gm-65.pdf`

[122] **Y Meyer**, *Wavelets and Operators*, volume 37 of *Cambridge Studies in Advanced Mathematics*, Cambridge University Press (1992)

[123] **A N Michel**, **C J Herget**, *Applied Algebra and Functional Analysis*, Dover Publications, Inc. (1993)Original version published by Prentice-Hall in 1981

[124] **K E Müller**, *Abriss der Algebra der Logik (Summary of the Algebra of Logic)*, B. G. Teubner (1909)“bearbeitet im auftrag der Deutschen Mathematiker-Vereinigung” (produced on behalf of the German Mathematical Society). “In drei Teilen” (In three parts)

[125] **M Müller-Olm**, *2. Complete Boolean Lattices*, from: “Modular Compiler Verification: A Refinement-Algebraic Approach Advocating Stepwise Abstraction'', Lecture Notes in Computer Science 1283, Springer (1997) 9–14Chapter 2

[126] **J R Munkres**, *Topology*, 2 edition, Prentice Hall, Upper Saddle River, NJ (2000)

[127] **M Nakamura**, *The permutability in a certain orthocomplemented lattice*, Kodai Math. Sem. Rep. 9 (1957) 158–160

[128] **L Narens**, *Theories of Probability: An examination of logical and qualitative foundations*, volume 2 of *Advanced series on mathematical psychology*, World Scientific (2007)

[129] **L Narens**, *Alternative Probability Theories for Cognitive Psychology*, Topics in Cognitive Science 6 (2014) 114–120

[130] **J von Neumann**, *Continuous Geometry*, Princeton mathematical series, Princeton University Press, Princeton (1960)

[131] **H T Nguyen**, **E A Walker**, *A First Course in Fuzzy Logic*, 3 edition, Chapman & Hall/CRC (2006)

[132] **V Novák**, **I Perfilieva**, **J Močkoř**, *Mathematical Principles of Fuzzy Logic*, The Springer International Series in Engineering and Computer Science, Kluwer Academic Publishers, Boston (1999)

[133] **O Ore**, *On the Foundation of Abstract Algebra. I*, The Annals of Mathematics 36 (1935) 406–437

[134] **O Ore**, *Remarks on structures and group relations*, Vierteljschr. Naturforsch. Ges. Zürich 85 (1940) 1–4

[135] **J Packer**, *Applications of the Work of Stone and von Neumann to Wavelets*, from: "Operator Algebras, Quantization, and Noncommutative Geometry: A Centennial Celebration Honoring John Von Neumann and Marshall H. Stone : AMS Special Session on Operator Algebras, Quantization, and Noncommutative Geometry, a Centennial Celebration Honoring John Von Neumann and Marshall H. Stone, January 15-16, 2003, Baltimore, Maryland'', (R S Doran, R V Kadison, editors), Contemporary mathematics—American Mathematical Society 365, American Mathematical Society, Baltimore, Maryland (2004) 253–280

[136] **R Padmanabhan**, **S Rudeanu**, *Axioms for Lattices and Boolean Algebras*, World Scientific, Hackensack, NJ (2008)

[137] **E Pap**, *Null-Additive Set Functions*, volume 337 of *Mathematics and Its Applications*, Kluwer Academic Publishers (1995)

[138] **A Papoulis**, *Probability, Random Variables, and Stochastic Processes*, 3 edition, McGraw-Hill, New York (1991)

[139] **C Peirce**, *On the Algebra of Logic*, American Journal of Mathematics 3 (1880) 15–57

[140] **W J Perschbacher** (editor), *The New Analytical Greek Lexicon*, Hendrickson Publishers, Peabody, Mass. (1990)

[141] **M A Pinsky**, *Introduction to Fourier Analysis and Wavelets*, Brooks/Cole, Pacific Grove (2002)

[142] **J C Pommerville**, *Fundamentals of Microbiology*, Jones & Bartlett Publishers (2013)

[143] **S Qian**, **D Chen**, *Joint time-frequency analysis: methods and applications*, PTR Prentice Hall (1996)

[144] **E Renedo**, **E Trillas**, **C Alsina**, *On the law $(a \cdot b')' = b + a' \cdot b'$ in de Morgan algebras and orthomodular lattices*, Soft Computing 8 (2003) 71–73

[145] **J Riečan**, *K axiomatike modulárnych svázov*, Acta Fac. Rer. Nat. Univ. Comenian 2 (1957) 257–262

[146] **S Roman**, *Lattices and Ordered Sets*, 1 edition, Springer, New York (2008)

[147] **G-C Rota**, *The Number of Partitions of a Set*, The American Mathematical Monthly 71 (1964) 498–504

[148] **G-C Rota**, *The many lives of lattice theory*, Notices of the American Mathematical Society 44 (1997) 1440–1445

[149] **e a Runtao He**, *Analysis of multimerization of the SARS coronavirus nucleocapsid protein*, Biochemical and Biophysical Research Communications 316 476–483

[150] **e a S G Gregory**, *The DNA sequence and biological annotation of human chromosome 1*, Nature: International Weekly Journal of Science 441 (2006) 315–321

[151] **V N Saliĭ**, *Lattices with Unique Complements*, volume 69 of *Translations of mathematical monographs*, American Mathematical Society, Providence (1988)Translation of *Reshetki s edinstvennymi dopolneniíam̀i*

[152] **U Sasaki**, *Orthocomplemented lattices satisfying the exchange axiom*, Journal of Science of the Hiroshima University 17 (1954) 293–302

[153] **E Schröder**, *Vorlesungen über die Algebra der Logik: Exakte Logik*, volume 1, B. G. Teubner, Leipzig (1890)

[154] **A Shen**, **N K Vereshchagin**, *Basic Set Theory*, volume 17 of *Student mathematical library*, American Mathematical Society, Providence (2002)Translated from Russian

[155] **R P Stanley**, *Enumerative Combinatorics*, volume 49 of *Cambridge studies in advanced mathematics*, 1 edition, Cambridge University Press (1997)

[156] **L A Steen**, **J A Seebach**, *Counterexamples in Topology*, 2, revised edition, Springer-Verlag (1978)A 1995 "unabridged and unaltered republication" Dover edition is available.

[157] **M Stern**, *Semimodular Lattices: Theory and Applications*, volume 73 of *Encyclopedia of Mathematics and its Applications*, Cambridge University Press, Cambridge (1999)

[158] **G Strang**, **T Nguyen**, *Wavelets and Filter Banks*, Wellesley-Cambridge Press, Wellesley, MA (1996)

[159] **L Straßburger**, *What is Logic, and What is a Proof?*, from: "Logica Universalis: Towards a General Theory of Logic'', (J-Y Beziau, editor), Mathematics and Statistics, Birkhäuser (2005) 135–145

[160] **N K Thakare**, **M M Pawar**, **B N Waphare**, *A structure theorem for dismantlable lattices and enumeration*, Journal Periodica Mathematica Hungarica 45 (2002) 147–160

[161] **A Troelstra**, **D van Dalen**, *Constructivism in Mathematics: An Introduction*, volume 121 of *Studies in Logic and the Foundations of Mathematics*, North Holland/Elsevier, Amsterdam/New York/Oxford/Tokyo (1988)

[162] **A de Vries**, *Algebraic hierarchy of logics unifying fuzzy logic and quantum logic* (2007)The registered submission date for this paper is 2007 July 14, but the date appearing on paper proper is 2009 December 6. The latest year in the references is 2006 Available at `http://arxiv.org/abs/0707.2161`

[163] **D F Walnut**, *An Introduction to Wavelet Analysis*, Applied and numerical harmonic analysis, Springer (2002)

[164] **J D Watson**, **F H C Crick**, *Genetical Implications of the Structure of Deoxyribonucleic Acid*, Nature (1953) 964–967

[165] **J D Watson**, **F H C Crick**, *Molecular Structures of Nucleic Acids: A Struture for Deoxyribose Nucleic Acid*, Nature (1953) 737–738

[166] **J Weickert**, *Linear Scale-Space has First been Proposed in Japan*, Journal of Mathematical Imaging and Vision 10 (1999) 237–252

[167] **A N Whitehead**, *A Treatise on Universal Algebra with Applications*, volume 1, University Press, Cambridge (1898)

[168] **P Wojtaszczyk**, *A Mathematical Introduction to Wavelets*, volume 37 of *London Mathematical Society student texts*, Cambridge University Press (1997)

*Telecommunications Engineering Department, National Chiao-Tung University, Hsinchu, Taiwan*

dgreenhoe@gmail.com

[135] This document was typeset using XƎLaTeX, which is part of the TeX family of typesetting engines, which is arguably the greatest development since the Gutenberg Press. Graphics were rendered using the pstricks and related packages and LaTeX graphics support. The main roman, *italic*, and **bold** font typefaces used are all from the Heuristica family of typefaces (based on the Utopia typeface, released by Adobe Systems Incorporated). The math font is XITS from the XITS font project. The font used for the title in the footers is Adventor (similar to Avant-Garde) from the TeX-Gyre Project. The font used for the version number in the footers is LIQUID CRYSTAL (Liquid Crystal) from FontLab Studio. This handwriting font is Lavi from the Free Software Foundation.